\documentclass[11pt]{article}

\usepackage{amsmath}
\usepackage{amssymb}
\usepackage{amsthm}
\usepackage{latexsym}
\usepackage{color}
\usepackage{graphicx}
\usepackage{appendix}
\usepackage{color} 
\usepackage{enumerate}
\usepackage[T1]{fontenc}
\usepackage[english]{babel}
\usepackage[table]{xcolor}
\usepackage{geometry}
\DeclareSymbolFont{calletters}{OMS}{cmsy}{m}{n}
\DeclareSymbolFontAlphabet{\mathcal}{calletters}

\def\be{\begin{eqnarray}}
\def\ee{\end{eqnarray}}

\def\b*{\begin{eqnarray*}}
\def\e*{\end{eqnarray*}}

\newtheorem{Theorem}{Theorem}[part]

\newtheorem{Proposition}{Proposition}[part]

\newtheorem{Assumption}{Assumption}[part]
\newtheorem{Lemma}{Lemma}[part]

\newtheorem{Remark}{Remark}[part]

\makeatletter \@addtoreset{equation}{section}

\@addtoreset{Definition}{section}

\@addtoreset{Theorem}{section}

\@addtoreset{Proposition}{section}

\@addtoreset{Property}{section}

\@addtoreset{Assumption}{section}

\@addtoreset{Corollary}{section}

\@addtoreset{Lemma}{section}

\@addtoreset{Remark}{section}

\@addtoreset{Example}{section}

\@addtoreset{Condition}{section}

\newcommand{\abs}[1]{\left|#1\right|}     

\def \D{{\bf{D}}}
\def \E{\mathbb{E}}
\def \F{\mathbb{F}}

\def \N{\mathbb{N}}
\def \P{\mathbb{P}}
\def \Q{\mathbb{Q}}
\def \R{\mathbb{R}}

\def\Ac{{\cal A}}

\def\Lc{{\cal L}}
\def\Mc{{\cal M}}

\def\Ut{{\widetilde U}}

\def\cbf{{\bf c}}
\def\ybf{{\bf y}}

\def\Tr#1{{\rm Tr}\left[#1\right]}

\def \Frac{\displaystyle\frac}

\def\x{\times}

\def\={\;=\;}
\def\.{\;.}

\def\eps{\varepsilon}

\def\reff#1{{\rm(\ref{#1})}}

\def\1{{\bf 1}}
\def \ep{\hbox{ }\hfill{ ${\cal t}$~\hspace{-5.1mm}~${\cal u}$   } }

\def \proof{{\noindent \bf Proof. }}
\def\ep{\epsilon}
\def\cbf{{\mathbf{c}}}

\def\ep{\epsilon}
\def\eps{\epsilon}

\def\b*{\begin{eqnarray*}}
\def\e*{\end{eqnarray*}}

 \def\normeL2#1{\left\|{#1}\right\|_{L^2}}
 
 \title{General indifference pricing with small transaction costs}

 \author{Dylan {\sc Possama\"{i}} \footnote{CEREMADE, Universit\'e Paris Dauphine, possamai@ceremade.dauphine.fr.}
 \and Guillaume {\sc Royer}\footnote{CMAP, Ecole Polytechnique Paris, guillaume.royer@polytechnique.edu.}}
 
             \date{\today}

\begin{document}

 \maketitle

\begin{abstract}
We study the utility indifference price of a European option in the context of small transaction costs. Considering the general setup allowing consumption and a general utility function at final time $ T$, we obtain an asymptotic expansion of the utility indifference price as a function of the asymptotic expansions of the utility maximization problems with and without the European contingent claim. We use the tools developed in \cite{st} and \cite{pst} based on homogenization and viscosity solutions to characterize these expansions. Finally we study more precisely the example of exponential utilities, in particular recovering under weaker assumptions the results of \cite{bichuch}.

\vspace{5mm}

\vspace{5mm}


\end{abstract}

\section{Introduction}
It is a widely known result in the finance literature that in any complete market, an investor who sold a contingent claim can replicate it perfectly by continuously trading a portfolio consisting of cash and the underlying risky asset. However, the corresponding strategies generally lead to portfolio rebalanced in continuous time, and which therefore are generically of unbounded variation. Thus, as soon as any transaction costs are introduced in the market, such strategies have exploding costs and cannot therefore be of any use. A possible way out of this is for the investor to search for super-replicating portfolios, instead of replicating ones. However, it turns out that in a market with transaction costs, the simple problem of super-replicating a Call option can only be solved by using the trivial buy-and-hold strategy, therefore leading to prohibitive costs. These types of results have been first conjectured by Davis and Clark \cite{dc}, and proved under more and more general frameworks (see among others \cite{bouch2,bt00,cs,ck,cpt,jk,kaba,kl,ks,kpt1,kpt2,ssc,touzi}). A rather natural alternative approach has been first proposed by Hodges and Neuberger \cite{hn}, and basically states that the price of a given contingent claim should be equal to the minimal amount of money that an investor has to be offered so that he becomes indifferent (in terms of utility) between the situation where he has sold the claim and the one where he has not. Such an approach is therefore very naturally linked to the general problem of investment and consumption under transaction costs, which has received a lot of attention since the seminal papers by Magill and Constantinides \cite{mc} and Constantinides \cite{c}. 

\vspace{0.25em}
Following these two works which rather concentrated on the numerical aspects of the problem, but contained already the fundamental insight that the no-transaction region is a wedge, Taksar, Klass and Assaf \cite{tka} studied an ergodic version of the maximization, before the classical paper of Davis and Norman put the problem into the modern framework of singular stochastic control theory. Building upon these works, Soner and Shreve \cite{ss} proposed a comprehensive analysis of the one-dimensional case (that is to say when there is only one risky asset in the market), using the dynamic programming approach as well as the theory of viscosity solutions (see also the earlier work of Dumas and Luciano \cite{dl} in this direction). Their approach was then extended to the case of several risky assets by Akian, Menaldi and Sulem \cite{ams} (see also \cite{ass}). Starting from there, an important strand of literature concerned itself with the problem of option pricing under transaction costs via utility maximization. 


\vspace{0.25em}
The first important result in this direction was obtained by Davis, Panas and Zariphopoulou \cite{dpz}, where they showed that the problem of pricing an European option in a market with proportional transaction costs boiled down to solve two stochastic optimal control problems, whose value functions were shown to be the unique viscosity solutions of quasi-linear variational inequalities. 
Then, starting with the work of Barles and Soner \cite{bs}, where they derived rigorously the limiting behavior of the aforementioned value function as both the transaction costs and the investor risk tolerance go to $0$, many papers studied practically relevant limiting regimes. Indeed, the quasi-variational inequalities derived in \cite{dpz} are difficult to handle numerically, especially in high dimensions, which make asymptotic expansions a lot more tractable\footnote{Let us nonetheless mention that starting with the paper by Dai and Yi \cite{dy}, who showed that the solution to the problem of optimal investment could be written as parabolic double obstacle problem, other numerical approaches are available (see also \cite{djly} and \cite{dz}).}. Thus, Whalley and Wilmott \cite{ww} considered a formal asymptotic expansion for small transaction costs, while Bouchard \cite{bou} and Bouchard, Kabanov and Touzi \cite{bkt} showed rigorously that as the risk-aversion of the investor goes to infinity, the utility indifference price of a contingent claim goes to the corresponding super-replication price. The first rigorous proof of the result in \cite{ww} was obtained, in one dimension, in the appendix of \cite{ss}. Since then, several rigorous results \cite{bichuch2,gms,js1,js2,r} (still in the one dimensional case), as well as formal ones \cite{am,go,km1,km2} were also derived. Let us also point out the work \cite{bichuch} by Bichuch, where the author obtained an asymptotic expansion of the utility indifference price of smooth contingent claims in a market with small transaction costs and on risky asset having geometric Brownian motion dynamics. However, the multidimensional problem, which presents intriguing free boundary problems for which regularity results remain scarce (see for instance \cite{ss2,ss3} for results in a related setting or the recent paper by Chen and Dai \cite{cd} which studies rigorously the shape of the no-transaction region), has remained out of reach until the paper by Bichuch and Shreve \cite{bichs}, where they treated the case of two risky assets which diffuse as arithmetic Brownian motions. Nonetheless, their method of proof requires not only to construct sharp sub and super-solutions to the dynamic programing equation satisfied by the value function, they also need to do a lengthy coefficient matching for the formal expansion, which has to be done on 9 regions of the space, which would become $3^n$ regions for $n$ risky assets, and is therefore ill-suited for arbitrary dimensions. The breakthrough for the treatment of the general problem was achieved by Soner and Touzi \cite{st}, where they connected the asymptotic expansion for small transaction costs to the theory of homogenization. Hence, the first order term in their expansion is shown to be written in terms of the so-called eigenvalue associated to the dynamic programming equation of an ergodic stochastic control problem. This identification allows them to construct a rigorous proof similar to the ones in homogenization theory, even though the problem at hand is not the typical incarnation of such a type of problem, notably because the "oscillatory" (or "fast") variable only appears after a change of variables and is not directly modeled in the original equations. Their approach, limited to dimension one in \cite{st}, was robust enough to cover general Markovian dynamics for the risky asset, as well as general utility functions (whereas the previous literature was limited to power utility and geometric Brownian motion dynamics), was then extended to an arbitrary number of risky assets in the model of Kabanov \cite{kaba} by Possama\"{i}, Soner and Touzi \cite{pst}. Since then, their method received a lot of attention, enabling notably Altarovici, Muhle-Karbe and Soner \cite{altams} to treat the fixed cost problem, Bouchard, Moreau and Soner \cite{bms} an hedging problem under expected loss constraints and Moreau, Muhle-Karbe and Soner \cite{mms} a price impact model. Finally, we would also like to mention the very recent paper by Kallsen and Li \cite{kalli}, which uses convex duality technics to prove rigorously the expansion for small costs in possibly non-Markovian, one-dimensional models.

\vspace{0.25em}
Our paper remains in the context of the general approach initiated by \cite{st}, and our main goal is to provide rigorous asymptotic expansions of the utility indifference price of European contingent claims in general Markovian, multidimensional models and with general utility functions. To the best of our knowledge, the only related papers in the literature are \cite{bichuch} and the very recent manuscript \cite{bms}. However, the level of generality we consider is new, in particular since both these works are restricted to the one dimension case. Furthermore, \cite{bichuch} is restricted to exponential utilities, because their scaling properties allow to deduce directly and completely explicitly the price from the value function of the control problem. Hence, it suffices to obtain the expansion for the value function to obtain the expansion for the price, whereas in our case, even though we follow the same approach, the expansion for the price cannot be deduced so easily. Moreover, our method of proof allows to weaken strongly the assumptions made in \cite{bichuch}, since, for instance, we roughly only need to assume $C^1$-regularity of the option payoffs we consider, while \cite{bichuch} needed $C^4$ regularity. When compared with \cite{bms}, even though we think that their approach could reasonably be extended to the same multi-dimensional setting as ours, the methods with which they approach the problem is different from ours, since they attack directly the expansion for the price, while we first start with an expansion for the value function. Besides, the set of assumptions under which their result for the utility indifference price holds true also implies strong regularity for the payoff functions, which makes our result more general in this regard.

\vspace{0.25em}
The rest of the paper is organized as follows. In Section \ref{sec.gen}, we present succinctly the markets we consider, with and without frictions, and we follow the general approach of \cite{st} to give formal asymptotics for both the value function and the utility indifference price. Section \ref{sec.main} is then devoted to the main results of the paper, as well as the general assumptions under which we will be working and the proof of the expansion for the price. Then, in Section \ref{sec:disc}, we discuss the particular example of exponential utility and compare our result with the existing literature. Finally, Section \ref{sect: proof of th convergence} provides the proofs of all the technical results of the paper.

\vspace{0.25em}
{\textbf{Notations:}} Throughout the paper, we will denote, for any $n\in \N^*$ and any $(x,y)\in\R^n\times\R^n$, by $x\cdot y$ the usual inner product on $\R^n$. For any finite-dimensional normed vector space $E$, and for any $(n,p)\in\N^*\times\N^*$, $\mathcal M_{n,p}(E)$ will denote the space of $n$ by $p$ matrices whose elements take value in $E$. When $n=p$, we simplify the notation $\Mc_n(E):=\Mc_{n,n}(E)$. As usual, we identify $\Mc_{n,1}(E)$ and $E^n$. By a slight abuse of notation, the norm associated to elements in $\R^n$ or $\Mc_n(E)$ will be denoted by $\abs{\cdot}$. The transpose of a matrix $M\in\Mc_{n,p}(E)$ will be denoted by $M^T$. We will also denote by ${\bf 1}_n$ a vector in $\R^n$ whose entries are all equal to $1$.

\section{General setting}\label{sec.gen}
In this section we describe the problem and recall how to obtain formal asymptotics.

\subsection{Financial market with frictions}
We work on a given probability space $(\Omega,\mathcal F,\P)$ on which is defined a $d$-dimensional Brownian motion $W$. For a fixed time horizon $T>0$, and for any $t\in[0,T]$, we define the filtration $\mathbb F^t:=(\mathcal F^t_s)_{t\leq s \leq T}$ to be the completed natural filtration of the process $W^t$, defined by $W^t_s:=W_s-W_t,\ s\in[t,T].$ For notational simplicity, we let $\F:=\F^0$. The financial market consists of a non-risky asset $S^0$ and $d$ risky assets with price process $\{S_t=(S^1_t, \ldots,S^d_t)^T,t\in [0,T]\}$ 
given by the stochastic differential equations (SDEs),
 $$
 \frac{dS^0_t}{S^0_t}
 =
 r(t,S_t)dt,~~
 \frac{dS^i_t}{S^i_t}
 =
 \mu^i(t,S_t) dt + \sum_{j=1}^d \sigma^{i,j}(t,S_t) dW^j_t,
 ~~1\le i\le d,
 $$
where $r: [0,T]\times\R^d\longrightarrow \R_+$ is the instantaneous interest rate and $\mu: [0,T]\times\R^d \longrightarrow \R^d$, $\sigma: [0,T] \times\R^d\longrightarrow\Mc_{d}(\R)$ are the coefficients of instantaneous mean return and volatility, satisfying the standing assumptions that $r, \mu, \sigma \mbox{ are bounded and Lipschitz, and }
 (\sigma \sigma^{T})^{-1}~~\mbox{is bounded.}$
In particular, this guarantees the existence and
the  uniqueness of a strong solution to SDEs.

\vspace{0.25em}
The portfolio of an investor is represented by the dollar value $X$ invested in the non-risky asset, the vector process $Y=(Y^1,\ldots,Y^d)^T$ of 
the value of the positions in each risky asset, and a  short position in a European option represented by some payoff function $ g:\R^d \longrightarrow \R_+$, that he has to hold until the final time $T$. Starting from any time $t\in[0,T]$, these state variables are controlled by the choices of the total amount of transfers $L^{i,j}_s$, $0\le i,j\le d$, from the $i$-th to the $j$-th asset cumulated between time $t$ and $s$. Naturally, the control processes $\{L^{i,j}_s,s\ge t\}$ are defined
as RCLL, nondecreasing, $\F^t$-progressively measurable processes with $L^{i,j}_{t^-}=0$ and $L^{i,i}\equiv 0$. 

\vspace{0.25em}
In addition to the trading activity, the investor consumes between time $t$ and $T$ at a rate determined by a nonnegative $\F^t$-progressively measurable process $\{c_s,t\leq s\le T\}$. Here $c_s$ represents the rate of consumption in terms of the non-risky asset $S^0$, which means that the investor can only consume from the bank account. Such a pair $\nu:=(c,L)$ is called a {\em consumption-investment strategy}. For any $t\in[0,T]$ and any initial position $(X_{t^-},Y_{t^-})=(x,y)\in\R\x\R^{d}$, the portfolio positions of the investor are given by the following state equation
 $$
 dX_u
 =
 \big(r(u,S_u)X_u-c_u\big)du
 +\mathbf{R}^0(dL_u),
 ~~\mbox{and}~~
 dY^i_u
 =
  Y^i_u\;\frac{dS^i_u}{S^i_u}
 +\mathbf{R}^i(dL_u),
 ~~i=1,\ldots,d,
 $$
where 
$$
\mathbf{R}^i(\ell) 
:= 
\sum_{j=0}^{d}
 \big(\ell^{j,i}-(1+\eps^3\lambda^{i,j})\ell^{i,j}\big),
 ~~i=0,\ldots,d,
 \text{ for all } 
 \ell\in\Mc_{d+1}(\R_+),
$$
is the change of the investor's position in the $i-$th asset induced by a transfer policy $\ell$, given a structure of proportional transaction costs $\eps^3\lambda^{i,j}$ for any transfer from asset $i$ to asset $j$. Here, $\eps>0$ is a small parameter, $\lambda^{i,j}\ge 0$, $\lambda^{i,i}=0$, for all $i,j=0,\ldots,d$, and the scaling $\ep^3$ is chosen to state the expansion results simpler. In some instances, we may forbid transactions between certain assets by setting the corresponding transaction costs to $+\infty$, however we will always allow transactions from and to the bank account, that is to say that we always assume
$$\lambda^{i,0}+\lambda^{0,i}<+\infty,\ i=1,\ldots,d.$$
For simplicity, we will also denote $\mathcal I:=\{(i,j)\in\left\{0,\ldots,d\right\}^2,\lambda^{i,j}<+\infty\}.$ For given initial positions $S_{t}=s \in \R_+^d$, $X_{t^-}=x \in \R$, $Y_{t^-}=y \in \R^d$ and a given {\it consumption-investment} strategy $\nu$, we denote by $S^{0,t,s}$, $S^{t,s}$, $X^{t,s,x,\nu}$ and $Y^{t,s,y,\nu}$ the corresponding prices and state processes. Following \cite{bkt}, the strategy $\nu$ will be said to be {\em admissible} for the initial position $(t,s,x,y)$, if the induced state process is well defined and satisfies the following solvency condition:
$$\text{ $\exists\delta^\nu>0$ s.t. $((X^{t,s,x,\nu}_r,(Y^{t,s,y,\nu}_r)^T)^T+\delta^\nu (S^{0,t,s}_r,(S^{t,s}_r)^T)^T)\in K_\eps,$ for all $r\in [t,T]$, $\P-$a.s.,}$$ where the solvency region $K_\eps$ is defined by:
 \begin{align*}
 K_\eps
 :=
 \left\{(x,y)\in\R \times \R^{d}: \ (x,y^T)^T+\mathbf{R}(\ell) \in \R_+^{1+d} \text{ for some }\ell\in\Mc_{d+1}(\R_+) \right\}.
 \end{align*}
The corresponding set of admissible strategies is denoted by $\Theta^\ep(t,s,x,y)$.  The consumption-investment problem is then the following maximization problem,
 \begin{align}\label{def:ve}
 v^{\epsilon,g}(t,s,x,y)
 :=
 \sup_{\nu\in \Theta^\ep(t,s,x,y)}\
 \E_t\left[\int_t^T \kappa e^{-\int_t^\xi k(\iota,S_\iota^{t,s}) d\iota}\ U_1(c_\xi)d\xi+e^{-\int_t^T k(\xi,S_\xi^{t,s})d\xi}\ \mathcal U_2^{\eps,g} \right],
 \end{align}
where $\kappa\in\{0,1\}$ is here so that we can consider simultaneously the problems with or without consumption and where $k:[0,T]\times\R^d$, $U_1:\R\longmapsto \R$ and $\mathcal U_2$ is defined by
$$\mathcal U_2^{\eps,g}:=U_2\left(\ell^\eps\left(X^{t,s,x,\nu}_T,Y^{t,s,y,\nu}_T\right)-g(S_T^{t,s})\right),$$
for some function $U_2:\R\longmapsto\R$ and the liquidation function $\ell^\eps:\R\times\R^d\longmapsto\R$ defined by
$$\ell^\eps(x,y):=\sup\left\{p\in\R,\ \left((x,y^T)^T-p(1,{\bf0}_d^T)^T\right)\in K_\eps\right\}.$$

Such a liquidation function was considered (among others) in \cite{bou,bou2,bkt} and represents the maximum amount, in terms of cash, that the investor can achieve by liquidating its positions in all the risky assets. It is also proved in \cite{bou,bou2,bkt}, that the function $\ell^\eps$ can be rewritten as
\begin{equation}\label{liquidation}
\ell^\eps(x,y)=\underset{r\in K_\eps^*}{\inf}(x,y^T)^T\cdot r,
\end{equation}
where 
$$K_\eps^*:=\left\{r\in\R_+^{d+1},\ r^0=1\text{ and }r^j-(1+\eps^3\lambda^{i,j})r^i\leq 0,\ 0\leq i,j\leq d\right\}.$$

Moreover we assume that $U_1 $ and $ U_2$ are utility functions which are $C^2$, increasing and strictly concave. We also denote the convex conjugate of $ U_1$ by,
 \b*
 \Ut_1(\tilde c)
 &:=& 
 \sup_{c>0} \big\{U_1(c)-c\tilde c\big\},
 \qquad
 \tilde c\in\R,
 \e*
 and by $\text{Supp}(\tilde U_1)$ its support, that is to say the points $\tilde c\in\R$ such that $\tilde U_1(\tilde c)<+\infty.$ Finally, for future reference, we have the following expansion for the function $\ell^\eps$
\begin{Lemma}
We have, uniformly on compact sets,  $0\leq \underset{\eps\downarrow 0}{\overline{\lim}}\ \frac{x+y\cdot{\bf 1}_d-\ell^\eps(x,y)}{\eps^3}<+\infty.$
\end{Lemma}

\textbf{Proof}
It is easy to see using \reff{liquidation} that for any $(x,y)\in\R^{d+1}$
$$0\leq x+y\cdot{\bf 1}_d-\ell^\eps(x,y)=\underset{r\in K_\eps^*}{\sup}\sum_{i=1}^dy^i(1-r^i)\leq \underset{r\in K_\eps^{*,0}}{\sup}\sum_{i=1}^dy^i(1-r^i),$$
where $K_\eps^{*,0}$ contains $K^*_\eps$ and is defined by
$$K_\eps^{*,0}:=\left\{r\in\R_+^{d+1},\ r^0=1\text{ and }\frac{1}{1+\eps^3\lambda^{j,0}}\leq r^j\leq (1+\eps^3\lambda^{0,j}),\ 0\leq j\leq d\right\}.$$
Then, the following immediate consequence finishes the proof
$$\underset{r\in K_\eps^{*,0}}{\sup}\sum_{i=1}^dy^i(1-r^i)=\eps^3\sum_{i=1}^d\left((y^i)^+\frac{\lambda^{i,0}}{1+\eps^3\lambda^{i,0}}+(y^i)^-\lambda^{0,i}\right).$$

\vspace{-2em}

\begin{flushright}$\Box$\end{flushright}
 
\subsection{The Merton problem without frictions}
 The Merton value function $v^g:=v^{0,g}$ corresponds to the limiting case $ \epsilon=0$, where there is no transaction costs. In this case, there is no longer any need to keep track of the transfers between the different assets, and we can take as a state variable the total wealth obtained by aggregating the positions on all the assets. We therefore define $Z:=X+Y\cdot{\bf1}_d$. The dynamics of $Z$ is
 \begin{equation}\label{eq:dynZ}
dZ_t=\left(r(t,S_t)Z_t-c_t\right)dt+\sum_{i=1}^dY_t^i\left(\frac{dS_t^i}{S_t^i}-r(t,S_t)dt\right).
\end{equation}

For given initial positions $S_{t}=s \in \R_+^d$, $Z_{t}=z \in \R_+$ and a given {\it consumption-investment} strategy $\mathfrak{v}:=(c,Y)\in\R_+\times\R^d$, we denote by $Z^{t,s,z,\mathfrak v}$ the corresponding wealth process. In this context, for any $(t,s,z)\in[0,T]\times\R^d\times \R_+$, the set of admissible investment-consumption strategies starting from time $t$ consists of the controls $\mathfrak v$ s.t. $Z^{t,s,z,\mathfrak v}$ is well-defined and there exists $\delta^\mathfrak v>0$ satisfying 
$$Z_r^{t,s,z,\mathfrak v}\geq -\delta^{\mathfrak v}(S^{0,t,s}_r,(S^{t,s}_r)^T)^T\cdot {\bf 1}_{d+1},\ r\in[t,T],\ \P-a.s.$$
We denote this set by $ \Theta^{0}(t,s,z)$. The value function of the Merton problem is then 
\begin{align}\label{def:v0}
 v^{g}(t,s,z)
 :=
 \sup_{\mathfrak{v} \in \Theta^0(t,s,z)}\
 \E_t\left[\int_t^T \kappa e^{-\int_t^\xi k(\iota,S_\iota) d\iota} U_1(c_\xi)d\xi+e^{-\int_t^T k(\xi,S_\xi)d\xi}\ \mathcal U_2^{0,g} \right],
 \end{align}
 where we have defined $ \mathcal U_2^{0,g}:=U_2(Z_T^{t,s,z,\mathfrak{v}}-g(S_T^{t,s})).$ We assume in the following that $v^g$ is smooth, and the unique classical solution of the HJB equation
\begin{align}
\label{def_v0}
\nonumber&kv^g- rz v^g_z - \Lc^0 v^g - \kappa\Ut_1(v^g_z)
 -\underset{y\in\R^d}{\sup}\left\{y\cdot\left((\mu-r\1_d) v^g_z
                                             +\sigma \sigma^{\rm T}\D_{sz} v^g
                                       \right)
                              +\frac12|\sigma^{\rm T}y|^2 v^g_{zz}
                        \right\}
 =
 0 \\
 &v^g(T,s,z)=U_2(z-g(s)),
\end{align}

 where $\D_{sz}:=\frac{\partial}{\partial z}\D_s$, and
 \begin{equation}\label{Lc0}
 \Lc^0
 := \frac{\partial}{\partial t}+
 \mu\cdot\D_s
 +\frac12\mbox{Tr}\big[\sigma \sigma^{\rm{T}} \D_{ss}\big],
 \end{equation}
 with for $i,j=1,\ldots,d$,
 \begin{equation*}
\D_s^i:=s^i\Frac{\partial}{\partial s^i},\ 
\D_{ss}^{i,j}:=s^is^j\Frac{\partial^2}{\partial s^i\partial s^j}, \  \D_s=(\D_{s}^i)_{1\le i\le d}, \text{ and }\D_{ss}:=(\D_{ss}^{i,j})_{1\le i,j\le d}.
 \end{equation*}


Moreover, for a smooth scalar functions $(t,s,x,y)\in [0,T]\times\R^d_+\times\R\times\R^d\longmapsto \psi(t,s,x,y)$ and $(t,s,z)\in [0,T]\times\R^d_+\times\R_+\longmapsto\varphi(t,s,z)$ we set
 \b*
  \psi_{x}:=\frac{\partial\psi}{\partial x}\;\in\R,
\
 \psi_{y}:=\frac{\partial\psi}{\partial y}\;\in\R^d, \ 
  \varphi_{z}:=\frac{\partial\varphi}{\partial z}\;\in\R.
 \e*

The optimal consumption and positioning in the various assets are defined by the functions $\cbf^g(t,s,z)$ and $\ybf^g(t,s,z)$ 
defined, for any $s\in\R^d_+$ and any $z\in\R_+$ by
 \begin{align}\label{def of c}
 \cbf^g(t,s,z) &:= -\kappa\Ut_1^{\prime}\left(v^g_z(t,s,z)\right)
 =\kappa\left(U_1^{\prime}\right)^{-1}\left(v^g_z(t,s,z)\right)\\
\label{def of y opt}
 -v^g_{zz}(t,s,z)\sigma\sigma^{\rm T}(t,s)\ybf^g(t,s,z)
& =
 (\mu-r\1_d)(t,s) v_z^g(t,s,z)
 +\sigma \sigma^{\rm T}(t,s)\D_{sz} v^g(t,s,z).
 \end{align} 

\subsection{The utility indifference price}
We are interested in the so-called utility indifference price of the European option $g$, in both models with or without frictions. They are defined respectively by:
\begin{align}\label{indifeps}
p^{\epsilon,g}(t,s,x) &:= \inf \big\{ p\in\R : v^{\epsilon,g}(t,s,x+p,0) \geq v^{\epsilon,0}(t,s,x,0) \big\}\\
 p^{g}(t,s,x)& := \inf \big\{ p\in\R : v^{g}(t,s,x+p) \geq v^{0}(t,s,x) \big\},
 \label{indif}
 \end{align}
where $v^{\epsilon,0}$ and $v^0$ correspond respectively to the value functions of the problems \eqref{def:ve} and \reff{def_v0} without the option, that is to say when $g=0$. Notice also that we consider here that the initial endowments of the investor are in cash only. This is purely for simplicity and all our results could be easily generalized if we allow the investor to have a non-zero position on the risky assets for the problem with frictions.

\subsection{Dynamic programming}
The dynamic programming equation corresponding to the singular stochastic control problem $v^{\eps,g}$ involves the following differential operators. Let:
\be\label{Lc}
 \Lc
 &:=& \frac{\partial}{\partial t}
+ \mu\cdot\left(\D_s +\D_y\right)
 +r\D_x
 +\frac12\mbox{Tr}\left[\sigma \sigma^{\rm T}\left(\D_{yy}+\D_{ss}+2 \D_{sy}\right) \right],
 \ee
and for $i,j=1,\ldots,d$,
 \begin{align*}
 \D_x&:=x\Frac{\partial}{\partial x},
  \ \D_y^i:=y^i\Frac{\partial}{\partial y^i},\ \D_{yy}^{i,j}:=y^iy^j\Frac{\partial^2}{\partial y^i\partial y^j},\ \D_{sy}^{i,j}:=s^iy^j\Frac{\partial^2}{\partial s^i\partial y^j}\\
 & \D_y=(\D_y^i)_{1\le i\le d},\ \D_{yy}:=(\D_{yy}^{i,j})_{1\le i,j\le d},\ \D_{sy}:=(\D_{sy}^{i,j})_{1\le i,j\le d}.
 \end{align*}

\begin{Theorem}\label{th: viscosity eq veps}
Assume that $v^{\epsilon,g} $ is locally bounded, then it is a viscosity solution of
\begin{equation} \label{e.dpp}
\begin{cases}
\underset{( i,j)\in \mathcal I}{\min}
 \left\{ kv^{\eps,g}- \Lc v^{\eps,g} - \kappa\Ut_1(v^{\eps,g}_{x}) , 
        \ \Lambda_{i,j}^\eps\cdot(v^{\eps,g}_x,(v^{\eps,g}_y)^T)^T 
 \right\} 
 =
 0,  \ (t,s,x,y)\in[0,T) \times \R_+^d \times K_\ep\\[0.5em]
 v^{\eps,g}(T,s,x,y)=U_2\left(\ell^\eps(x,y)-g(s)\right), \ (s,x,y)\in\R_+^d\times K_\eps,
 \end{cases}
 \end{equation}
where $\ \Lambda_{i,j}^\eps
 :=
 e_i-e_j+\eps^3\lambda^{i,j}\;e_i, \ 0\le i,j\le d
.$
Moreover $ v^{\eps,g}$ converges to the Merton value function $ v^g$, as $ \epsilon$ tends to zero.
\end{Theorem}

Let us point out that the result as stated above does not seem to be present in the literature (at least as far as we know) on the subject. Several related results, can be found however, for instance with infinite time-horizon and without consumption (see Kabanov and Safarian \cite{ks2009}), or when consumption and transfers between the risky assets are not allowed (see Akian, Menaldi and Sulem \cite{ams} or Akian, S\'equier and Sulem \cite{ass}). Nonetheless, this is a classical result and does not lie at the heart of our analysis. We will therefore refrain from writing its proof.
\subsection{Formal Asymptotics for the value function}
\label{sect: formal asymptotics for the value function}

Based on \cite{st} and \cite{pst}, we postulate the following expansion for $(t,s,x,y)\in[0,T)\times\R_+^d\times K_\eps$
\begin{align}\label{expansion ve}
v^{\eps,g}(t,s,x,y)=v^g(t,s,x,y)-\epsilon^2 u^g(t,s,z)-\epsilon^4 w^g(t,s,z,\xi)+\circ(\epsilon^2),
\end{align}

where we recall that $z=x+y\cdot{\bf 1}_d$ and we define the "fast" variable $\xi\in\R^d$ by
$$ \xi^i:=\xi^i_\epsilon(t,s,x,y)=\frac{y^i-{\bf{y}}^{g,i}(t,s,z)}{\epsilon},\ 1\leq i\leq d,$$
with the additional useful convention $ \xi^0=0$. We now derive the key equations verified by $ u^g$ and $ w^g$, from the dynamic programming equation \eqref{e.dpp}. The easiest part corresponds to the gradient constraint in \eqref{e.dpp}. By straightforward formal calculations, we have for $ (i,j) \in\mathcal I$
\begin{align*}
\Lambda_{i,j}^\eps \cdot(v^{\eps,g}_x,(v^{\eps,g}_y)^T)^T &= \eps^3 \left( \lambda^{i,j} v^g_z+(e_i-e_j)\cdot D_\xi w^g \right) + \circ(\eps^3)= \eps^3 \left( \lambda^{i,j} v^g_z+w^g_{\xi_i}-w^g_{\xi_j} \right) + \circ(\eps^3).
\end{align*}

We now explore the drift condition in \eqref{e.dpp}. Thank to the linearity of $ \mathcal{L}$, we decompose the calculation in several parts. First of all we have using \eqref{def:v0} and \eqref{def of y opt}
\begin{align}\label{eq:vvg}
\nonumber k v^g-\mathcal{L}  v^g-\kappa\tilde{U}_1( v^g_x ) = &\ k v^g- v_t^g-\mu \cdot \D_s v^g-r z v^g_z + y \cdot \left( r \1_d-\mu \right) v^g_z\\
\nonumber &-\frac{1}{2} \text{Tr} \left[\sigma \sigma^T\left( \D_{yy} v^g+\D_{ss} v^g+2\D_{sy}v^g \right) \right]-\kappa\tilde{U}_1\left( v^g_x \right)\\
\nonumber = &\ (\ybf^g-y) \cdot \left( \left( \mu- r \1_d \right) v^g_z +\sigma \sigma^T \D_{sz} v^g \right) +\frac{1}{2} |\sigma^T \ybf^g |^2 v^g_{zz}-\frac{1}{2} | \sigma^T y|^2 v^g_{zz}\\
 = &\ - \frac{1}{2} \left| \sigma^T (\ybf^g-y) \right|^2 v^g_{zz} =- \frac{\eps^2}{2} \left| \sigma^T \xi\right|^2 v^g_{zz}.
\end{align}
Similarly, we obtain by straightforward but tedious calculations that
\begin{align*}
\eps^4 \left( k w^g-\mathcal{L} w^g\right) &=\frac{\eps^4}{2} \Tr{\D_{yy}w^g+\D_{ss}w^g+2\D_{sy} w^g}+\circ(\eps^2)= \frac{\eps^2}{2} \Tr{\alpha^g (\alpha^{g})^Tw^g_{\xi \xi}} +\circ(\eps^2),
\end{align*}

where the diffusion coefficient is given by
\begin{align}\label{def alphag}
\alpha^g(t,s,z):=\left[\left(I_d-\ybf^g_z(t,s,z)\1_d^T\right)
{\rm{diag}}[\ybf^g(t,s,z)]-(\ybf^g_s)^T(t,s,z){\rm{diag}}[s]\right]\sigma(t,s).
\end{align}

This calculation highlights the role played by the so-called fast variable $ \xi$. Indeed any of the second order derivatives of $w^g$ with respect to $s$ or $ y$, corresponds to a second-order derivative of $ \hat{w}^g$ scaled by $1/\eps^2$. These terms are then exactly of the same order as the one obtained above. Finally it is obvious that, using the definition of $ \cbf^g$ in \eqref{def of c}:
\begin{align*}
\kappa\tilde{U}_1(v^{\eps,g}_x)-\kappa\tilde{U}_1(v^g_x)+\kappa\cbf^g (v^{\eps,g}_x-v^g_x)&=\kappa\tilde{U}_1(v^{\eps,g}_x)-\kappa\tilde{U}_1(v^g_x)-\eps^2\kappa \cbf^g u^g_z +\circ(\eps^2)=\circ(\eps^2).
\end{align*}

Combining these approximations and putting them into the drift condition of \eqref{e.dpp}, we obtain that $ u^g$ must be solution of the second corrector equation:
\begin{equation}\label{eq:corr2}\begin{cases}
\displaystyle\Ac^g u^g
=
a^g(t,s,z), \ (t,s,z)\in[0,T)\times(0,+\infty)^{d+1}\\
u^g(T,s,z)=0,\ (s,z)\in(0,+\infty)^{d+1},
\end{cases}\end{equation}
where the differential operator $\Ac^g$ is defined by
$$\Ac^g u
:= ku-\Lc^0 u
 - \big(rz + \ybf^g\cdot(\mu-r\1_d) - \kappa\cbf^g\big) 
   u_z
 -\frac12 |\sigma^{\rm T} \ybf^g|^2 \,u_{zz}
 -\sigma\sigma^{\rm T}\ybf^g\cdot\D_{sz} u,$$

and the function $a^g$ is the second component of the solution $(w^g,a^g)$ of the first corrector equation:
\begin{align}\label{firstcor}
\nonumber &\underset{ (i,j)\in\mathcal I}{\max}\max\left\{\frac{\abs{ \sigma^T(t,s)\xi}^2}{2}v_{zz}^g(t,s,z)
-\frac12\Tr{\alpha^g(\alpha^g)^T(t,s,z) w^g_{\xi\xi}(t,s,z,\xi)}+a^g(t,s,z)\ , \right.\\
&\hspace{5.9em}\left. -\lambda^{i,j}v^g_z(t,s,z)+ \frac{\partial w^g}{\partial \xi^i}(t,s,z,\xi)
- \frac{\partial w^g}{\partial \xi^j}(t,s,z,\xi) \right\}=0,\ \xi\in\R^d.
\end{align}

\begin{Remark}\label{rem.dig}
Notice that we consider naturally \eqref{firstcor} only on $ [0,T)\times \R_+ \times \R_+$, because since the value function is known at time $T$, its expansion is trivial. 
Since we enforce that the function $u^g$ solution of the second corrector equation \reff{eq:corr2} is null at time $T$, it would seem reasonable to think that the expansion \eqref{expansion ve} also holds at time $T$. However, as we will see in our proofs, this will usually only be true if the Merton value function and the corresponding optimal strategy are smooth enough at time $T$. If explosions are allowed at time $T$ {\rm(}which, as pointed out in Section \ref{sec:disc}, can happen for the derivatives of $\ybf^g$ if $g$ is a Call option{\rm)}, then the remainder in the expansion \reff{expansion ve} can become unbounded near $T$. In the previous works by Bichuch \cite{bichuch} and Bouchard, Moreau and Soner \cite{bms}, strong regularity on $v^g$ up to time $T$ was assumed {\rm(}which implies then that the payoff $g$ has to be regular{\rm)}, in order to prevent $\ybf^g$ and several of its derivatives from exploding at $T$. With our method however, this is no longer needed. We refer the reader to Section \ref{sec.assump} for more details on our assumptions.
\end{Remark}

Finally, we recall from \cite{st} and \cite{pst} the following normalization. Set
$$
\eta^g(t,s,z):=-\frac{v^g_z(t,s,z)}{v^g_{zz}(t,s,z)},\ 
\rho:=\frac{\xi}{\eta^g(t,s,z)},\ \overline{w}^g(t,s,z,\rho)
:=\frac{w^g(t,s,z,\eta^g(t,s,z)\rho)}{\eta^g(t,s,z)v^g_z(t,s,z)},
$$ 
$$
\overline{a}^g(t,s,z):=\frac{a^g(t,s,z)}{\eta^g(t,s,z)v^g_z(t,s,z)}, \ 
\bar{\alpha}^g(t,s,z):=\frac{\alpha^g(t,s,z)}{\eta^g(t,s,z)},
$$
so that the corrector equations with variable 
$\rho \in \R^d$ have the form,
\begin{align}
\begin{cases}
\displaystyle\displaystyle\underset{(i,j)\in\mathcal I}{\max}\max\left\{-\frac{\abs{ \sigma^T(t,s)\rho}^2}{2}-
\frac12\Tr{\bar\alpha^g(\bar\alpha^g)^T(t,s,z) \overline{w}^g_{\rho\rho}(t,s,z,\rho)}+\overline a^g(t,s,z)\ ;\right.
\\
\left.\displaystyle\hspace{5.9em}-\lambda^{i,j}+ \frac{\partial \overline{w}^g}{\partial \rho_i}(t,s,z,\rho)- 
\frac{\partial \overline{w}^g}{\partial \rho_j}(t,s,z,\rho) \right\}=0
\label{eq:corrector}\\[0.8em]
\mathcal A^gu^g(t,s,z)=v^g_z(t,s,z)\eta^g(t,s,z)\overline{a}^g(t,s,z).
\end{cases}
\end{align}
 We emphasize that the first corrector equation \reff{eq:corrector} is an equation for the variable $\xi$, $(t,s,z)$ are only parameters. Moreover, the wellposedness of this equation has been obtained in \cite{pst}. We recall below the properties of $\overline{w}^g$ that we will use. Before stating the result, let us define the following closed convex subset of $\mathbb R^d$, and the corresponding support function, with the convention that $\rho_0=0$
 \b*
 C
 :=
 \left\{\rho\in\mathbb R^d:
        -\lambda^{j,i}\leq \rho_i-\rho_j\leq \lambda^{i,j}, 
        \ (i,j)\in\mathcal I
 \right\},
 &\delta_C(\rho):=\underset{u\in C}{\sup}\ u \cdot \rho,& 
 \rho\in\mathbb R^d.
\e*

 \begin{Proposition}\label{prop:wg}
Assume that $\bar{\alpha}^g(t,s,z)$ is non-degenerate for any $(t,s,z)\in[0,T)\times\R_+^d\times\R_+$. Then there exists a unique solution $(\overline{w}^g,\overline{a}^g)$ of the 
equation \reff{eq:corrector}, such that $\overline{a}^g\in\R_+$, $\rho\longmapsto \overline{w}^g(t,s,z,\rho)$ is $C^1$ in $\R^d$ with a Lipschitz gradient, such that the following growth condition holds $$\underset{|\rho|\to+\infty}{\lim}\frac{\overline{w}^g(t,s,z,\rho)}{\delta_C(\rho)}=1,$$ and such that $\overline{w}^g(\cdot,0)=0$. Moreover, for any $(t,s,z)\in[0,T)\times\R_+^d\times\R_+$, $\overline{w}^g(t,s,z,\cdot)$ is convex, the set $\mathcal O_0^g(t,s,z):=\left\{\rho\in\mathbb R^d,
	\ \overline{w}^g_\rho(t,s,z,\rho)\in\text{$\rm{int}$}(C)\right\}$
	 is open and bounded, the map $\rho\longmapsto\overline{w}^g(t,s,z,\rho)$ is $C^\infty$ on $\mathcal O^g_0(t,s,z)$,
 $\overline{w}^g(t,s,z,\rho)$ attains its minimum in $\rho$ at some point $\rho^*(t,s,z)$ in $\mathcal O^g_0(t,s,z)$, and there is a constant $M>0$ s.t. $0\leq \overline{w}^g_{\rho\rho}(t,s,z,\rho)\leq M{\bf1}_{\overline{\mathcal O}^g_0(t,s,z)}(\rho)$ for a.e. $\rho\in\mathbb R^d$.
 \end{Proposition}
 
Of course, under suitable regularity assumptions on $\eta^g$ and $v^g_z$, $w^g$ satisfies similar properties.

\subsection{Formal asymptotics for the utility indifference price}
We now develop an expansion for $ p^{\epsilon,g}$, using the expansion of $ v^{\eps,g}$ defined in \eqref{expansion ve}. We first recall that, at least formally, for $\vartheta=0$ or $g$
\begin{align*}
v^{\eps,\vartheta}(t,s,x,y) &= v^\vartheta(t,s,x,y)-\epsilon^2 u^\vartheta(t,s,z)+\circ(\epsilon^2).
\end{align*}
Then, at least if $v^{\eps,g}$ is increasing with respect to $x$, $ p^{\epsilon,g}(t,s,x)$ should be such that:
$$ v^{\eps,g}(t,s,x+p^{\epsilon,g}(t,s,x),0)=v^{\eps,0}(t,s,x,0).$$

We conjecture (and will prove under natural assumptions) that $p^{\eps,g}$ satisfies the following expansion
\begin{equation}\label{asymp_prix}
p^{\eps,g}(t,s,x) = p^{g}(t,s,x) + \eps^2h^g(t,s,x)+\circ(\eps^2),
\end{equation}
for some function $h^g$ to be determined. Using \eqref{expansion ve}, we obtain formally
\begin{align*}
v^0(t,s,x)-\epsilon^2 u^0(t,s,x)+\circ(\epsilon^2) =&\ v^g(t,s,x+p^g(t,s,x))+\eps^2v^g_x(t,s,x+p^g(t,s,x))h^g(t,s,x)\\
&-\epsilon^2 u^g(t,s,x+p^g(t,s,x))+\circ(\epsilon^2).
\end{align*}
Since by definition we have $v^0(t,s,x,0)=v^g(t,s,x+p^g(t,s,x))$, we deduce:
$$h^g(t,s,x)=\frac{u^g(t,s,x+p^g(t,s,x))-u^0(t,s,x)}{v^g_x(t,s,x+p^g(t,s,x))}.$$

\section{Main results}\label{sec.main}
We recall from \cite{st} the following notations. For any $ f(s,x,y)$, we define the change of variable:
$$ \hat f(t,s,z,\xi) := f\big(t,s,z-\eps \xi\cdot \1_d-\ybf^g(t,s,z)\cdot\1_d,\eps\xi+\ybf^g(t,s,z)\big).$$

We then define
\be
\bar{u}^{\eps,g}(t,s,x,y)
&:=&
\frac{v^g(t,s,z)-v^{\eps,g}(t,s,x,y)}{\eps^2},\ s\in\mathbb R_+^d, ~~(x,y)\in K_\eps,
\ee
and its relaxed semi-limits:
$$
u^{g,*}(t,s,x,y)
:=
\underset{(\eps,t',s',x',y')\longrightarrow (0,t,s,x,y)}{\overline{\lim}}\bar{u}^{\eps,g}(t',s',x',y'), $$
$$ u_*^g(t,s,x,y)
:=\underset{(\eps,t',s',x',y')\longrightarrow (0,t,s,x,y)}{\underline{\lim}}\bar{u}^{\eps,g}(t',s',x',y').$$
Finally, we introduce, $u^{\eps,g}(t,s,x,y):=\bar{u}^{\eps,g}(t,s,x,y)-\eps^2w^g(t,s,z,\xi), \ s\in\mathbb R_+^d,\ (x,y)\in K_\eps.$

%

\subsection{Assumptions}\label{sec.assump}
In all the following, we consider payoff functions $g$ and functions $r$, $\mu$ and $\sigma$ such that the following four assumptions hold.

\begin{Assumption}[Smoothness of $ \ybf^g$, $ v^g$, $ \ybf^0$ and $ v^0$]\label{assumption 1 smoothness}
For $\vartheta=0\text{ or }g,$ we have

\vspace{0.25em}
{\rm(i)} The map $v^\vartheta(t,s,z)$ is $C^{1,2,2}$ in $[0,T)\times(0,+\infty)^{d+1}$ and $C^{0,0,0}$ in $[0,T]\times(0,+\infty)^{d+1}$. Moreover, for any $(t,s)\in[0,T]\times(0,+\infty)^d$, the map $z\longmapsto v^\vartheta(t,s,z)$ is $C^1$ in $(0,+\infty)$ and we have
\begin{align*}
v^\vartheta_z(t,s,z)&>0,\ (t,s,z)\in[0,T]\times(0,+\infty)^{d+1},\\
\abs{v^\vartheta_{zz}}(t,s,z)&\leq \frac{C(s,z)}{(T-t)^{1-\mu}}, \ (t,s,z)\in[0,T)\times(0,+\infty)^{d+1},
\end{align*}
for some continuous function $C$ and some $\mu\in(0,1]$.

\vspace{0.25em}
{\rm(ii)} The map $\ybf^\vartheta(t,s,z)$ is $C^{1,2,2}$ in $[0,T)\times(0,+\infty)^{d+1}$ and $C^{0,0,0}$ in $[0,T]\times(0,+\infty)^{d+1}$. Moreover, for any $(t,s)\in[0,T]\times(0,+\infty)^d$, the map $z\longmapsto \ybf^\vartheta(t,s,z)$ is $C^1$ in $(0,+\infty)$ and there exist some constants $(c_0,c_1,\eta)\in(0,+\infty)\times(0,+\infty)\times(0,1]$ such that for any $(t,s,z)\in[0,T]\times(0,+\infty)^{d+1}$
$$\ybf^{\vartheta,i}_z(t,s,z)>0,\ c_0\leq \ybf^\vartheta_z(t,s,z)\cdot \1_d
	\text{ and }
	\left[\alpha^\vartheta(\alpha^\vartheta)^{\rm T}\right](t,s,z)\ge c_1 I_d,
	\ 1\leq i\leq d,$$
and for any $(t,s,z)\in[0,T)\times(0,+\infty)^{d+1}$
$$\left[\abs{\ybf^\vartheta_t}+\abs{\ybf^\vartheta_s}+\abs{\ybf^\vartheta_{zz}}+\abs{\ybf^\vartheta_{sz}}+\abs{\ybf^\vartheta_{ss}}\right](t,s,z)\leq\frac{C(s,z)}{(T-t)^{1-\eta}},$$
for some continuous function $C$.
\end{Assumption}

\begin{Remark}\label{rem.ygz}
It can be readily checked that if it happens that $\ybf^\vartheta$ does not depend on $z$, then even though Assumption \ref{assumption 1 smoothness}{\rm(}ii{\rm)} does not hold {\rm(}since $\ybf^\vartheta_z=0${\rm)}, all our subsequent proofs still go through. It will be important for us later on when we treat the case of exponential utility in Section \ref{sec.exp}.
\end{Remark}


\begin{Remark}
We assumed here that the first-order derivatives of $v^\vartheta$ and $\ybf^\vartheta$ with respect to $z$ are well defined at $T$, unlike the other derivatives which may not exist at $T$. This is basically due to the so-called remainder estimate that we obtain in Lemma \ref{lem remainder estimate}, since these terms are the only ones which appear in conjunction with $\tilde U_1$ and its derivatives. We may have let them explode at time $T$ with a certain speed, but we would then have needed to control the growth at infinity of $\tilde U_1$ and its derivatives. The above assumptions being already technical, we refrained from doing so, but we insist on the fact that in particular examples, our general conditions may be readily improved simply by looking at the remainder estimate obtained and using it in the proof of the viscosity subsolution property at the boundary in Section \ref{sec.boundary}.
\end{Remark}


We now state an assumption on the regularity of the solution of the first corrector equation with respect to the parameters $(t,s,z)$.
\begin{Assumption}[First corrector equation: regularity on the parameters]\label{assumption 4 bound on derivatives}
For $\vartheta=0\text{ or }g,$ the set $ \mathcal{O}^\vartheta_0(t,s,z)$ {\rm(}see Proposition \ref{prop:wg}{\rm)} as well as $ a^\vartheta(t,s,z)$ and $\rho^*(t,s,z)$ are continuous in $ (t,s,z)\in[0,T)\times(0,+\infty)^{d+1}$. Moreover, both $w^\vartheta$ and $\tilde w^\vartheta(\cdot,\xi):=w^\vartheta(\cdot,\xi)-w^\vartheta(\cdot,\eta^\vartheta(\cdot)\rho^*(\cdot))$ are $ C^{1,2,2}$ in $[0,T)\times(0,+\infty)^{d+1}$ and satisfy for any $(t,s,z,\xi)\in[0,T)\times(0,+\infty)^{d+1}\times\R^d$
\begin{align}
& \left( | \varpi_t|+|\varpi_s|+|\varpi_{ss}|+|\varpi_z|+|\varpi_{sz}|+|\varpi_{zz}|  \right)(t,s,z,\xi)\le C(t,s,z) \left( 1+|\xi| \right)\\
&\left( |\varpi_\xi|+|\varpi_{s\xi}|+|\varpi_{z\xi}|  \right)(t,s,z,\xi) \le C(t,s,z),
\end{align}
for $\varpi=w^\vartheta\text{ or }\tilde w^\vartheta$ and for some continuous function $ C(t,s,z)$.
\end{Assumption}
The above assumption can be readily verified in dimension $d=1$ for which the functions $w^g$ and $a^g$ are given explicitly in terms of the Merton value function and its derivatives. However, it would be a very difficult task to verify it in the general framework considered here. Our intention is simply to state directly what are the kind of regularity we must assume to recover the expansions, and then these can be checked on particular examples. For further reference, we also insist on the fact that by definition, the function $\tilde w^g$ is non-negative.

\vspace{0.3em}
A fundamental step in any homogenization proof is to show that the correctors are uniformly locally bounded. In our context, this means that we need to show that $\bar{u}^{\eps,g}$ is locally uniformly bounded. Since by definition it is a positive quantity, we only need an upper bound. We put this as an assumption. 
\begin{Assumption}[Local bound of $ \bar{u}^g$]\label{assumption 2 local bound}
The family of functions $ \bar{u}^{\eps,g}$ is locally uniformly bounded from above.
\end{Assumption}

Of course, one could argue that we are avoiding a major problem here. However, exactly as for the previous assumption, given the level of generality we are working with, verifying that it holds for generic models goes beyond the scope of this paper. However, we will show later on in Section \ref{sec:export} that when utilities are exponential and under an additional regularity assumption (which is always satisfied when $d=1$), Assumption \ref{assumption 2 local bound} is satisfied. Let us nonetheless sketch the arguments of an alternative approach which would not require in general this additional regularity assumption. First of all, notice that we can without loss of generality only consider the case where all the $\lambda^{i,j}=+\infty$ for $1\leq i,j\leq d$ and where $\kappa=0$ ({\it i.e.} no consumption allowed). Indeed, the corresponding value function is clearly smaller than $v^{\eps,g}$, and thus the corresponding $\bar u^{\eps,g}$ is greater than the one for which we want to find an upper bound. Hence, it suffices to consider this case.

\vspace{0.25em}
The first step is then to construct a regular viscosity sub-solution to the dynamic programming equation \reff{e.dpp} which has the form
$$V^{\eps,K}(t,s,z,\xi):=v^g(t,s,z)-\eps^2Ku^g(t,s,z)-\eps^4w^g(t,s,z,\xi),$$
where $u^g$ and $w^g$ are the solutions to the corrector equations and where $K$ is a large constant. Indeed, using comparison for \reff{e.dpp}, this would then imply that $V^{\eps,K}\leq v^{\eps,g}$, from which we can immediately deduce the required upper bound for $\bar{u}^{\eps,g}$.

\vspace{0.25em}
Of course, the first problem would then be that we are not sure that $u^g$ and $w^g$ are smooth. For $u^g$, since it is the solution to a linear PDE, this could be readily checked as soon as we have enough regularity on $a^g$. However, for $w^g$ as soon as $d>1$, we cannot reasonably expect it to be more than $C^1$ in $\xi$, since variational inequalities with gradient constraints in dimension greater than $1$ are generally not $C^2$. Nonetheless, this issue can easily be solved by replacing $w^g$ by a function $W^g(t,s,z,\xi)=\widehat{w}^g(\eta^g(t,s,z)\rho)/(\eta^g(t,s,z)v_z^g(t,s,z))$, where $\widehat w^g$ is the first component of the solution $(\widehat w^g,\widehat a^g)$ to the following equation,
\begin{equation}
\label{e.what}
\underset{0\leq i\leq d}
{\max}\max\left\{-\frac{c^*_1\abs{\rho}^2}{2}
-\frac{c^*_2}{2}\Delta \widehat{w}(\rho)+\widehat{a}\ ,\ -\widehat{\lambda}^{i}
+ \frac{\partial \widehat{w}}{\partial \rho_i}(\rho)\ ,\ -\tilde{\lambda}^{i}
- \frac{\partial \widehat{w}}{\partial \rho_i}(\rho)
\right\}=0,
\end{equation}
with the normalization $\widehat{w}(0)=0$ and where the positive constants
$c^*_1, c^*_2, \widehat{\lambda}^{i}, \tilde{\lambda}^{i}$
are s.t. for some constant $M>0$
$$
c^*_1 I_d \ge \sigma \sigma^T, \qquad
c^*_2 I_d \ge \bar{\alpha}^g (\bar{\alpha}^T)^g,
\qquad
\widehat{\lambda}^i=
\tilde{\lambda}^i= M \overline \lambda:=M
\underset{(i,j)\in\mathcal I}{\max}\lambda^{i,j}.
$$
Then, as shown in \cite{pst}, the unique solution $\widehat{w}^g$ 
is given as, $
\widehat{w}^g(\rho) = \sum_{i=1}^{d}\ \widetilde{w}^g_i(\rho_i),
$
where $\widetilde{w}^g_i$ is the explicit solution of 
the corresponding 1-d problem, which is known explicitly and is $C^2$.

\vspace{0.25em}
To prove the viscosity sub solution property, one can then argue exactly as in the proof of Lemma $3.1$ in \cite{pst}. This proof is made under assumptions ensuring homotheticity in $z$ of the functions appearing, but the general approach will be valid in other cases as well, albeit with more complicated computations. For instance, in the case where $U_2$ is an exponential utility, and the frictionless market is the Black-Scholes model, the dependence in $(t,s,z)$ of all quantities involved is known explicitly (see the formulas in Section \ref{sec:disc} for details). Basically, one has to check that for $M$ large enough the gradient constraints
\begin{align*}&\Lambda^\eps_{i,0}\cdot(V^{\eps,K}_x, V^{\eps,K}_{y_i}) \le 0,
\quad
{\mbox{ holds whenever}} 
\quad
M \overline \lambda
+ \frac{\partial \widehat{w}}{\partial \rho_i}(\rho) \le 0,\\
&\Lambda^\eps_{0,i}\cdot(V^{\eps,K}_x, V^{\eps,K}_{y_i}) \le 0,
\quad
{\mbox{holds whenever}} 
\quad
-M \overline \lambda
+ \frac{\partial \widehat{w}}{\partial \rho_i}(\rho) \le 0.
\end{align*}
Then, by choosing $K$ large enough, one has to show that the diffusion operator in \reff{e.dpp} applied to $V^{\eps,K}$ is a non-positive quantity, which would then give the desired result. 

\vspace{0.25em}
Since we assumed that $\bar u^{\eps,g}$ is uniformly locally bounded, we can define for $ (t_0,s_0,x_0,y_0)\in [0,T] \times (0,\infty)^d \times \R \times \R^d$ with $ x_0+y_0 \cdot \1_d >0$
\begin{align}
b(t_0,s_0,x_0,y_0):=\sup \big\{ u^{\eps,g}(t,s,x,y) \ : (t,s,x,y) \in B_{r_0}(t_0,s_0,x_0,y_0), \ \eps \in (0,\eps_0] \big\}.
\end{align} 
Then using the continuity of $ w^g$, there exists $ r_0(t_0,s_0,x_0,y_0)>0 $ and $ \eps_0(t_0,s_0,x_0,y_0)>0$ such that $ b(t_0,s_0,x_0,y_0)<\infty $.

\vspace{0.25em}
Our final assumption ensures that we have a comparison theorem for the second corrector equation.
\begin{Assumption}[Second corrector equation: comparison]\label{assumption 3 comparison}
For $\vartheta=g$ or $0$, there exists a set of functions $\mathcal C$ which contains $u^{*,\vartheta}$ and $u_*^\vartheta$ and such that $u_1\geq u_2$ on $[0,T]\times (0,+\infty)^{d+1}$, whenever $u_1$ {\rm(}resp. $u_2${\rm)} is a l.s.c. {\rm(}resp. u.s.c.{\rm)} viscosity super-solution {\rm(}resp. sub-solution{\rm)} of \reff{eq:corr2} in $\mathcal C$.
\end{Assumption}

Once again, we will not attempt to verify this assumption. Nonetheless, we insist on the fact that the PDE \reff{eq:corr2} is linear, so that we can reasonably expect that a comparison theorem on the class of functions with polynomial growth will hold as soon as $a^g$ itself has polynomial growth.

\subsection{The results}

\begin{Theorem}[Convergence of $ u^{\epsilon,g}$]\label{th convergence de ueps}
Under assumptions \ref{assumption 1 smoothness}, \ref{assumption 4 bound on derivatives}, \ref{assumption 2 local bound} and \ref{assumption 3 comparison}, the sequence $\bar{u}^{\eps,g}$ converges locally uniformly to a function $u^g$ depending only on $ (t,s,z)$ and which is the unique viscosity solution of \reff{eq:corr2}.
\end{Theorem}

The proof is relegated to Section \ref{sect: proof of th convergence}.

\vspace{0.25em}

\begin{Theorem}[Expansion of the utility indifference price]

Under assumptions \ref{assumption 1 smoothness}, \ref{assumption 4 bound on derivatives}, \ref{assumption 2 local bound} and \ref{assumption 3 comparison}, we have for all $(t,s,x)$:

$$\frac{p^{\eps,g}(t,s,x)-p^g(t,s,x)-\eps^2 h^g(t,s,x)}{\eps^2}\longrightarrow1, \text{ locally uniformly as $\eps\longrightarrow0$},$$ 
where $ h^g(t,s,x):=\frac{u^g(t,s,x+p^g(t,s,x))-u^0(t,s,x)}{ v^g_z(t,s,x+p^g(t,s,x))}.$
\end{Theorem}

\vspace{0.25em}

\textbf{Proof} \textit{Step 1:} We first show that $ p^{g}$ is continuous in $(t,s,z)$. Since $ v^g $ and $ v^0$ are $ C^2$ and $ v^g$ is partially strictly concave w.r.t. $z$, we have that $ v^g_z>0$. The continuity of $ p^g$ follows easily. Indeed assume on the contrary that there exists $ (t_0,s_0,x_0)$, $ \eps>0$ and $ (t_n,s_n,x_n)\longrightarrow (t_0,s_0,x_0)$ s.t. $$ \left| p^g(t_n,s_n,x_n)-p^g(t_0,s_0,x_0)\right|>\eps.$$ W.l.o.g., we can assume that $ p^g(t_n,s_n,x_n)>p^g(t_0,s_0,x_0)+\eps$. Then we have by definition of $ p^g$ and the fact that $ v^g$ is increasing w.r.t. the $x$-variable that for all $ n \geq 0$
\begin{align}\label{DL eq interm 1}
v^g(t_n,s_n,x_n+p^g(t_0,s_0,x_0)+\eps)&<v^0(t_n,s_n,x_n),\\
\label{DL eq interm 2}
v^g(t_0,s_0,x_0+p^g(t_0,s_0,x_0)+\eps)&>v^0(t_0,s_0,x_0).
\end{align}

Then by continuity of $ v^g$ and $ v^0$ we obtain from \eqref{DL eq interm 1} that $ v^g(t_0,s_0,x_0+p^g(t_0,s_0,x_0)+\eps) \geq v^0(t_0,s_0,x_0),$ which contradicts \eqref{DL eq interm 2}.

\vspace{0.25em}
\textit{Step 2:} Let $(t_0,s_0,x_0)\in [0,T]\times \R^d_+ \times \R_+$. We consider $ r>0$ s.t. on $\bar{B}_r(t_0,s_0,x_0)$, the quantity $ u^{\eps,g}(t,s,x+p^g(t,s,x))$ and $ u^{\eps,0}(t,s,x)$ converge uniformly to, respectively, $ u^{g}(t,s,x+p^g(t,s,x))$ and $ u^0(t,s,x)$. Notice that the existence of $r$ is guaranteed by the result of Theorem \ref{th convergence de ueps}, together with the fact that $p^g$ is continuous. We use the notations $ p^g$ (resp. $h^g$) for $ p^g(t,s,x)$ (resp. $ h^g(t,s,x) $) for simplicity. For any $ \delta \in(-1,1)$, we have uniformly on $ \bar{B}_r(t_0,s_0,x_0)$:
\begin{align*}
v^{\eps,g}(t,s,x+p^g+\eps^2 h^g+\eps^2 \delta,0) =&\ v^g(t,s,x+p^g)+ \eps^2 \left( h^g+\delta \right) v^g_x(t,s,x+p^g)\\
& -\eps^2 u(t,s,x+p^g)+ \circ(\eps^2)\\
=&\ v^{\eps,0}(t,s,x)+ \eps^2 \delta v^g_z(t,s,x+p^g) +\circ(\eps^2).
\end{align*}

Hence, the following holds uniformly on $ \bar{B}_r(t_0,s_0,x_0)$
\begin{align}\label{DL eq interm 3}
\frac{v^{\eps,g}(t,s,x+p^g+\eps^2 h^g+\eps^2 \delta)-v^{\eps,0}(t,s,x)}{\eps^2}=\delta v^g_z(t,s,x+p^g) +\circ(1).
\end{align}

We now claim that for any $ \delta>0$, there is some $\eps^*(\delta)$ such that we have for $ \eps \le \eps^*(\delta) $ that on $ \bar B_r(t_0,s_0,x_0)$:
$$ p^g(t,s,x)+\eps^2 h^g(t,s,x)-\eps^2 \delta \le p^{\eps,g}(t,s,x) \le p^g(t,s,x)+\eps^2 h^g(t,s,x)+\eps^2\delta,$$
which implies directly the required result. It remains to prove the claim. Assume on the contrary that we have $ \delta>0$ and $ (\eps_n,t_n,s_n,x_n)$, where for all $ n \geq 0$, $ (t_n,s_n,x_n)\in \bar{B}_r(t_0,s_0,x_0) $ and $ \eps_n \longrightarrow 0$, such that, for example, for all $ n\geq 0$
$$ p^{\eps_n,g}(t_n,s_n,x_n)>p^g(t_n,s_n,x_n)+\eps_n^2h^g(t_n,s_n,x_n)+\eps_n^2 \delta.$$

Then by definition of $ p^{\eps,g}$, we have that
$$ v^{\eps_n,g}(t_n,s_n,x_n+p^g+\eps_n^2 h^g+\eps_n^2 \delta)-v^{\eps_n,0}(t_n,s_n,x_n)\leq 0,$$

which contradicts \eqref{DL eq interm 3} for $ n$ large enough. The other inequality can be shown similarly.

\vspace{-2em}

\begin{flushright}$\Box$\end{flushright}

\vspace{-1.5em}

\section{Examples and applications}\label{sec:disc}

In this Section we will specialize our discussion to a simpler case, in order to highlight how our method allows not only to recover existing results but to go beyond them. Throughout the section, we place assume a Black-Scholes dynamic for the risky asset, that is $\mu$, $\sigma$ and $r$ are constant. The investor also aims at solving the following versions of the stochastic control problems \reff{def:ve} and \reff{def:v0}
 \begin{align}\label{def:veBS}
 v^{\epsilon,g}(t,s,x,y)
& :=
 \sup_{(c,L) \in \Theta^\ep(t,s,x,y)}\
 \E_t\left[\int_t^T \kappa U_1(c_\xi)d\xi+U_2\left(\ell^\eps\left(X^{t,s,x,y}_T,Y^{t,s,x,y}_T\right)-g(S_T)\right) \right],\\
 \label{def:v0BS}
 v^{g}(t,s,z)
&:=
 \sup_{(c,\theta) \in \Theta^0(t,s,z)}\
 \E_t\left[\int_t^T \kappa U_1(c_\xi)d\xi+U_2\left(Z_T^{\theta,t,s,z}-g(S_T)\right) \right],
 \end{align}
 corresponding to the case $k=0$ in \reff{def:ve}. We will now show, for a particular choice of utility functions, that we can calculate almost explicitly all the quantities involved in the asymptotic expansion \reff{asymp_prix}, as well as check that all our assumptions hold under certain explicit conditions. For further reference and use, we recall that in this setting, the so-called Black-Scholes price of the claim $g$, denoted by $V^g$ is given by
 $$V^g(t,s):=\mathbb E_t^\Q\left[e^{-r(T-t)}g(S_T^{t,s})\right],$$
 where $\Q$ is the so-called risk neutral pricing measure defined by
 $$\frac{d\Q}{d\P}=\mathcal E\left(-\sigma^{-1}(\mu-r{\bf 1}_d)\cdot W_T\right).$$
 
 Moreover, it is a well known result that as soon as $g$ has sub-exponential growth at infinity, $V^g$ is also the unique (classical) solution to the following PDE  \begin{align}\label{eq:BS}
 -V^g_t-r{\bf 1}_d\cdot \D_sV^g_s-\frac12\Tr{\sigma\sigma^T\D_{ss}V^g}+rV^g=0,\ (t,s)\in[0,T)\times(0,+\infty)^d,\ V^g(T,\cdot)=g(\cdot).
 \end{align}

We also recall that in the special case where $d=1$, \cite{st} gave an explicit solution to the 2nd corrector equation (see their Section 4.1), 
\begin{align}\label{eq:www}
\nonumber&w^g(t,s,z,\xi)=\\
&\begin{cases}\displaystyle \left(v^g_z\left[-\frac{\sigma^2}{12(\eta^g)^2(\alpha^g)^2}\xi^4+\frac{\sigma^2}{2(\eta^g)^2(\alpha^g)^2}\xi_0^2\xi^2+\frac{\lambda^{1,0}-\lambda^{0,1}}{2}\xi\right]\right)(t,s,z),\text{ if $\abs{\xi}\leq\xi_0$}\\[0.8em]
\displaystyle v^g_z(t,s,z)\left[-\frac{3}{16}\left(\lambda^{1,0}+\lambda^{0,1}\right)\xi_0-\lambda^{0,1}\xi\right],\text{ if $\xi\leq-\xi_0$}\\[0.8em]
\displaystyle v^g_z(t,s,z)\left[-\frac{3}{16}\left(\lambda^{1,0}+\lambda^{0,1}\right)\xi_0+\lambda^{1,0}\xi\right],\text{ if $\xi\geq\xi_0$},
\end{cases}
\end{align}
where 
$$\xi_0:=\xi_0(t,s,z):=\eta^g(t,s,z)\left(\frac34\frac{(\alpha^g)^2(t,s,z)}{(\eta^g)^2(t,s,z)\sigma^2}\left(\lambda^{1,0}+\lambda^{0,1}\right)\right)^{1/3},$$
which in turn allows us to have an explicit form for the function $a^g(t,s,z)$ in terms of the Merton value function
\begin{equation}\label{eq:a}
a^g(t,s,z)=\frac{\sigma^2v^g_z(t,s,z)}{2\eta^g(t,s,z)}\xi_0^2(t,s,z).
\end{equation}

Let us now assume throughout this section that $U_1(x)=U_2(x) = -e^{-\gamma x},\text{ for some $\gamma>0$.}$

\subsection{Derivation of the expansion}\label{sec.exp}

We start by giving the solution to the Merton problem corresponding to $\eps=0$. In the case $\kappa=0$ (no consumption), the solution when $d=1$ can be found for instance in \cite{bichuch} (see also the references therein). The generalization to the consumption and multidimensional case is an easy (but lengthy) exercise, so that we omit its proof.

\begin{Proposition}\label{prop:expBS}
The value function for the stochastic control problem \reff{def:v0BS} is given for any $(t,s,z)\in[0,T]\times\R_+\times\R_+$
by
$$v^g(t,s,z)=-\exp\big(-\gamma v_1(t)\left(z-V^g(t,s)\right)+v_2(t)\big),$$
provided that $\kappa v_z^g(t,s,z)\leq \gamma$, and where
\begin{align*}
v_1(t):=&\ \frac{r}{\kappa+e^{-r(T-t)}(r-\kappa)}\\
v_2(t):=&\ \frac{1}{\kappa e^{r(T-t)}+r-\kappa}\left[\frac12(\sigma\sigma^T)^{-1}(\mu-r{\bf 1}_d)\cdot(\mu-r{\bf 1}_d)(\kappa-r)(T-t)-\kappa r(T-t)e^{r(T-t)}\right.\\
&\left.\hspace{8em}-\kappa\left(e^{r(T-t)}-1\right)\left(\frac{1}{2r}(\sigma\sigma^T)^{-1}(\mu-r{\bf 1}_d)\cdot(\mu-r{\bf 1}_d)+\log(r)-2\right)\right.\\
&\left.\hspace{8em}+\kappa\frac{\log^2\left(\kappa e^{r(T-t)}+r-\kappa\right)-\log^2\left(r\right)}{2}\right].
\end{align*}
Moreover, the optimal trading strategy and consumption are given by
\begin{align*}
&\ybf^g(t,s):=\D_sV^g(t,s)+\frac1\gamma(\sigma\sigma^T)^{-1}(\mu-r{\bf 1}_d)v_1^{-1}(t)\\
 &\cbf^g(t,s,z):=\kappa\left(-\frac1\gamma\left(\log(v_1(t))+v_2(t)\right)+v_1(t)(z-V^g(t,s))\right).
 \end{align*}
\end{Proposition}

\begin{Remark}
It can be checked directly that when $\kappa=0$, the above reduces to the formula given in Remark 3.4 of \cite{bichuch}. Moreover, the condition $\kappa v^g_z\leq \gamma$ is here to ensure that $\cbf^g\geq 0$. When $\kappa=1$, it can be verified directly that it is for instance satisfied if $z\geq V^g(t,s)$ and $r$ is small enough.
\end{Remark}

Using the above proposition, we recover the expected result that the utility indifference price $p^g(t,s,z)$ does not depend on $z$ in this case, and is simply given by the Black-Scholes price $V^g(t,s)$ of the contingent claim $g$. We refer the reader to Theorem 1 and Section 3 in \cite{dpz} for further details on this general result. Next, we deduce immediately that the matrix $\alpha^g$ is independent of $z$ and that $\eta^g$ is independent of $(s,z)$ and are given by
$$\alpha^g(t,s)=\left(\frac{1}{\gamma v_1(t)}{\rm diag}\left[(\sigma\sigma^T)^{-1}(\mu-r{\bf 1}_d)\right]-{\rm diag}[s]V_{ss}^g{\rm diag}[s]\right)\sigma,\ \ \eta^g(t)=\frac{1}{\gamma v_1(t)}.$$

Hence, it is clear from the first corrector equation in \reff{eq:corrector} that $\overline w^g$ does not depend on $z$, so that we have the factorization
$$a^g(t,s,z)=-v^g(t,s,z)\overline{a}^g(t,s).$$
We then naturally expect to be able to write the solution to the second corrector equation as $u^g(t,s,z)=-v^g(t,s,z)\tilde u^g(t,s)$, where, after easy calculations using the PDE \reff{firstcor}, $\tilde u^g$ must satisfy
\begin{equation}\label{PDEutilde}
\begin{cases}
\displaystyle -\tilde u^g_t-r{\bf 1}_d\cdot \D_s\tilde u^g_s-\frac12\Tr{\sigma\sigma^T\D_{ss}\tilde u^g}+\kappa v_1(t)\tilde u^g=\overline{a}^g(t,s)\\[0.5em]
 \tilde u^g(T,\cdot)=0.
\end{cases}
\end{equation}
By the classical Feynman-Kac formula, we deduce immediately that 
$$\tilde u^g(t,s)=\mathbb E^\Q\left[\int_t^Te^{-\kappa\int_t^uv_1(w)dw}\overline a^g(u,S_u^{t,s})du\right].$$
Finally, the expansion \reff{asymp_prix} takes the form
$$p^{\eps,g}(t,s)=V^g(t,s)+\frac{\eps^2}{\gamma v_1(t)}(\tilde u^g(t,s)-\tilde u^0(t,s))+\circ(\eps^2).$$

Furthermore, when $d=1$, we can use \reff{eq:a} to deduce that 
$$\overline{a}^g(t,s)= \gamma^2v_1^2(t)\frac{\sigma^2}{2}\left(\frac34\left(\lambda^{0,1}+\lambda^{1,0}\right)\left(\frac{\mu-r}{\sigma^2}-\gamma v_1(t)s^2V_{ss}^g(t,s)\right)^2\right)^{\frac23},$$
which, when $\kappa=0$, gives us exactly the expansion proved in Corollary 3.8 of \cite{bichuch}.

\vspace{0.25em}
Let us now give sufficient conditions under which all the above calculations are rigorous and under which Assumptions \ref{assumption 1 smoothness}, \ref{assumption 4 bound on derivatives} and \ref{assumption 3 comparison} are satisfied. Concerning Assumption \ref{assumption 2 local bound}, we verify it in the next section under additional assumptions, by constructing a nearly optimal strategy (note that this approach was used by Bichuch \cite{bichuch} and Bouchard, Moreau and Soner in \cite{bms}). Moreover, we also recall that if we assume in addition that $\kappa=0$ and $d=1$, the expansion obtained in \cite{bichuch} also allows us to obtain immediately that Assumption \ref{assumption 2 local bound} is satisfied. We give a result for $d=1$ because it allows us to give precise conditions, since we now everything explicitly in this case.

\begin{Proposition}
In the framework of this section, fix $d=1$. If we assume that for $\vartheta=g\text{ or }0$

\vspace{0.25em}
{\rm(i)} There exists a constant $c_0>0$ such that
$$\abs{\frac{\mu-r}{\gamma\sigma^2v_1(t)}-s^2V_{ss}^\vartheta(t,s)}\geq c_0,\ (t,s)\in[0,T)\times(0,+\infty).$$

{\rm(ii)} $V^\vartheta$ is $C^{1,4}$ in $[0,T)\times(0,+\infty)$ and continuous on $[0,T]\times(0,+\infty)$ and there exists $\eta\in(0,1]$ such that
$$\left[\abs{V^\vartheta_s}+\abs{V^\vartheta_{ts}}+\abs{V^\vartheta_{sss}}\right](t,s)\leq \frac{C(s)}{(T-t)^{1-\eta}},\ (t,s)\in[0,T)\times(0,+\infty),$$
and there exists $\nu\in(1/4,1]$ such that
$$\abs{V^\vartheta_{ss}}(t,s)\leq \frac{C(s)}{(T-t)^{1-\nu}}, \ (t,s)\in[0,T)\times(0,+\infty),$$
for some continuous function $C$. Then Assumptions \ref{assumption 1 smoothness}, \ref{assumption 4 bound on derivatives} and \ref{assumption 3 comparison} are satisfied
\end{Proposition}

\proof
We start with Assumption \ref{assumption 1 smoothness}. First of all, it is clear that $v^\vartheta$ is $C^{1,2,2}$ in $[0,T)\times(0,+\infty)^2$ and continuous in $[0,T]\times(0,+\infty)^2$, since $V^\vartheta$ is and $v^1$ and $v^2$ are $C^\infty$ on $[0,T]$. Moreover, we have that $v^\vartheta$ is actually $C^\infty$ in $z\in(0,+\infty)$ for every $(t,s)\in[0,T]\times(0,+\infty)$. In particular, $v^\vartheta_{zz}$ is bounded on $(t,s,z)\in[0,T]\times(0,+\infty)^2$ and
$$v^\vartheta_z(t,s,z)=\gamma v^1(t)\exp\left(-\gamma v_1(t)\left(z-V^\vartheta(t,s)\right)+v_2(t)\right)>0,\ (t,s,z)\in[0,T]\times(0,+\infty)^2.$$

Next, notice that $\ybf^\vartheta$ does not depend on $z$ (see Remark \ref{rem.ygz}) and that 
$$\alpha^\vartheta(t,s)=\sigma\left(\frac{\mu-r}{\gamma\sigma^2v_1(t)}-s^2V_{ss}^\vartheta(t,s)\right),$$
so that we clearly have $(\alpha^\vartheta)^2(t,s)\geq c_1$ for some $c_1>0$. Then, the estimates on the derivatives of $\ybf^\vartheta$ are immediate consequences of the assumed estimates on the derivatives of $V^\vartheta$. Hence Assumption \ref{assumption 1 smoothness} is satisfied. Let us now look at Assumption \ref{assumption 4 bound on derivatives}. We have in this framework
$$a^\vartheta(t,s,z)=\frac{\sigma^2v^\vartheta_z(t,s,z)}{2}\gamma v^1(t)\left(\frac{3\left(\lambda^{1,0}+\lambda^{0,1}\right)}{4\gamma^2(v^1)^2(t)}\left(\frac{\mu-r}{\gamma\sigma^2v_1(t)}-s^2V_{ss}^\vartheta(t,s)\right)^2\right)^{2/3},$$
which implies that $a^\vartheta$ is continuous in $[0,T)\times(0,+\infty)^2$. Then, using \reff{eq:www}, the required estimates and regularity in $(t,s,z)$ for $w^\vartheta$ and $\mathcal O^\vartheta$ are direct consequences of the fact that $V^\vartheta$ is $C^{1,4}$ in $[0,T)\times(0,+\infty)^2$ and that $v^\vartheta$ is $C^\infty$ in $z\in(0,+\infty)$ for every $(t,s)\in[0,T]\times(0,+\infty)$. Next, $\rho^*(t,s,z)$ is a solution to a cubic equation so that it has the same regularity as its coefficients, which then implies that $\tilde w^\vartheta$ also satisfies the required regularity and estimates. Hence Assumption \ref{assumption 4 bound on derivatives} is satisfied.

\vspace{0.25em}
Finally, concerning Assumption \ref{assumption 3 comparison}, as mentioned before, obtaining a comparison theorem for viscosity solutions with polynomial growth is a classical result. Moreover, in this particular case, it is easy to check using Feynman-Kac formula that the PDE \reff{PDEutilde} has a unique smooth solution which admits the following probabilistic representation
$$u^g(t,s,z)=-Bv^g(t,s,z)\mathbb E^\Q\left[\int_t^Tv_1^2(u)e^{-\kappa\int_t^uv_1(\xi)d\xi}\left(\frac{\mu-r}{\sigma^2}-\gamma v_1(u)(S^{t,s}_u)^2V^g_{ss}(u,S^{t,s}_u)\right)^{\frac43}du\right],$$
where $B:=\frac{\gamma^2\sigma^2}{2}\left(\frac34(\lambda^{0,1}+\lambda^{1,0})\right)^{\frac23}.$ When $g=0$, this can actually be further simplified to obtain
$$u^0(t,s,z)=-Bv^0(t,s,z)\left(\frac{\mu-r}{\sigma^2}\right)^{\frac43}\int_t^Tv_1^2(u)e^{-\kappa\int_t^uv_1(\xi)d\xi}du.$$

Of course, for all this to be meaningful, the above expectation should be finite, which is once again an implicit assumption on the payoff $g$. It is easy to show that a sufficient condition for this to be true is that there exist some $\beta\in(0,3/4)$ such that
$$\abs{V_{ss}^\vartheta}(t,s)\leq \frac{C(s)}{(T-t)^{\beta}}.$$

\vspace{-2em}

\begin{flushright}$\Box$\end{flushright}
\subsubsection{A nearly optimal strategy}\label{sec:export}

In this section\footnote{The approach followed here has been suggested to us by Mete Soner and Nizar Touzi in private communications. We would like to thank them deeply.}, we derive, under additional regularity assumptions, an asymptotically optimal strategy for the problem \reff{def:ve}, which incidentally shows that Assumption \ref{assumption 2 local bound} holds. Such an approach has also been used in \cite{bichuch} and \cite{bms}.

\vspace{0.3em}
We will assume here that $\xi\longmapsto w^g(t,s,\xi)$ is $C^2$ on $\mathbb R^d$, for every $(t,s)\in[0,T]\times(0,+\infty)^d$. We can then define the following sets, for any $(t,s)\in[0,T]\times(0,+\infty)^{d}$
$$\mathcal T^g(t,s):=\left\{\xi\in\mathbb R^d,\ \frac12|\sigma^T\xi|^2\gamma^2v^2_1(t)-\frac12\Tr{\overline\alpha^g(\overline\alpha^g)^T(t,s)\overline w^g_{\rho\rho}(t,s,\gamma v_1(t)\xi)}+\overline a^g(t,s)=0\right\}.$$
We expect that, at the first order in $\eps$, the no-transaction region looks like
$${\rm NT}^{\eps,g}(t,s):=\left\{(x,y)\in\mathbb R\times\mathbb R^d,\ y-y^g(t,s)\in\eps\mathcal T^g(t,s)\right\}.$$
Indeed, according to the expansion we obtained, the region where the gradient constraint for $w^g$ is not binding is a natural candidate for being a first-order approximation of the region where the gradient constraints for $v^{\eps,g}$ are not binding either, which justifies the introduction of the set $\mathcal T^g(t,s)$. Moreover, we also remind the reader that the fast variable $\xi$ was defined as $(y-\ybf^g(t,s,z))/\eps$, which explains the introduction of the set ${\rm NT}^{\eps,g}(t,s).$ Then, according to the classical results on utility maximization with transaction costs (see for instance \cite{ss}, the optimal strategy usually consists in doing nothing while on the interior of the so-called "no-transaction" region, and making transactions when on its boundary. Thus, for any $(i,j)\in\mathcal I$, at the first order, we expect that transactions between assets $i$ and $j$ only occur on the set
$$\left(\partial{\rm NT}^{\eps,g}\right)^{i,j}(t,s):=\left\{(x,y)\in\mathbb R\times\mathbb R^d,\ y-y^g(t,s)\in\eps\left(\partial\mathcal T^g\right)^{i,j}(t,s)\right\},$$
where
$$\left(\partial\mathcal T^g\right)^{i,j}(t,s):=\left\{\xi\in\partial \mathcal T^g(t,s),\ -\lambda^{i,j}+\overline w^g_{\rho^i}(t,s,\gamma v_1(t)\xi)-\overline w^g_{\rho^j}(t,s,\gamma v_1(t)\xi)=0\right\}.$$
This naturally leads us to consider the following Markovian consumption-investment strategy, defined, for any $(t,s,x,y)\in[0,T]\times(0,+\infty)^d\times\R\times\R^d$ by
$$\nu^\eps_u:=(\cbf^g(u,S^{t,s}_u,Z^{\eps,t,s,x,y}_u),L^\eps_u),$$
where $Z^{\eps,t,s,x,y}_u:=X^{\eps,t,s,x}_u+Y_u^{\eps,t,s,y}\cdot{\bf 1}_d$, with
\begin{align}\label{skohorod}
\begin{cases}
X^{\eps,t,s,x}_u=\displaystyle x+\int_t^u\left(rX_w^{\eps,t,s,x}-\cbf^{\eps,g}(w,S_w^{t,s},Z^{\eps,t,s,x,y}_w)\right)dw+\int_t^u\mathbf{R}^0(dL_w^{\eps}),\ t\leq u\leq T,\\[0.5em]
\displaystyle Y^{\eps,t,s,y}_u=y+\int_t^u{\rm diag}[Y^{\eps,t,s,y}_w]\left(\mu dw +\sigma dW_w\right)+\int_t^u\mathbf{R}(dL_u^{\eps}),\ t\leq u\leq T\\
 \displaystyle (X^{\eps,t,s,x}_u,Y^{\eps,t,s,y}_u)\in {\rm NT}^{\eps,g}(t,s),\ u\in[t,T],\\
\displaystyle L_u^{\eps,i,j}=\int_t^u{\bf 1}_{(X^{\eps,t,s,x}_w,Y^{\eps,t,s,y}_w)\in\left(\partial{\rm NT}^{\eps,g}\right)^{i,j}(t,s)}dL^{\eps,i,j}_w, u\in[t,T],\ (i,j)\in\mathcal I.
\end{cases}
\end{align}

The system \reff{skohorod} is a classical Skorohod problem with reflection on the boundary of ${\rm NT}^{\eps,g}(t,s)$. It is well-known (see for instance \cite{lions}) that it will admit a solution if this boundary is smooth. Though, as seen in the above calculations, we know in the present framework that this assumption is satisfied when $d=1$, it is not clear at all in higher dimensions. We nonetheless do not try to address this difficult problem here and simply assume the following.
\begin{Assumption}\label{assump sko}
The first corrector equation admits a $C^2$ solution $w^g(t,s,\cdot)$ for every $(t,s)\in[0,T]\times(0,+\infty)^d$, and the Skorohod problem \eqref{skohorod} also admits a solution $(X^{\eps,t,s,x},Y^{\eps,t,s,y})$. Furthermore, $\tilde u^g$ is $C^{1,2}$ on $[0,T]\times(0,+\infty)^d$, we can take $\eta=1$ in Assumption \ref{assumption 1 smoothness}$(ii)$, $\overline a^g$ as well as all the functions of $s$ in Assumption \ref{assumption 4 bound on derivatives} have polynomial growth in $s$, uniformly in $t$.
\end{Assumption}
We emphasize again that when $d=1$, the first part of Assumption \ref{assump sko} is immediately satisfied. The regularity of $\tilde u^g$ is an implicit assumption on the regularity of $\overline a^g$, which can be checked directly when $d=1$, as all the quantities are explicit. Besides, we take $\eta=1$ in Assumption \ref{assumption 1 smoothness} for simplicity, since it prevents all the quantities appearing from exploding at $T$. The polynomial growth assumption is simply an assumption on $\overline a^g$, which can also be verified easily when $d=1$. 

\vspace{0.3em}
Let us now define
\begin{align*}
\mathcal V^{\eps,g}(t,s,x,y):&=v^g(t,s,z)-\eps^2u^g(t,s,z)-\eps^4w^g(t,s,z,\xi)\\
&=v^g(t,s,z)\left(1+\eps^2\tilde u^g(t,s)+\eps^4\overline{w}^g(t,s,\gamma v_1(t)\xi)\right).
\end{align*}
Under Assumption \ref{assump sko}, we can apply the remainder estimate of Section \ref{sect: remainder estimate} below as well as the fact that $u^g$ and $w^g$ solve the corrector equations to obtain for any $(x,y)\in{\rm NT}^{\eps,g}(t,s)$
\begin{equation}\label{eq:rem}
\mathcal L\mathcal V^{\eps,g}(t,s,x,y)-\cbf^{g}(t,s,z)V^{\eps,g}_x(t,s,x,y)=U_1(\cbf^g(t,s,z))-\eps^2\mathcal R^\eps(u^g,w^g)(t,s,z,\xi).
\end{equation}
Moreover, using Assumption \ref{assumption 1 smoothness}, and given that we can factor out $v^g$ in all the above functions, it can be readily verified that
$$\abs{\mathcal R^\eps(u^g,w^g)}(t,s,z,\xi)\leq \abs{v^g(t,s,z)}\left(\mathcal R_1+\mathcal R_2+\mathcal R_3\right), $$
with
\begin{align*}
\mathcal R_1&\leq C(\eps|\xi|+\eps^2|\xi|^2)(1+|\ybf^g|)(1+|\tilde u^g|+|\tilde u^g_s|),\\
\mathcal R_2&\leq \eps^4C(s)(1+|\ybf^g|+|\ybf^g|^2)(1+\eps|\xi|+\eps^2|\xi|^2)(1+|\xi|+\eps^{-1}),\\
\mathcal R_3&\leq C(s)\eps^4\left(1+|\xi|)(1+\abs{z-V^g}\right)+\frac{C(s)}{\abs{v^g}}\eps^4(|\tilde u^g|^2+\eps^4(1+|\xi|^2)).
\end{align*}
Notice also that we know from the solution of the second corrector equation that the set $\mathcal T^g(t,s)$ is bounded for fixed $(t,s)$. Therefore, by definition, if $(x,y)\in{\rm NT}^{\eps,g}(t,s)$, then $\xi$ is bounded, uniformly in $\eps$, by some continuous function of $s$ (remember that $t$ lives in a compact set). Therefore, using Assumption \ref{assump sko}, the above estimates can be rewritten
\begin{align*}
\mathcal R_1\leq \eps C(s),\ \mathcal R_2\leq \eps^3C(s),\ \mathcal R_3\leq \eps^4C(s)(1+\abs{z-V^g})+\frac{C(s)}{\abs{v^g}}\eps^4.
\end{align*}
Besides, it can checked immediately that there is some $\eps_0$ such that for any $\eps\leq \eps_0$, we have
\begin{equation}\label{eq:skosko}
\Lambda^\eps_{i,j}\cdot (\mathcal V^{\eps,g}_x,(\mathcal V^{\eps,g}_y)^T)^T(t,s,x,y)\geq 0, \text{ on }\left(\partial{\rm NT}^{\eps,g}\right)^{i,j}(t,s).
\end{equation}

Next, we apply It\^o's formula to $\mathcal V^{\eps,g}$ between $t$ and some stopping time $\tau_n$ which localizes the local martingales appearing. We have, using \eqref{eq:rem}, \eqref{eq:skosko} and the fact that $(X^{\eps,t,s,x},Y^{\eps,t,s,y})$ solves \eqref{skohorod},
\begin{align*}
&\mathbb E\left[\int_t^{\tau_n}U_1(\cbf^g(u,S_u^{t,s}, Z^{\eps,g,t,s,x,y}_u))du+U_1\left(\ell^\eps(X^{\eps,t,s,x}_T,Y_T^{\eps,t,s,y})-g(S_T^{t,s})\right)\right]\\
&\geq V^{\eps,g}(t,s,x,y)+\mathbb E\left[U_1\left(\ell^\eps(X^{\eps,t,s,x}_T,Y_T^{\eps,t,s,y})-g(S_T^{t,s})\right)-\mathcal V^{\eps,g}(\tau_n, S^{t,s}_{\tau_n},Z^{\eps, t,s,x,y}_{\tau_n})\right]\\
&\hspace{0.9em}-\eps^2\mathbb E\left[\int_t^{\tau_n}\abs{\mathcal R^\eps(u^g,w^g)}(u,S_u^{t,s},Z^{\eps,t,s,x,y}_u,\xi)du\right]\\
&\geq V^{\eps,g}(t,s,x,y)+\mathbb E\left[U_1\left(\ell^\eps(X^{\eps,t,s,x}_T,Y_T^{\eps,t,s,y})-g(S_T^{t,s})\right)-\mathcal V^{\eps,g}(\tau_n, S^{t,s}_{\tau_n},Z^{\eps, t,s,x,y}_{\tau_n})\right]\\
&\hspace{0.9em}-\eps^2\mathbb E\left[\int_t^{T}\abs{\mathcal R^\eps(u^g,w^g)}(u,S_u^{t,s},Z^{\eps,t,s,x,y}_u,\xi)du\right].
\end{align*}
Then, we use the fact that $U_1$ and $\mathcal V^{\eps, g}$ are non-positive to take the $\overline{\lim}$ on the left-hand side, the $\underline{\lim}$ on the right-hand side, so that, applying Fatou's lemma, we obtain
\begin{align*}
&\mathbb E\left[\int_t^{T}U_1(\cbf^g(u,S_u^{t,s}, Z^{\eps,g,t,s,x,y}_u))du+U_1\left(\ell^\eps(X^{\eps,t,s,x}_T,Y_T^{\eps,t,s,y})-g(S_T^{t,s})\right)\right]\\
&\geq V^{\eps,g}(t,s,x,y)+\mathbb E\left[U_1\left(\ell^\eps(X^{\eps,t,s,x}_T,Y_T^{\eps,t,s,y})-g(S_T^{t,s})\right)-\mathcal V^{\eps,g}(T, S^{t,s}_{T},Z^{\eps, t,s,x,y}_{T})\right]\\
&\hspace{0.9em}-\eps^2\mathbb E\left[\int_t^{T}\abs{\mathcal R^\eps(u^g,w^g)}(u,S_u^{t,s},Z^{\eps,t,s,x,y}_u,\xi)du\right]\\
&\geq \mathbb E\left[e^{-\gamma(Z^{\eps,t,s,x,y}_T-g(S_T^{t,s}))}\left(1-e^{\gamma(Z_T^{\eps,t,s,x,y}-\ell^\eps(X_T^{\eps,t,s,x},Y_T^{\eps,t,s,y}))}+\eps^4\overline{w}^g(T,S_T^{t,s},\gamma v_1(T)\xi)\right)\right]\\
&\hspace{0.9em}-\eps^2\mathbb E\left[\int_t^{T}\abs{\mathcal R^\eps(u^g,w^g)}(u,S_u^{t,s},Z^{\eps,t,s,x,y}_u,\xi)du\right]+V^{\eps,g}(t,s,x,y)\\
&\geq V^{\eps,g}(t,s,x,y)+\mathbb E\left[e^{-\gamma(Z^{\eps,t,s,x,y}_T-g(S_T^{t,s}))}\left(-C\eps^3\abs{Y_T^{\eps,t,s,y}}e^{C\eps^3\abs{Y_T^{\eps,t,s,y}}}-\eps^4(1+C(S_T^{t,s}))\right)\right]\\
&\hspace{0.9em}-\eps^2\mathbb E\left[\int_t^{T}\abs{\mathcal R^\eps(u^g,w^g)}(u,S_u^{t,s},Z^{\eps,t,s,x,y}_u,\xi)du\right].
\end{align*}
Arguing exactly as in the proof of Lemma 3.2 in \cite{bichuch}, it can be verified, using in particular the polynomial growth in $s$ of the functions appearing in the estimates, that 
$$\E\left[\abs{U_1(Z^{\eps,t,s,x,y}_T-g(S_T^{t,s}))}+\abs{v^g(u,S_u^{t,s},Z^{\eps,t,s,x,y}_u)}+\abs{Z^{\eps,t,s,x,y}_u-V^g(u,S_u^{t,s})}\right]<+\infty.$$
We therefore have
\begin{align*}
&v^{\eps,g}(t,s,z)\geq\mathbb E\left[\int_t^{T}U_1(\cbf^g(u,S_u^{t,s}, Z^{\eps,g,t,s,x,y}_u))du+U_1\left(\ell^\eps(X^{\eps,t,s,x}_T,Y_T^{\eps,t,s,y})-g(S_T^{t,s})\right)\right]\\
&\geq V^{\eps,g}(t,s,x,y)-C\eps^3=v^g(t,s,z)-\eps^2u^g(t,s,z)-C\eps^3-\eps^4\overline w^g(t,s,\gamma v_1(t)\xi),
\end{align*}
for some $C>0$. Hence, the strategy we have exhibited coincides with the value function $v^{\eps,g}$ up to the order $\eps^2$, and is "nearly" optimal. Furthermore, the above inequality implies immediately that Assumption \ref{assumption 2 local bound} is satisfied.

\subsubsection{Discussion on the Assumptions in this setting}
As we have seen above, the fact that the diffusion coefficient $\alpha^g$ should not be equal to $0$ translates directly in our setting into
$$\frac{\mu-r}{\gamma\sigma^2}v_1^{-1}(t)-s^2V_{ss}^g(t,s)  \neq 0.$$

This is an implicit assumption on the payoff $g$, which may not be satisfied if $s^2V^g_{ss}$ can become arbitrarily big, which would be the case for a Call option for instance (for which $g(x)=(x-K)^+$), since this quantity explodes to $+\infty$ as $t$ goes to $T$, when we are at the money forward ({\it i.e.} $s=Ke^{-r(T-t)}$). This condition also naturally appears in the recent work of Bouchard, Moreau and Soner \cite{bms}, and under a stronger form in \cite{bichuch} (see Assumption 3.2). However, we would like to insist on the fact that in our approach, we do not need to assume regularity on the payoff $g$ directly (except continuity) but on its Black-Scholes price which is much more regular in general. Hence, our assumptions are less restrictive than the ones in \cite{bms} and \cite{bichuch}. For instance, it can be checked directly that for any $\eps>0$, the payoffs of power-calls $g(x):=((x-K)^+)^{1+\eps}$ are covered by our results. Since those payoffs are $C^1$ and not $C^2$, they could not be treated using the previous literature (which required $C^4$ regularity of $g$). We would also like to point out that the quantity of interest here is then $s^2V_{ss}^g(t,s)$, which is the so-called {\it activity rate} of {\it portfolio Gamma} which plays a central role in the formal asymptotics obtained by Kallsen and Muhle-Karke in \cite{km1,km2}.  

\vspace{0.25em}
Notice also that a Call option does not satisfy the assumption that the third order derivative of its Black-Scholes price does not explode at time $T$ at a speed strictly less than $(T-t)^{-1}$, however, we believe that this condition can be improved by maybe using other test functions in our proof of the sub solution property at the boundary in Section \ref{sec.boundary}. This, as pointed out in \cite{bichuch}, leads to conjecture that the expansion should also hold in the case of Call options. We leave this problem for future research. However, if one considers a Digital option $g(s)={\bf 1}_{s\geq K}$, then one can readily check that the function $\tilde u^g$ becomes infinite, which shows that the expansion cannot hold in this case, and that the corresponding first order term, (if it exists) goes to $0$ more slowly than $\eps^2$. We emphasize that the exact same phenomenon was already highlighted by Possama\"{i}, Soner and Touzi \cite{pst2} in a market where the frictions came from the absence of infinite liquidity. Moreover, the techniques of proof used in this paper to show the expansion for Call option can certainly be adapted in our setting.

\vspace{-1em}

\section{Proof of Theorem \ref{th convergence de ueps}}

\vspace{-0.6em}

We would like to point out immediately to the reader that several of the proofs below (especially the proofs of the viscosity sub and super-solution properties inside the domain) are very close to the ones given in \cite{pst}. Nonetheless, they also provide some corrections to small gaps that we identified in \cite{pst}, and are made under assumptions which are a little bit more general (in particular, we no longer require the upper bound for $\ybf^g$ in their Assumption $3.1$) and we therefore think that they can be of interest. However, the proof of the viscosity sub-solution property at the boundary is new, and the derivation of the remainder estimate has to be done with a lot more precision than in their case, because of the possible explosions at the boundary.  
\label{sect: proof of th convergence}

\vspace{-1em}

\subsection{First properties and derivatives estimates}
Denote by $ L$ the upper bound of the set $ C$, we define, $ \bar{\lambda}:=\max_{(i,j)\in\mathcal I} \lambda^{i,j}, \ \ \underline{\lambda}:= \min_{(i,j)\in\mathcal I} \lambda^{i,j}.$ We would like to mention that for notational simplicity, we state all the results of this section for $u^{\eps,g}$ and $v^{\eps,g}$, but they of course still hold true for $u^{\eps,0}$ and $v^{\eps,0}$. That being said, we have first the following easy result, whose proof can be found in \cite{pst} for instance
\begin{Lemma}\label{lem: ueps estimate}
Let $ (t,s,x,y)\in [0,T) \times  (0,\infty)^d\times K_\eps $. Then 
$$ u^{\eps,g}(t,s,x,y) \geq -\eps L v_z^g(t,s,z) \left| y-\ybf^g(t,s,z)\right|,$$
so that under Assumption \ref{assumption 2 local bound} we obtain that, $ 0 \le u_*^g(t,s,x,y) \le u^{g*}(t,s,x,y) <\infty.$
\end{Lemma}

We start with a technical lemma, which will be used in the proof of Lemma \ref{Lem: u*g dep of (t,s,z)}. The proof follows exactly the same arguments as the ones given in \cite{pst}, with some modifications due to the fact that, unlike in \cite{pst}, we do not assume any upper bound for $\ybf^g_z$. We therefore provide them for the sake of completeness.

\begin{Lemma}\label{Lem estimate of hat v}
Under assumption \ref{assumption 1 smoothness}, \ref{assumption 4 bound on derivatives} and \ref{assumption 2 local bound}, the gradient of  $\hat{v}^{\epsilon,g}$ exists almost everywhere and there exists a universal constant $ A$ such that for all $(t,s,z,\xi)\in[0,T)\times(0,+\infty)^{d+1}\times\R^d$, we can find some $\eps^*:=\eps^*(t,s,z)>0$ such that
$$ |\hat{v}^{\epsilon,g}_\xi |(t,s,z,\xi) \le A \epsilon^4 | \hat{v}^{\epsilon,g}|(t,s,z,\xi),\text{ for $\eps\leq \eps^*$,}\text{ and }\hat{v}^{\epsilon,g}_z(t,s,z,\xi) \le \gamma^\eps (t,s,z,\xi), \text{ $\forall\eps>0$},$$
where
\begin{align*}
\gamma^\eps (t,s,z,\xi) :=&\ D(t,s,z)v^g_z(t,s,z-\epsilon)+ \epsilon \abs{1-\ybf^g_z(t,s,z)\cdot{\bf 1}_d}u^{\epsilon,g}(t,s,x-\epsilon,y)\\
&+\eps\sum_{i=1}^d\ybf_z^{g,i}(t,s,z)u^{\epsilon,g}(t,s,x,y-\epsilon e_i)+\epsilon^3 C(t,s,z-\eps) \left( 1+|\xi|\right)D(t,s,z)\\
&+\epsilon^3 C(t,s,z-\eps) \abs{1-\ybf^g_z(t,s,z)\cdot{\bf 1}_d}\frac{| \ybf^g(t,s,z)-\ybf^g(t,s,z-\epsilon)|}{\epsilon}\\
& +\epsilon^3 C(t,s,z-\eps) \sum_{i=1}^d\ybf^{g,i}_z(t,s,z)\frac{ | \ybf^g(t,s,z)-\ybf^g(t,s,z-\eps)-\eps e_i |}{\eps},
\end{align*}
where $C$ is the function appearing in Assumption \ref{assumption 4 bound on derivatives} and where
$$D(t,s,z):=\abs{1-\ybf^g_z(t,s,z)\cdot{\bf 1}_d}+\ybf^g_z(t,s,z)\cdot{\bf 1}_d.$$
\end{Lemma}

\textbf{Proof}
\textit{Step 1: first estimate.}  By Theorem \ref{th: viscosity eq veps}, we have for all $ 1 \le i \le d$ in the viscosity sense that
\begin{align*}
\Lambda^\epsilon_{i,0} \cdot (v^{\epsilon,g}_x,v^{\epsilon,g}_y)\geq 0 \ \text{and} \ \Lambda^\epsilon_{0,i} \cdot (v^{\epsilon,g}_x, v^{\epsilon,g}_y) \geq 0.
\end{align*}
We deduce immediately from the definition of $ \hat{v}^{\epsilon,g}$ that for all $ 1 \le i \le d$
\begin{align}\label{eq interm 1}
\frac{\epsilon^4 \lambda^{i,0}}{1+\epsilon^3 \lambda^{i,0}} \hat{v}^{\epsilon,g}_z(t,s,z,\xi)-\frac{\epsilon^3\lambda^{i,0}}{1+\epsilon^3 \lambda^{i,0}} \ybf^g_z(t,s,z) \cdot \hat{v}^{\epsilon,g}_\xi(t,s,z,\xi)+ \hat{v}^{\epsilon,g}_{\xi^i}(t,s,z,\xi) \geq 0,\\
\label{eq:2}
 \eps^4 \lambda^{0,i} \hat{v}^{\epsilon,g}_z(t,s,z,\xi)-\eps^3 \lambda^{0,i} \ybf^g_z(t,s,z) \cdot \hat{v}^{\epsilon,g}_\xi(t,s,z,\xi)-\hat{v}^{\epsilon,g}_{\xi^i}(t,s,z,\xi) \geq 0.
 \end{align}
Now since we have by Assumption \ref{assumption 1 smoothness} that for all $ 1\le i \le d$, $ \ybf^{g,i}_z(t,s,z)>0$, we have, by multiplying \eqref{eq interm 1} by $\ybf_z^{g,i}$and summing for all $ 1\le i \le d$ that in the viscosity sense
\begin{align}\label{eq interm 2}
\left(1-\eps^3 \sum_{i=1}^d\frac{\ybf^{g,i}_z(t,s,z)\lambda^{i,0}}
{1+\eps^3\lambda^{i,0}}\right)\ybf^g_z(t,s,z).\hat{v}^{\eps,g}_\xi(t,s,z,\xi)
\geq -\eps^4\sum_{i=1}^d
\frac{\lambda^{i,0}\ybf^{g,i}_z(t,s,z)}{1+\eps^3\lambda^{i,0}}\hat{v}^{\eps,g}_z(t,s,z,\xi).
\end{align}
Now, we know that there exists a $\eps^*(t,s,z)$ such that
\begin{equation}\label{epsstar}
1-\eps^3 \sum_{i=1}^d\frac{\ybf^{g,i}_z(t,s,z)\lambda^{i,0}}{1+\eps^3\lambda^{i,0}}\geq 0,\text{ for $\eps\leq\eps^*(t,s,z)$},
\end{equation}
so that in the viscosity sense, we have for $\eps\leq\eps^*(t,s,z)$
\begin{equation}\label{eq:3}
\ybf^g_z(t,s,z).\hat{v}^{\eps,g}_\xi(t,s,z,\xi)\geq 
- \frac{\sum_{i=1}^d\frac{\lambda^{i,0}\ybf^{g,i}_z(t,s,z)}{1+\eps^3\lambda^{i,0}}}{1-\eps^3
 \sum_{i=1}^d\frac{\ybf^{g,i}_z(t,s,z)\lambda^{i,0}}{1+\eps^3\lambda^{i,0}}}\eps^4\hat{v}^{\eps,g}_z(t,s,z,\xi).
\end{equation}
Using this estimate in \reff{eq:2}, we deduce
\begin{align*}
\hat v^{\eps,g}_{\xi^i}(t,s,z,\xi)&\leq \lambda^{0,i}\eps^4
\left(1+\frac{\eps^3\sum_{i=1}^d\frac{\lambda^{i,0}\ybf^{g,i}_z(t,s,z)}{1+\eps^3\lambda^{i,0}}}{1-\eps^3 
\sum_{i=1}^d\frac{\ybf^{g,i}_z(t,s,z)\lambda^{i,0}}{1+\eps^3\lambda^{i,0}}}\right)\hat v^{\eps,g}_z(t,s,z,\xi)
\\
&\leq \eps^4\overline\lambda\left(1+\frac{1-c_0}{c_0}\right)\hat v^{\eps,g}_z(t,s,z,\xi)=\eps^4\frac{\overline\lambda}{c_0}\hat v^{\eps,g}_z(t,s,z,\xi), \text{ for $\eps\leq\eps^*(t,s,z)$},
\end{align*}
where we used Assumption \ref{assumption 1 smoothness} and the fact that the map $x\longmapsto x/(1-x)$ is non-decreasing. Similarly, using \reff{eq:3} in \reff{eq interm 1} leads to
\begin{align*}
\hat v^{\eps,g}_{\xi^i}(t,s,z,\xi)&\geq -\frac{\eps^4\lambda^{i,0}}{1+\eps^3\lambda^{i,0}}\left(1+\frac{\eps^3\sum_{i=1}^d\frac{\lambda^{i,0}\ybf^{g,i}_z(t,s,z)}{1+\eps^3\lambda^{i,0}}}{1-\eps^3 
\sum_{i=1}^d\frac{\ybf^{g,i}_z(t,s,z)\lambda^{i,0}}{1+\eps^3\lambda^{i,0}}}\right)\hat v^{\eps,g}_z(t,s,z,\xi)\\
&\geq -\eps^4\frac{\overline\lambda}{c_0}\hat v^{\eps,g}_z(t,s,z,\xi), \text{ for $\eps\leq\eps^*(t,s,z)$}.
\end{align*}
Now since by the concavity of $ v^{\epsilon,g}$ in $ (x,y)$, we know that its gradient exists almost everywhere and since by Assumption \ref{assumption 1 smoothness}, $ \ybf^g$ is twice continuously differentiable, we have that $ \hat{v}^{\epsilon,g}_z$ exists almost everywhere and we have for $\eps\leq\eps^*(t,s,z)$, $|\hat{v}^{\eps,g}_\xi|\leq A\eps^4 \hat{v}^{\eps,g}_z, \ \text{where} \ A
 :=\overline\lambda/c_0.$
 
 \vspace{0.3em}
\textit{Step 2: second estimate.} We now estimate $ \hat{v}^{\eps,g}_z$. We first notice that, remembering that $v^{\eps,g}$ is clearly non-decreasing with respect to $x$ and to $y^i$ for $i=1,\ldots,d$
\begin{align}\label{eq interm 4}
\nonumber\hat{v}^{\eps,g}_z(t,s,z,\xi) &=(1-\ybf^g_z(t,s,z)\cdot{\bf 1}_d)v^{\eps,g}_x(t,s,x,y)+\ybf_z^g(t,s,z)\cdot v^{\eps,g}_y(t,s,x,y)\\
&\leq \abs{1-\ybf^g_z(t,s,z)\cdot{\bf 1}_d}v^{\eps,g}_x(t,s,x,y)+\sum_{i=1}^d\ybf_z^{g,i}(t,s,z)v^{\eps,g}_{y^i}(t,s,x,y).
\end{align}

Then by concavity of $ v^{\eps,g}$ in $ x$ and of $ v^{g}$ in $ z$ and since $ v^{\eps,g} \le v^{g}$, we have:
\begin{align*}
v^{\eps,g}_x(t,s,x,y) &\le  v^{g}_z(t,s,z-\eps)+ \frac{ v^{g}(t,s,z-\eps)-v^{\eps,g}(t,s,x-\eps,y)}{\eps}.
\end{align*}
Then by definition of $ u^{\eps,g}$, we have:
$$  v^{\eps,g}_x(t,s,x,y) \le v^{g}_z(t,s,z-\eps)+\eps \left( u^{\eps,g}(t,s,x-\eps,y)+\eps^2 w^g(t,s,z-\eps,\xi_\eps) \right), $$
where $$ \xi_\eps:= \frac{y-\ybf^g(t,s,z-\eps}{\eps}=\xi+\frac{\ybf^g(t,s,z)-\ybf^g(t,s,z-\eps) }{ \eps}.$$ 
Then we recall from the estimate of $ w^g$ given by Assumption \ref{assumption 4 bound on derivatives} that:
\begin{align*}
 | w^g(t,s,z-\eps,\xi_\eps)| &\le C(t,s,z-\eps) (1+ |\xi_\eps|) \\
 &\le C(t,s,z) \left(1+ | \xi|+ \eps^{-1} | \ybf^g(t,s,z)-\ybf^g(t,s,z-\eps)|  \right), 
 \end{align*}
for some continuous positive function $C$. Hence, we deduce
\begin{align*}
v^{\eps,g}_x(t,s,x,y) \le &\  v^g_z(t,s,z-\eps) + \eps u^{\eps,g}(t,s,x-\eps,y)\\
&+\eps^3 C(t,s,z-\eps) \left(   1+ | \xi|+ \frac{ | \ybf^g(t,s,z)-\ybf^g(t,s,z-\eps)|}{\eps}  \right).
\end{align*}
Now following the same arguments, we also have for all $ 1 \le i \le d$:
\begin{align*}
v^{\eps,g}_{y^i}(t,s,x,y) \le & \ v^g_z(t,s,z-\eps) + \eps u^{\eps,g}(t,s,x,y-\eps e_i)\\
&+\eps^3 C(t,s,z-\eps) \left(   1+ | \xi|+ \frac{ | \ybf^g(t,s,z)-\ybf^g(t,s,z-\eps)-\eps e_i |}{\eps}  \right).
\end{align*}
Plugging the estimates for $ v^{\eps,g}_x$ and $ v^{\eps,g}_{y^i}$ in \eqref{eq interm 4}, we obtain $ \hat{v}^{\eps,g}_z(t,s,z,\xi) \le \gamma^\eps(t,s,z,\xi).$

\vspace{-2em}

\begin{flushright}$\Box$\end{flushright}

\begin{Lemma}\label{Lem: u*g dep of (t,s,z)}
Under assumption \ref{assumption 1 smoothness}, \ref{assumption 4 bound on derivatives} and \ref{assumption 2 local bound}, $ u^{*,g}$ and $ u_*^g$ are only functions of $ (t,s,z)$. Furthermore, we have:
\b*
 u_*^{g}(t,s,z)
 &=&
 \underset{(\eps,t',s',z')\longrightarrow (0,t,s,z)}{\underline{\lim}}
 \bar{u}^{\eps,g}\big(t',s',z'-\ybf^g(t',s',z')\cdot \1_d,\ybf^g(t',s',z')\big)
 \\
 u^{*,g}(t,s,z)
 &=&
 \underset{(\eps,t',s',z')\longrightarrow (0,t,s,z)}{\overline{\lim}}
 \bar{u}^{\eps,g}\big(t',s',z'-\ybf^g(t',s',z')\cdot \1_d,\ybf^g(t',s',z')\big).
 \e*
\end{Lemma}

\textbf{Proof} We split the proof in two parts:

\vspace{0.25em}
\textit{Step 1.} We first show the lemma for $ t\in [0,T)$. 
The result is a consequence of the gradient constraints in \eqref{def:v0} thanks to which we obtained the estimates of Lemma \ref{Lem estimate of hat v}. By definition of $ \hat{u}^{\epsilon,g}$, we have that for every $(t,s,z,\xi)$, there exists $\eps^*(t,s,z)$ such that for any $\eps\leq\eps^*(t,s,z)$
\begin{align*}
|\hat{u}^{\epsilon,g}_\xi|(t,s,z,\xi) &\le \epsilon^{-2} | \hat{v}^{\epsilon,g}_\xi(t,s,z,\xi) |+\epsilon^2 | w^g_\xi (t,s,z,\xi) |  \le  \epsilon^2  \left( A\gamma^\eps(t,s,z,\xi) +C(t,s,z) \right),
\end{align*}
where the second inequality (and the constant $ A$) comes from Lemma \ref{Lem estimate of hat v}. Then for any $ \xi_0 \in\R^d$ such that $ 1-\sum_{i=1}^d \xi^i_0=0$, we have for $\eps\leq\eps^*(t,s,z)$
$$ \abs{\left(\sum_{i=1}^d \xi^i_0 e_i-e_0\right)\cdot (u^{\eps,g}_x,u^{\eps,g}_y)}=\frac{1}{\eps} \abs{\xi_0 \cdot \hat{u}^{\eps,g}_\xi} \le \epsilon |\xi_0| \left( A\gamma^\eps(t,s,z,\xi)  +C(t,s,z) \right).$$
Next, we remind the reader that $ u^{\eps,g}$ is locally bounded.  Fix therefore some $(t_0,s_0,x_0,y_0)$, a $r_0>0$ small such that $u^{\eps,g}$, and the continuous functions $\gamma^\eps$ and $C$ are bounded uniformly on $B_{r_0}(t_0,s_0,x_0,y_0)$. Now recall also that $\eps^*(t_0,s_0,z_0)$ is defined (see \reff{epsstar}) such that
$$1-\eps^3 \sum_{i=1}^d\frac{\ybf^{g,i}_z(t_0,s_0,z_0)\lambda^{i,0}}{1+\eps^3\lambda^{i,0}}\geq 0, \text{ for $\eps\leq\eps^*(t_0,s_0,z_0)$}.$$
However, since the left-hand side above goes to $1$ as $\eps$ goes to $0$ and since it is continuous in $(t,s,z)$, then reducing $\eps$ if necessary, this inequality will also hold for any $(t,s,x,y)\in B_{r_0}(t_0,s_0,x_0,y_0)$. Therefore, we can find a constant $ K$ independent of $ \eps$ and large enough such that for all $ \xi_0\in\R^d$ such that $ 1-\sum_{i=1}^d \xi^i_0=0$, the maps
$$ t \longmapsto u^{\eps,g}(t,s,x-t,y+t\xi^0)+\eps K t \ \text{and} \  t \longmapsto -u^{\eps,g}(t,s,x-t,y+t\xi^0)+\eps K t ,$$
are non-decreasing. Then by definition, we obtain that  $ u^{*,g}$ and $ u^g_*$ are independent of the $ \xi$-variable.

\vspace{0.25em}
\textit{Step 2.} The previous proof does not hold at $ t=T$ because the gradient constraints verified by $w^g$ may not hold at $T$, since $w^g$ may not be defined there. By definition of the relaxed semi limit, we have, $ u^g_*(T,s_0,x_0,y_0)=l_1(s_0,x_0,y_0) \wedge l_2(s_0,x_0,y_0)$,
where
\begin{align*}l_1(s_0,x_0,y_0)&:= \liminf_{(\epsilon,s,x,y)\longrightarrow (0,s_0,x_0,y_0)}\bar{u}^{\epsilon,g}(T,s,x,y)\\
 l_2(s_0,x_0,y_0)&:=\liminf_{(\eps,t,s,x,y)\longrightarrow (0,T,s_0,x_0,y_0), t\neq T} \bar{u}^{\eps,g}(t,s,x,y).
\end{align*}
We consider separately these two terms. Freezing the variable $ t=T$, we obtain that
\begin{align*}
l_1(s_0,x_0,y_0) &= \liminf_{(\epsilon,s,x,y)\longrightarrow (0,s_0,x_0,y_0)} \bar{u}^{\eps,g}(t,s,x,y)\\
&=\lim_{(\epsilon,s,x,y)\longrightarrow (0,s_0,x_0,y_0)}  \frac{U_2(z-g(s))-U_2\left(\ell^\eps(x,y)-g(s)\right)}{\eps^2}=0.
\end{align*}

Then by Step 1, we know that 
\begin{align*}
 l_2(s_0,x_0,y_0) =\liminf_{(\eps,t,s,x,y)\longrightarrow (0,T,s_0,x_0,y_0), t\neq T} \bar{u}^{\eps,g}(t,s,x,y)&=\liminf_{(t,s,x,y)\longrightarrow (T,s_0,x_0,y_0), t\neq T} u^{g}_*(t,s,x,y)\\
 &=\liminf_{(t,s,z)\longrightarrow (T,s_0,z_0), t\neq T} u^{g}_*(t,s,z),
\end{align*}
so that we obtain the required result for $ u^g_*$. The same arguments lead to the result for $ u^*(T,s,x,y)$.

\vspace{-1em}

\begin{flushright}$\Box$\end{flushright}

\subsection{The remainder estimate}
\label{sect: remainder estimate}

We now isolate an important estimate introduced in \cite{st} and \cite{pst}, which will be of crucial importance in the proofs of sub and super-solutions properties below. Following the seminal work of Evans \cite{evans} on the perturbed test function technique, it will be convenient for us to consider, for a test function $\phi$ of the second corrector equation \reff{eq:corr2}, potential test functions $ \psi$ for \reff{e.dpp} of the form
$$v^g(t,s,z)-\eps^2 \tilde{\phi}^\eps(t,s,z)-\eps^4\varpi(t,s,z,\xi),$$
where $ \tilde{\phi}^\eps$ will be a perturbation of $ \phi$, and $ \varpi $ a smooth function close to $w^g$. The aim of the following Lemma is to provide a detailed estimate of the remainder terms in the expansion of the parabolic part of \reff{e.dpp} when applied to such a function, which was formally obtained in Section \ref{sect: formal asymptotics for the value function}. We emphasize here that unlike in \cite{pst}, we want to have a very precise estimate, in particular when it comes to the derivatives of $\ybf^g$ which appear. Indeed, as mentioned in Remark \ref{rem.dig}, these derivatives may explode at time $T$, which will cause some difficulties when proving viscosity solution properties at the terminal time in the subsequent sections. Such a problem was not present in \cite{pst} which considered only the infinite horizon case.

\begin{Lemma}\label{lem remainder estimate}
Let $ \Psi^\eps(t,s,x,y):=v^g(t,s,z)-\eps^2 \phi(t,s,z)-\eps^4\hat\varpi(t,s,z,\xi)$, with smooth $\phi$ and such that $\varpi$ satisfies the same estimates as $ w^g$ in Assumption \ref{assumption 4 bound on derivatives}. We then have
\begin{align*}
\mathcal I(\Psi^\eps)(t,s,x,y):=&\ \left(k(t,s) \Psi^\eps-\mathcal L\Psi^\eps-\tilde U_1(\Psi^\eps_x)\right)(t,s,x,y)\\
=&\ \eps^2\Big[-\frac12\abs{\sigma^T(t,s)\xi}^2v^g_{zz}(t,s,z)+\frac12\Tr{\alpha^g(\alpha^g)^T(t,s,z)
\hat\varpi_{\xi\xi}(t,s,z,\xi)} \\
&\hspace{1.7em} -\mathcal A^g\phi(t,s,z)+\mathcal R^\eps(\phi,\hat\varpi)(t,s,z,\xi)\Big],
\end{align*} 
where $ \mathcal R^\eps(\phi,\hat{\varpi}):=\mathcal R^\eps(\phi,\hat\varpi)(t,s,z,\xi)$ verifies
\begin{align*}
&\abs{ \mathcal R^\eps(\phi,\hat{\varpi})} \le\  \Big[K\left(\eps|\xi|+\eps^2|\xi|^2\right)\left(1+|\ybf^g|+|\ybf^g|^2\right)\left(1+|\ybf^g_{t}|+|\ybf^g_{s}|+|\ybf^g_{z}|+|\ybf^g_{zz}|+|\ybf^g_{sz}|+|\ybf^g_{ss}|\right)\\
 &\hspace{0.9em}\times\left(1+|\phi_z|+|\phi_{zz}|+|\phi_{sz}|+\eps^4\mathfrak R^\eps(\varpi)\right)+\eps^4K\left(1+\zeta^\eps(t,s,z,\xi)\right)\left(1+|\ybf^g_z|+|\ybf^g_z|^2\right)\\
&\hspace{0.9em}\times\left(1+|\hat{\varpi}_z|+|\phi_z|^2+\eps^4|\hat\varpi_z|^2+\eps^2|\hat\varpi_{\xi}|^2+\eps^{-1}|\hat\varpi_{\xi}|\right)\Big](t,s,z),
\end{align*}
where $K$ is a positive continuous function which depends only on $r$, $\mu$, $\sigma$, $\tilde U_1$ and $v^g$ and where the quantities $\mathfrak R^\eps(\varpi)$ and $\zeta^\eps$ are defined in the proof.
\end{Lemma}

\vspace{0.25em}

\textbf{Proof} For notational simplicity, we will omit the dependence of the coefficients in the parameters. We have:
\begin{align*} 
\mathcal I(\Psi^\eps)(t,s,x,y)= &\ k \Psi^\eps- \mathcal{L} \Psi^\eps -\kappa\tilde{U}_1 \left( \Psi^\eps_x \right) \\
= & \ k v^g-\mathcal{L}  v^g-\kappa\tilde{U}_1\left( v^g_x \right)-\eps^2 \left( k \phi - \mathcal{L} \phi \right) - \eps^4\left( k \hat\varpi - \mathcal{L} \hat\varpi \right) + \kappa\left[ \tilde{U}_1\left( v^g_x \right)-\tilde{U}_1 \left( \Psi^\eps_x \right) \right].
\end{align*}
We now consider separately every term. We first recall from Section \ref{sect: formal asymptotics for the value function} that:
$$ k v^g-\mathcal{L}  v^g-\kappa\tilde{U}_1\left( v^g_x \right)= - \frac{1}{2} \left| \sigma^T (\ybf^g-y) \right|^2 v^g_{zz}.$$

Similarly to the previous calculations, we have
\begin{align*}
k\phi-\mathcal{L} \phi =&\ k\phi- \mathcal{L}^0 \phi - r z \phi_z - y \cdot \left[  \left( \mu -r \1_d \right) \phi_z +\sigma \sigma^T \D_{sz} \phi \right]-\frac{1}{2} | \sigma^T y|^2 \phi_{zz}\\
=&\ \mathcal{A}^g \phi-\kappa\cbf^g \phi_z + (\ybf^g-y)\cdot  \left(\left( \mu -r \1_d \right) \phi_z+ \sigma \sigma^T \D_{sz}\phi\right)-\frac{1}{2} \phi_{zz} \left( | \sigma^T y|^2-|\sigma^T \ybf^g|^2 \right).
\end{align*}
Define then
$$ \mathcal{R}_\phi:=  (\ybf^g-y)\cdot  \left( \mu -r \1_d \right) \phi_z-\frac{1}{2} \phi_{zz} \left( | \sigma^T y|^2-|\sigma^T \ybf^g|^2 \right)-\left(y-\ybf^g \right) \cdot \sigma \sigma^T \D_{sz} \phi. $$
We clearly have that
\begin{align*}
 | \mathcal{R}_\phi(t,s,z,\xi)|& \le \eps \left( | \xi| \left| \mu -r \1_d \right| |\phi_z|+\frac{\abs{\sigma}^2}{2}\left(2 |\ybf^g| | \xi|+\eps | \xi|^2\right)|\phi_{zz}|+|\sigma|^2 |\xi| | \D_{sz} \phi | \right)\\
 &\leq  K_1(t,s,z)\left(\eps|\xi|+\eps^2|\xi|^2\right)\left(1+|\ybf^g|\right)\left(1+|\phi_z|+|\phi_{zz}|+|\phi_{sz}|\right),
 \end{align*}
where $K_1$ is a positive continuous function which depends only on $\sigma$, $r$ and $\mu$. The third term is more tedious. We sum up the calculations here
$$ |\varpi_y| \le |\hat\varpi_z|+\frac{1}{\eps}\left(1+\abs{\ybf^g_z}\right) |\hat\varpi_{\xi}|, \ \ |\varpi_x|\le  |\hat{\varpi}_z|+\frac{1}{\eps}|\ybf^g_z|| \hat{\varpi}_{\xi}|,$$
$$ |\varpi_s|\le  |\hat{\varpi}_s|+\frac{1}{\eps}|\ybf^g_s|| \hat{\varpi}_{\xi}|,\ \ |\varpi_t|\leq |\hat\varpi_t|+\frac1\eps|\ybf^g_t||\hat\varpi_\xi|,$$
$$\abs{\varpi_{yy}}\leq {\rm Const}\left(|\hat\varpi_{zz}|+\frac1\eps\left(\left(1+|\ybf^g_z|\right)|\hat\varpi_{z\xi}|+|\ybf^g_{zz}||\hat\varpi_\xi|\right)\right),$$
$$\abs{\varpi_{ys}}\leq {\rm Const}\left(|\hat\varpi_{sz}|+\frac1\eps\left(\left(1+|\ybf^g_z|\right)|\hat\varpi_{s\xi}|+|\ybf^g_{sz}||\hat\varpi_\xi|+|\ybf^g_s||\hat\varpi_{z\xi}|\right)\right),$$
$$\abs{\varpi_{ss}}\leq {\rm Const}\left(|\hat\varpi_{ss}|+\frac1\eps\left(\left(1+|\ybf^g_s|\right)|\hat\varpi_{s\xi}|+|\ybf^g_{ss}||\hat\varpi_\xi|\right)\right).$$
We deduce that
$$  \text{Tr}\left[\sigma \sigma^T\left( \D_{yy}+\D_{ss}+2\D_{sy}\right)\varpi \right]=\frac{1}{\eps^2}\text{Tr}\left[ \alpha^g ( \alpha^g)^T \hat{\varpi}_{\xi \xi} \right]  +\widetilde{\mathcal{R}}_2(\varpi), $$
where
\begin{align*}
 \widetilde{\mathcal{R}}_2(\varpi) \le&\ {\rm Const}\left(|\ybf^g|^2+\eps^2|\xi^2|\right)\left(|\hat\varpi_{zz}|+\frac1\eps\left(\left(1+|\ybf^g_z|\right)|\hat\varpi_{z\xi}|+|\ybf^g_{zz}||\hat\varpi_\xi|\right)\right)\\
 &+{\rm Const}|s|\left(|\ybf^g|+\eps |\xi|\right)\left(|\hat\varpi_{sz}|+\frac1\eps\left(\left(1+|\ybf^g_z|\right)|\hat\varpi_{s\xi}|+|\ybf^g_{sz}||\hat\varpi_\xi|+|\ybf^g_s||\hat\varpi_{z\xi}|\right)\right)\\
 &+{\rm Const}|s|^2\left(|\hat\varpi_{ss}|+\frac1\eps\left(\left(1+|\ybf^g_s|\right)|\hat\varpi_{s\xi}|+|\ybf^g_{ss}||\hat\varpi_\xi|\right)\right).
\end{align*}

We therefore deduce
$$ -\eps^4\left( k \varpi -\mathcal{L} \varpi \right)(t,s,x,y)=\frac{\eps^2}{2}\text{Tr}\left[ \alpha^g ( \alpha^g)^T \hat{\varpi}_{\xi \xi} \right]  +\eps^4\mathcal{R}_2(\varpi),  $$
where
\begin{align*}
\mathcal{R}_2(\varpi) \le&\ |k||\hat\varpi|+|\hat\varpi_t|+\frac1\eps|\ybf^g_t||\hat\varpi_\xi|+|\mu||s|\left( |\hat{\varpi}_s|+\frac{1}{\eps}|\ybf^g_s|| \hat{\varpi}_{\xi}|\right)\\
&+ |r|\left(1+|\ybf^g|+\eps|\xi|\right)\left(|\hat{\varpi}_z|+\frac{1}{\eps}|\ybf^g_z|| \hat{\varpi}_{\xi}|\right)\\
&\ +|\mu|\left(|\ybf^g|+\eps|\xi|\right)\left(|\hat\varpi_z|+\frac{1}{\eps}\left(1+\abs{\ybf^g_z}\right) |\hat\varpi_{\xi}|\right)\\
&\ +{\rm Const}\left(|\ybf^g|^2+\eps^2|\xi|^2\right)\left(|\hat\varpi_{zz}|+\frac1\eps\left(\left(1+|\ybf^g_z|\right)|\hat\varpi_{z\xi}|+|\ybf^g_{zz}||\hat\varpi_\xi|\right)\right)\\
 &+{\rm Const}|s|\left(|\ybf^g|+\eps |\xi|\right)\left(|\hat\varpi_{sz}|+\frac1\eps\left(\left(1+|\ybf^g_z|\right)|\hat\varpi_{s\xi}|+|\ybf^g_{sz}||\hat\varpi_\xi|+|\ybf^g_s||\hat\varpi_{z\xi}|\right)\right)\\
 &+{\rm Const}|s|^2\left(|\hat\varpi_{ss}|+\frac1\eps\left(\left(1+|\ybf^g_s|\right)|\hat\varpi_{s\xi}|+|\ybf^g_{ss}||\hat\varpi_\xi|\right)\right)\\
 \leq &\ K_2(t,s,z)\left(1+\eps|\xi|+\eps^2|\xi|^2\right)\left(1+|\ybf^g|+|\ybf^g|^2\right)\\
 &\times\left(1+|\ybf^g_{t}|+|\ybf^g_{s}|+|\ybf^g_{z}|+|\ybf^g_{zz}|+|\ybf^g_{sz}|+|\ybf^g_{ss}|\right)\mathfrak R^\eps(\varpi),
\end{align*}
where $K_2(t,s,z)$ is a positive continuous function which depends only on $r$ and $\mu$ and where
$$\mathfrak R^\eps(\varpi):=|\hat\varpi|+|\hat\varpi_{t}|+|\hat\varpi_{s}|+|\hat\varpi_{z}|+|\hat\varpi_{zz}|+|\hat\varpi_{sz}|+|\hat\varpi_{ss}|+\eps^{-1}\left(|\hat\varpi_{\xi}|+|\hat\varpi_{z\xi}|+|\hat\varpi_{s\xi}|\right).$$

Summarizing up, we have that the remainder $ \mathcal{R}^\eps(\phi,\hat\varpi)$ denoted $ \mathcal{R}$ for short here verifies:
\begin{align*}
| \mathcal{R}^\eps(\phi,\hat\varpi)|(t,s,z,\xi) \le \Big[ |\mathcal{R}_\phi|+\eps^4|\mathcal{R}_2(\varpi)|+|\tilde{U}_1\left( v^g_x \right)-\tilde{U}_1 \left( \psi^\eps_x \right)+\eps^2 \cbf^g \phi_z|\Big](t,s,z,\xi) .
\end{align*}
We now estimate the last term above. Recall that $ \cbf^g=-\widetilde{U}_1'(v^g_z(t,s,z))$. Hence, we have, omitting the dependence in $(t,s,z,\xi)$
$$\mathcal R_{\tilde U_1}:= \tilde{U}_1 \left( \Psi^\eps_x \right)-\tilde{U}_1\left( v^g_x \right)-\eps^2 \cbf^g \phi_z=\tilde{U}_1 \left( \Psi^\eps_x \right)-\tilde{U}_1\left( v^g_x \right)+(\Psi^\eps_x-v_z^g) \widetilde{U}_1'(v^g_z)+r_1, $$
where $ r_1:=\eps^4 \varpi_x \tilde{U}_1'(v^g_z)$ verifies $ |r_1| \le\eps^4|\tilde{U}_1'(v^g_z)| \left(|\hat{\varpi}_z|+\eps^{-1}|\ybf^g_z|| \hat{\varpi}_{\xi}|\right).$ Then we have, using that $ \tilde{U}_1$ is concave
\begin{align*}
\left|\tilde{U}_1 \left( \Psi^\eps_x \right)-\tilde{U}_1\left( v^g_x \right)+(\Psi^\eps_x-v^g_z) \tilde{U}_1'(v^g_z) \right| \le \left| \Psi^\eps_x-v^g_x \right| \left| \tilde{U}_1'( \Psi^\eps_x)-\tilde{U}_1'( v^g_x)\right|.
\end{align*}
Then since $ \tilde{U}_1$ is $ C^2$ we have that
$$ \left| \tilde{U}_1'( \Psi^\eps_x)-\tilde{U}_1'( v^g_x)\right| \le \left| \Psi^\eps_x-v^g_x \right| \zeta^\eps(t,s,z,\xi), \text{ where }\zeta^\eps(t,s,z,\xi):=\underset{m\in K^\eps(t,s,z,\xi)}{\sup}\abs{\tilde U_1^{''}(m)},$$
where $K^\eps(t,s,z,\xi):=\text{Supp}(\tilde U_1)\cap\left\{m\in\mathbb R,\ \abs{m+v^g_z}\leq \mathfrak H^\eps(t,s,z,\xi)\right\}$, with $$\mathfrak H^\eps(t,s,z,\xi):=|v^g_z|+\eps^2\left(|\phi_z|+\eps^2\left(|\hat{\varpi}_z|+\eps^{-1}|\ybf^g_z|| \hat{\varpi}_{\xi}|\right)\right).$$
Then we obtain 
\begin{align*}
\abs{\mathcal R_{\tilde U_1}}\leq&\ \eps^4\abs{\tilde{U}_1'(v^g_z)} \left(|\hat{\varpi}_z|+\frac{|\ybf^g_z|| \hat{\varpi}_{\xi}|}{\eps}\right)+\eps^4{\rm Const}\left(|\phi_z|^2+\eps^4\left(|\hat{\varpi}_z|^2+\frac{|\ybf^g_z|^2| \hat{\varpi}_{\xi}|^2}{\eps^2}\right)\right)\zeta^\eps(t,s,z,\xi)\\
\leq &\ \eps^4K_3(t,s,z)\left(1+\zeta^\eps(t,s,z,\xi)\right)\left(1+|\ybf^g_z|+|\ybf^g_z|^2\right)\\
&\times\left(1+|\hat{\varpi}_z|+|\phi_z|^2+\eps^4|\hat\varpi_z|^2+\eps^2|\hat\varpi_{\xi}|^2+\eps^{-1}|\hat\varpi_{\xi}|\right),
\end{align*}
where $K_3$ is a positive continuous function which depends only on $\tilde U_1$ and $v^g$. Finally, we have
\begin{align*}
\abs{\mathcal{R}} \le&\  K(t,s,z)\left(\eps|\xi|+\eps^2|\xi|^2\right)\left(1+|\ybf^g|+|\ybf^g|^2\right)\\
 &\times\left(1+|\ybf^g_{t}|+|\ybf^g_{s}|+|\ybf^g_{z}|+|\ybf^g_{zz}|+|\ybf^g_{sz}|+|\ybf^g_{ss}|\right)\left(1+|\phi_z|+|\phi_{zz}|+|\phi_{sz}|+\eps^4\mathfrak R^\eps(\varpi)\right)\\
 &+\eps^4K(t,s,z)\left(1+\zeta^\eps(t,s,z,\xi)\right)\left(1+|\ybf^g_z|+|\ybf^g_z|^2\right)\\
&\times\left(1+|\hat{\varpi}_z|+|\phi_z|^2+\eps^4|\hat\varpi_z|^2+\eps^2|\hat\varpi_{\xi}|^2+\eps^{-1}|\hat\varpi_{\xi}|\right),
\end{align*}
where $K$ is a positive continuous function which depends only on $r$, $\mu$, $\sigma$, $\tilde U_1$ and $v^g$.

\vspace{-2em}

\begin{flushright}$\Box$\end{flushright}

\subsection{Viscosity subsolution on $ [0,T)\times \R^d \times \R_+$}
\label{sect:subsol interior}

We focus here on the interior of the domain. Consider $ (t_0,s_0,z_0)\in [0,T)\times \R^d \times \R_+$ and $ \phi\in C^2([0,T)\times \R^d \times \R_+, \R)$ such that for all $ (t,s,z)\in [0,T)\times \R^d \times \R_+ \backslash\left\{ (t_0,s_0,z_0)\right\}$:
$$ 0=(u^{*,g}-\phi)(t_0,s_0,z_0)>(u^{*,g}-\phi)(t,s,z).$$
We want to show that $ \mathcal{A}^g\phi(t_0,s_0,z_0)-a^g(t_0,s_0,z_0) \le 0$. We separate the proof in 4 steps.

\vspace{0.25em}
\textit{Step 1:} By Lemma \ref{Lem: u*g dep of (t,s,z)}, there exists a sequence $ (t^\eps,s^\eps,z^\eps)\longrightarrow(t_0,s_0,z_0)$ when $ \eps\longrightarrow 0$ such that
$$ \hat{u}^{\eps,g}(t^\eps,s^\eps,z^\eps,0) \underset{\eps \longrightarrow 0}{\longrightarrow} u^{*,g}(t_0,s_0,z_0).$$
Then we have that $ l^\eps_*:= \hat{u}^{\eps,g}(t^\eps,s^\eps,z^\eps,0)-\phi(t^\eps,s^\eps,z^\eps) \longrightarrow 0$ and $(x^\eps,y^\eps)\longrightarrow(x_0,y_0)$ where
\begin{align*}
(x^\eps,y^\eps):=\left(z^\eps-\ybf^g(t^\eps,s^\eps,z^\eps)\cdot \1_d,\ybf^g(t^\eps,s^\eps,z^\eps)\right),\\(x_0,y_0):=\left(z_0-\ybf^g(t_0,s_0,z_0)\cdot \1_d,\ybf^g(t_0,s_0,z_0)\right).
\end{align*}

Now since $ u^\eps$ is locally bounded from above, there exists $ r_0>0$ and $ \eps_0>0$ depending on $ (t_0,s_0,x_0,y_0)$ such that:
 $$ b_*:=\sup\left\{u^{\eps,g}(t,s,x,y),\ (t,s,x,y)\in B_{r_0}(t_0,s_0,x_0,y_0),\ \eps\in(0,\eps_0]\right\}<+\infty,$$
where, reducing $r_0$ if necessary, the ball is strictly included in the interior of the domain, and where, reducing $\eps_0$ if necessary, we can assume w.l.o.g. that $(t^\eps,s^\eps,x^\eps,y^\eps)\in B_{r_0}(t_0,s_0,x_0,y_0)$ for $\eps\leq \eps_0$. We now build a test function from $ \phi$ for $ v^{\eps,g}$ in order to apply the PDE associated to $ v^{\eps,g}$. We define for $ (\eps,\delta)\in (0,1]^2$ the function $ \hat{\psi}^{\eps,\delta}$ and the corresponding $ \psi^{\eps,\delta}$ by:
$$ \hat{\psi}^{\eps,\delta}(t,s,z,\xi):= v^g(t,s,z)-\eps^2 \left( l^\eps_*+\phi(t,s,z)+\hat{\Phi}^\eps(t,s,z,\xi) \right)-\eps^4 (1+\delta) \tilde w^g(t,s,z,\xi), $$
 where $ \hat{\Phi}^{\eps}$ is defined by
 $$ \hat{\Phi}^{\eps}(t,s,z,\xi):=c \left( (t-t^\eps)^4+(s-s^\eps)^4+(z-z^\eps)^4+\eps^4 (\tilde w^g)^4(t,s,z,\xi) \right),$$
and $ c$ is a constant chosen large enough so that for $\eps\leq \eps_0$
\begin{align}\label{subsol eq interm 1}
\Phi^\eps\geq 1+b_*-\phi,\text{ on }B_{r_0}(t_0,s_0,x_0,y_0)\backslash B_{r_0/2}(t_0,s_0,x_0,y_0).
\end{align}

Notice that $ c_0$ is independent of $\eps$. The constant $ \delta$ will be fixed later. We also emphasize that by assumption, $ w^g$ and $\tilde w^g$ are only $ C^1$ in $\xi$ on the whole domain.

\vspace{0.3em}
\textit{Step 2:} We now show that for $ \epsilon$ and $ \delta$ small enough, the difference $ (v^{\epsilon,g}-\psi^{\eps,\delta})$ has a local minimizer in $ B_0:=B_{r_0}(t_0,s_0,x_0,y_0)$. Indeed it is sufficient to show that $ I^{\eps,\delta}$ has a local minimizer where:
\begin{align*}
I^{\eps,\delta}(t,s,x,y):=&\ \frac{v^{\eps,g}(t,s,x,y)-\psi^{\eps,\delta}(t,s,x,y)}{\eps^2}\\
=&-u^{\eps,g}(t,s,x,y)+l_*^\eps+\phi(t,s,z)+\Phi^\eps(t,s,x,y)+\eps^2\delta \tilde w^g(t,s,z,\xi)\\
&-\eps^2w^g(t,s,z,\eta^g(t,s,z)\rho^*(t,s,z)).
\end{align*}
Now since $ w^g$ and $\rho^*(t,s,z)$ are continuous, $\tilde w^g$ is non-negative and using \eqref{subsol eq interm 1}, for $ \delta>0 $ small enough and $\eps\leq\eps_0$,  we have for any $ (t,s,x,y)\in \partial B_0$:
 $$I^{\eps,\delta}(t,s,x,y)\geq -u^{\eps,g}(t,s,x,y)+l_*^\eps+1+b_*-\eps^2w^g(t,s,z,\eta^g(t,s,z)\rho^*(t,s,z))\geq \frac{1}{2}+l_*^\eps>0,$$
for $ \eps$ small enough. Now since $ I^{\eps,\delta}(t^\eps,s^\eps,x^\eps,y^\eps) \longrightarrow 0$ when $ \eps \longrightarrow 0$, this implies that for $ \eps$ small enough, $ I^{\eps,\delta}$ has a local minimizer $ (\tilde{t}^\eps,\tilde{s}^\eps,\tilde{x}^\eps,\tilde{y}^\eps)$ in $ B_0$ and we introduce:
$$ \tilde{z}^\eps:=\tilde{x}^\eps+\tilde{y}^\eps \cdot \1_d , \ \text{and} \ \tilde{\xi}^\eps:= \frac{ \tilde{y}^\eps-\ybf^g ( \tilde{t}^\eps,\tilde{s}^\eps,\tilde{z}^\eps)}{\eps}.$$
To summarize, we have:
$$ \min_{B_0} ( \hat{v}^{\eps,\delta}-\hat{\psi}^{\eps,\delta})=( \hat{v}^{\eps,\delta}-\hat{\psi}^{\eps,\delta})(\tilde{t}^\eps,\tilde{s}^\eps,\tilde{z}^\eps,\tilde{\xi}^\eps), \ \text{with} \  | \tilde{t}^\eps-t_0|+ |\tilde{s}^\eps-s_0|+ |\tilde{z}^\eps-z_0| \le r_0, \ \  \abs{\tilde{\xi}^\eps} \le \frac{r_1}{\eps,} $$
for some constant $ r_1$. Now since $ \psi^{\eps,\delta}$ is at least $C^1$, we have that by the dynamic programming equation verified by $ v^{\eps,g}$ that:
\begin{align}\label{subsol eq interm 2}
\Lambda^\eps_{i,j}
 \cdot \big(\psi_x^{\eps,\delta},\psi_y^{\eps,\delta}\big)
 (\tilde{t}^\eps,\tilde{s}^\eps,\tilde{x}^\eps,\tilde{y}^\eps)
 \geq 0,  \text{ for }
 (i,j)\in\mathcal I.
\end{align}
 
 \textit{Step 3:} Our aim in this step is to show that for $ \eps$ small enough, $ \psi^{\eps,\delta}$ is actually $ C^{2}$ in $\xi$. Thank to Proposition \ref{prop:wg}, it is enough to show that for $ \eps$ small enough we have:
 $$ \tilde{\rho}^\eps :=
\tilde{\xi}^\eps/\eta^g(\tilde{t}^\eps,\tilde{s}^\eps,\tilde{z}^\eps) \ \in \ 
\mathcal O^g_0(\tilde{t}^\eps,\tilde{s}^\eps,\tilde{z}^\eps),
$$
where $ \mathcal O^g_0(t,s,z) $ is the open set introduced in Proposition \ref{prop:wg}. Assume on the contrary that there exists $ \eps_n\longrightarrow 0$ such that for $n$ large enough $ \tilde{\rho}^{\eps_n} \notin  \mathcal O^g_0(\tilde{t}^{\eps_n},\tilde{s}^{\eps_n},\tilde{z}^{\eps_n})$. Then since $ \bar{w}^g$ is $ C^{1}$ and thanks to \eqref{eq:corrector}, we have:
$$ -\lambda^{i_0^n,j_0^n}v_z^g(\tilde{t}^{\eps_n},\tilde{s}^{\eps_n},\tilde{z}^{\eps_n})
 +\big(\tilde{w}_{\xi^{i_0^n}}^g-\tilde{w}^g_{\xi^{j_0^n}}\big)
  (\tilde{t}^{\eps_n},\tilde{s}^{\eps_n},\tilde{z}^{\eps_n},\tilde{\xi}^{\eps_n})=0 \  \text{for some}  \
 (i_0^n,j_0^n)\in\mathcal I. $$
We obtain then by boundedness of $\big(\tilde{t}^{\eps_n},\tilde{s}^{\eps_n},\tilde{z}^{\eps_n},\eps_n\tilde{\xi}^{\eps_n}\big)_n$, \eqref{subsol eq interm 2} and using Assumption \ref{assumption 1 smoothness} (and in particular the constant $c_0$ introduced there)
\begin{align*}
& -4c_0\eps_n^2(\eps_n \tilde w^g)^3(\tilde{t}^{\eps_n},\tilde{s}^{\eps_n},\tilde{z}^{\eps_n},\tilde{\xi}^{\eps_n})
\left(\tilde w^g_{\xi^{i_0^n}}-\tilde w^g_{\xi^{j_0^n}}\right)(\tilde{t}^{\eps_n},\tilde{s}^{\eps_n},\tilde{z}^{\eps_n},\tilde{\xi}^{\eps_n}) \\
&+\eps_n^3v_z^g(\tilde{t}^{\eps_n},\tilde{s}^{\eps_n},\tilde{z}^{\eps_n})\left[\lambda^{i_0^n,j_0^n}
-(1+\delta)(\overline{w}^g_{\xi^{i_0^n}}-\overline{w}^g_{\xi^{j_0^n}})(\tilde{t}^{\eps_n},\tilde{s}^{\eps_n},\tilde{z}^{\eps_n},
\tilde{\rho}^{\eps_n})\right]+\circ(\eps_n^3)\geq 0.
\end{align*}
And by positivity of $ \tilde w^g$, we have:
 \begin{align*}
 0
 &\le
 -4c_0\lambda^{i_0^n,j_0^n}v_z^g(\tilde{t}^{\eps_n},\tilde{s}^{\eps_n},\tilde{z}^{\eps_n})\eps_n^2
  (\eps_n\tilde w^g)^3(\tilde{t}^{\eps_n},\tilde{s}^{\eps_n},\tilde{z}^{\eps_n},\tilde{\xi}^{\eps_n})
 -\delta\lambda^{i_0^n,j_0^n}\eps_n^3
  v_z^g(\tilde{t}^{\eps_n},\tilde{s}^{\eps_n},\tilde{z}^{\eps_n})
 +\circ(\eps_n^3)
 \\
 &\le 
 -\delta\lambda^{i_0^n,j_0^n}\eps_n^3
  v_z^g(\tilde{t}^{\eps_n},\tilde{s}^{\eps_n},\tilde{z}^{\eps_n})
  +\circ(\eps_n^3),
 \end{align*}
which leads to a contradiction when $n$ goes to $+\infty$.

\vspace{0.25em}
\textit{Step 4:} Since $ \psi^{\eps,\delta}$ is smooth enough, we are now able to use it as a test function for the parabolic operator in \reff{e.dpp}. By the supersolution property of $ v^{\eps,g}$, we have:
$$  kv^{\eps,g}- \Lc \psi^{\eps,\delta} - \Ut_1(\psi^{\eps,\delta}_{x})(\tilde{t}^\eps,\tilde{s}^\eps,\tilde{x}^\eps,\tilde{y}^\eps) \geq 0.  $$
Since $ (t,s,z) \mapsto \mathcal{O}^g_0(t,s,z)$ is continuous by Assumption \ref{assumption 4 bound on derivatives} and since $ (\tilde{t}^\eps,\tilde{s}^\eps,\tilde{z}^\eps)$ is bounded, we know that $ ( \tilde{\xi}^\eps)_\eps$ is bounded. By standard results of the theory of viscosity solutions, we then have a sequence $ (\eps_n)_n$ such that $ \eps_n \longrightarrow 0$ and such that $$ (t_n,s_n,z_n,\xi_n):=(\tilde{t}^{\eps_n},\tilde{s}^{\eps_n},\tilde{z}^{\eps_n},\tilde{\xi}^{\eps_n}) \longrightarrow(t_0,s_0,z_0,\tilde{\xi}),$$ for some $\tilde\xi\in\R^d$. We then have 
\begin{align*}
-\frac12 v^{g}_{zz}(t_n,s_n,z_n)\abs{\sigma^T(t_n,s_n)\xi_n}^2+\frac12(1+\delta)
\Tr{\alpha^g(\alpha^g)^T(t_n,s_n,z_n)w^g_{\xi\xi}(s_n,z_n,\xi_n)}& \\
-\mathcal A^g\phi(t_n,s_n,z_n)-\mathcal{A}^g \Phi^{\eps_n}(t_n,s_n,x_n,y_n) +\mathcal R^{\eps_n}(\phi+\Phi^{\eps_n},(1+\delta)\tilde w^g)(t_n,s_n,z_n,\xi_n) &\geq 0,
\end{align*}
where the remainder term $\mathcal R^{\eps_n}(\phi+\Phi^{\eps_n},(1+\delta)\tilde w^g)(t_n,s_n,z_n,\xi_n)$ is controlled using the result of Lemma \ref{lem remainder estimate}. We know that $ w^g$ is $ C^2$ at the points $ (t_n,s_n,z_n,\xi_n)$ but not necessarily at $ (t_0,s_0,z_0,\tilde{\xi})$, which might be so that $\tilde\rho:=\tilde\xi/(\eta^g(t_0,s_0,z_0))\in\partial\mathcal O^g_0(t,s,z)$. Now we remind the reader that by definition of $ w^g$ and since $ \rho_{n} \in \mathcal{O}^g(t_n,s_n,z_n,\xi_n) $
$$ -\frac12 v^{g}_{zz}(t_n,s_n,z_n)\abs{\sigma^T(t_n,s_n)\xi_n}^2+\frac12
\Tr{\alpha^g(\alpha^g)^Tw^g_{\xi\xi}}(t_n,s_n,z_n,\xi_n)=-a^g(t_n,s_n,z_n),$$
so that:
\begin{align*}
&a^g(t_,s_n,z_n)-\mathcal A^g\phi(t_n,s_n,z_n)-\mathcal{A}^g \Phi^{\eps_n}(t_n,s_n,x_n,y_n) +\delta \Big( a^g(t_n,s_n,z_n)\\
&+\frac12 v^{g}_{zz}(t_n,s_n,z_n)\abs{\sigma^T(t_n,s_n)\xi_n}^2 \Big)+\mathcal R^{\eps_n}(\phi+\Phi^\eps_n,(1+\delta)\tilde w^g)(t_n,s_n,z_n,\xi_n) \geq 0.
\end{align*}
Therefore, $w^g_{\xi\xi}$ no longer appears directly in the above equation, except in the remainder $\mathcal R^{\eps_n}(\phi+\Phi^\eps_n,(1+\delta)\tilde w^g)$ for which it is implicitly understood that we do the same transformation. Now by continuity of the map $ (t,s,z) \longmapsto a^g(t,s,z)$ stated in Assumption \ref{assumption 4 bound on derivatives}, and since we clearly have that $ \mathcal{R}^{\eps_n}(\phi+\phi^{\eps_n},(1+\delta)\tilde w^g)(t_n,s_n,z_n,\xi_n) \longrightarrow 0 $ (recall that we are away from $T$ here, so that none of the quantities in the upper bound given in Lemma \ref{lem remainder estimate} can explode) and $\mathcal{A}^g \Phi^{\eps_n}(t_n,s_n,x_n,y_n)\longrightarrow 0 $ when $ n \longrightarrow \infty$, $\Phi^\eps$ and all its derivatives go to $0$. Finally, we obtain
$$ a^g(t_0,s_0,z_0)-\mathcal A^g\phi(t_0,s_0,z_0)+\delta\left(a^g(t_0,s_0,z_0)-\frac12v^g_{zz}(t_0,s_0,z_0)|\sigma^T(t_0,s_0)\tilde\xi|^2\right)\geq 0. $$
Recall that $ \tilde{\xi}$ may depend on $ \delta $ but is uniformly bounded. Then we can send $ \delta $ to $ 0 $ to obtain the required result.

\vspace{-2em}

\begin{flushright}$\Box$\end{flushright}

\subsection{Viscosity subsolution on $ \big\{T\big\}\times \R^d \times \R_+$}\label{sec.boundary}
\label{sect:subsol frontier}
In contrast with the previous section, the use of $ u^{g,\epsilon}$ is not necessary here, and we will therefore concentrate only on $\bar u^{g,\eps}$.

\vspace{0.25em}
Let $(s_0,z_0,\phi)\in(0,+\infty)^{d+1}
\times C^2\left((0,+\infty)^{d+1}\right)$ be such that
$$
0=(u^{g,*}-\phi)(T,s_0,z_0)>(u^{g,*}-\phi)(t,s,z),\ \forall(t,s,z)
\in [0,T] \times (0,+\infty)^{d+1}\backslash\left\{(T,s_0,z_0)\right\}.
$$
By definition of viscosity solutions, we want to deduce that $ \phi(T,s_0,z_0) \le 0$. Assume on the contrary that $ \phi(T,s_0,z_0)>2\delta$ for some $\delta>0$. Then we have for $ r_0>0$ small enough,
\begin{align}\label{eq phi>delta}
\phi(t,s,z)>\delta, \ \forall (t,s,z)\in [T-r_0,T] \times B_{r_0}(s_0,z_0).
\end{align}

Let us then consider a sequence $ (t_\epsilon,s_\epsilon,z_\epsilon)$ converging to $(T,s_0,z_0)$ such that $\hat{ \bar{u}}^{\epsilon,g}(t_\epsilon,s_\epsilon,z_\epsilon,0) \longrightarrow u^{*,g}(T,s_0,z_0)$. We introduce:
\begin{align}\label{eq leps converges}
 l_*^\epsilon:=\hat{\bar u}^{\epsilon,g}(t_\epsilon,s_\epsilon,z_\epsilon,0)-\phi(t_\epsilon,s_\epsilon,z_\epsilon) \underset{\epsilon \longrightarrow 0}{\longrightarrow} 0.
\end{align}
By assumption \ref{assumption 2 local bound}, there exists $ 0<r_1<r_0$ such that:
$$ b^*:=\sup \big\{ \hat{\bar{u}}^{g,\epsilon}(t,s,z,0),\ (t,s,z)\in [T-r_1,T] \times B_{r_1}(s_0,z_0) \big\}<+\infty.$$
We will denote for simplicity $ B_1:=[T-r_1,T] \times B_{r_1}(s_0,z_0) $. We split the proof in two parts.

\vspace{0.25em}
\textit{Step 1:} We first show that there is some $\eps_0$ such that $ t_\epsilon < T$ for any $ \epsilon\leq\eps_0  $. Assume on the contrary that we have a sequence $ \epsilon_n\longrightarrow 0$ such that $\hat{\bar{u}}^{\epsilon,g}(t_{\eps_n},s_{\epsilon_{n}},z_{\epsilon_n}) \longrightarrow u^{*,g}(T,s_0,z_0)$ and such that $t_{\eps_n}=T$ for countably many $n$. Extracting a further subsequence if necessary, we can assume without loss of generality that the sequence $(t_{\eps_n})_{n}$ is actually stationary at $T$. We then have
\begin{align*}
\bar{u}^{g,\epsilon}(T,s_{\epsilon_{n}},x_{\epsilon_{n}},y_{\epsilon_n})
&=\frac{U_2(z_{\epsilon_{n}}-g(s_{\epsilon_{n}}))-U_2(\ell^{\epsilon_{n}}(x_{\epsilon_{n}},y_{\epsilon_{n}})-g(s_{\epsilon_{n}}))}{\epsilon^2},
\end{align*}
where  $ (x_{\epsilon_{n}},y_{\epsilon_{n}}):=(z_{\epsilon_n}-\ybf^g(T,s_{\epsilon_n},z_{\epsilon_n}),\ybf^g(T,s_{\epsilon_n},z_{\epsilon_n}))$. Since  $ \ybf^g(T,\cdot,\cdot) $ is continuous by Assumption \ref{assumption 1 smoothness}, we have by definition of $ \ell^\epsilon$:
$$ \left| \frac{\ell^{\epsilon_n}(x_{\epsilon_{n}},y_{\epsilon_{n}})-z_{\epsilon_n}}{\epsilon_n^3} \right|=\left| \frac{1}{\epsilon_n^3}\sum_{j=1}^d \left(y_j\left(\frac{1}{1+\eps^3_n\lambda^{j,0}}-1\right){\bf 1}_{y_j\geq 0}+\eps_n^3\lambda^{0,j}y_j{\bf 1}_{y_j<0}\right)\right|\le C, $$
for some constant $ C$ independent of $ n$ and $ \epsilon$. Since $ U^2$ is $ C^1$, we deduce that
$$ \frac{U_2(z_{\epsilon_{n}}-g(s_{\epsilon_{n}}))-U_2(\ell^{\epsilon_{n}}(x_{\epsilon_{n}},y_{\epsilon_{n}})-g(s_{\epsilon_{n}}))}{\epsilon_n^2} \longrightarrow 0,$$
 as $ n \longrightarrow +\infty$, wich contradicts \eqref{eq phi>delta} and \eqref{eq leps converges}.

\vspace{0.25em}
\textit{Step 2:}
Similarly as in Section \ref{sect:subsol interior} and in \cite{st} and \cite{pst}, we build a test function $ \psi^\epsilon$ for $ v^{g,\epsilon}$.

Let $\overline{p}\in(0,1)$ be a constant which will be fixed later. We define $ \psi^\epsilon $ by
$$ \hat{\psi}^\epsilon (t,s,z,\xi) := v^g(t,s,z)-\epsilon^2 ( l_*^\epsilon + \phi(t,s,z)+ \Phi^\epsilon(t,s,z))-\eps^4\hat{\varpi}(\xi),$$
with $ \Phi^\epsilon(t,s,z):= l_0 (T-t)^{\overline{p}}+l_1 ( (s-s_\epsilon)^2+(z-z_\epsilon)^2), \ \hat{\varpi}(\xi):=\abs{\xi}^2,$
for some constants $l_1$ and $l_0$. By definition, we have $ \hat{\bar u}^{\epsilon,g}(t,s,z,0) \le b^* $ for all $(t,s,z)\in B_1$. We now choose $ l_1$ large enough and $ l_0$ so that on $ B_1 \backslash  B_2$ where $ B_2:=  [T-\frac{r_1}{2},T]\times B_{\frac{r_1}{2}}(s_0,z_0)$, we have
$$ \phi(t,s,z)+ \Phi^\eps(t,s,z)+\eps^2\xi^2 \geq 2+ b^*.$$

We then have that $  v^\epsilon-\psi^\epsilon$ has a local minimizer in $ B_1$. Indeed on $ \partial B_1$, for $ \epsilon$ small enough, since $ l_*^\epsilon \longrightarrow 0$, we have:
\begin{align*}
\frac{  v^{\epsilon,g}(t,s,x,y)-\psi^\epsilon(t,s,x,y)}{\epsilon^2} =&\ -u^{\epsilon,g}(t,s,x,y)+l^\epsilon_* +\phi(t,s,z)+\Phi^\epsilon(t,s,z)+(y-\ybf^g(t,s,z))^2\\
 \geq &\ -b^*+l_*^\epsilon + 2+b^*>0.
\end{align*}
Since $v^{\epsilon,g}(t_\epsilon,s_\epsilon,x_\epsilon,y_\epsilon)-\psi^\epsilon(t_\epsilon,s_\epsilon,x_\epsilon,y_\epsilon)=0 $, we then  have the existence of a local minimizer $ ( \tilde{t}_\epsilon,\tilde{s}_\epsilon,\tilde{x}_\epsilon, \tilde{y}_\epsilon) \in B_1$. We denote by $ (\tilde{t}_\eps,\tilde{s}_\eps,\tilde{z}_\eps,\tilde{\xi}_\eps)$ the corresponding minimizer after the usual change of variable. We also recall that by classical results on viscosity solutions, we have $(\tilde{t}_\eps,\tilde{s}_\eps,\tilde{z}_\eps)\longrightarrow (T,s_0,z_0)$ as $\eps$ goes to $0$. Now by the viscosity supersolution property of $ v^\epsilon$ at $ ( \tilde{t}_\epsilon,\tilde{s}_\epsilon,\tilde{x}_\epsilon, \tilde{y}_\epsilon)$, we have, since we recall that we do have $\tilde t_\eps<T$
\begin{align}\label{eq: viscos super sol}
\min_{(i,j)\in\mathcal I} \big\{ k\psi^\eps-\mathcal{L} \psi^\eps-\kappa \tilde{U}_1(\psi^\eps_x) , \ \Lambda^\eps_{i,j} \cdot \left(\psi^\eps_x,\psi^\eps_y \right) \big\} ( \tilde{t}_\epsilon,\tilde{s}_\epsilon,\tilde{x}_\epsilon, \tilde{y}_\epsilon) \geq 0.
\end{align}

\textit{Step 3:} We now show that there exists $ \hat{\eps}$ such that for $ \eps \le \hat{\eps}$ the sequence $ ( \tilde{\xi}_\eps)_{0< \eps \le \hat{\eps}}$ is bounded. Since the sequence $ (\tilde{t}_\eps,\tilde{s}_\eps,\tilde{z}_\eps,\eps \tilde{ \xi}_\eps)$ is bounded, we indeed easily compute that the gradient constraints in \eqref{eq: viscos super sol} implies for $ (i ,j)\in\mathcal I $
$$ \Lambda_{i,j}^\eps \cdot \left( \psi^\eps_x,\psi^\eps_y\right) \left( \tilde{t}_\eps,\tilde{s}_\eps,\tilde{z}_\eps,\tilde{\xi}_\eps\right)=\eps^3 \left(\lambda^{i,j} v^g_z(\tilde{t}_\eps,\tilde{s}_\eps,\tilde{z}_\eps)-2(e_i-e_j)\cdot \tilde{\xi}_\eps \right)+\circ(\eps^3) \geq 0. $$
Then for $ i=0$ and $ j \geq 1$, we obtain, since $\lambda^{0,j}\in\mathcal I$ for any $j\geq 1$, that for $\eps$ small enough
$$ \tilde{\xi}_\eps^j \geq -\lambda^{0,j}v^g_z(\tilde{t}_\eps,\tilde{s}_\eps,\tilde{z}_\eps)>-\text{Const},$$
where $ \text{Const}>0$ is uniform in $ \eps$. Then for $ i\geq 1$ and $ j=0$, we obtain that for $ \eps$ small enough
$$\tilde{\xi}_\eps^i \le \lambda^{i,0} v^g_z(\tilde{t}_\eps,\tilde{s}_\eps,\tilde{z}_\eps) < \text{Const}. $$
Hence $ (\tilde{\xi}^\eps)$ is bounded for $\eps$ small enough.

\vspace{0.25em}
\textit{Step 4:} We now deduce from \eqref{eq: viscos super sol} and Lemma \ref{lem remainder estimate} that at point $(\tilde{t}^\eps,\tilde{s}^\eps, \tilde{x}^\eps,\tilde{y}^\eps)$:
\begin{align}\label{eq remainder estimate interm}
\eps^2\left( -\frac{v_{zz}^g}{2} \left|\sigma^T \tilde{\xi}_\eps  \right|^2 +\Tr{ \alpha^g (\alpha^{g})^{T}}-\mathcal{A}^g (l_*^\eps+\phi+\Phi^\eps)+\mathcal{R}^\eps (l_*^\eps+\phi+\Phi^\eps,\hat\varpi) \right) \geq 0.
\end{align}
Since for $ \eps$ small, $ (\tilde{\xi}_\eps)$ is bounded, $\hat\varpi$ only depends on $\xi$ and $\Phi^\eps$ and all its derivatives with respect to $s$ and $z$ are bounded, we obtain by Lemma \ref{lem remainder estimate} and Assumption \ref{assumption 1 smoothness} that for some ${\rm Const}>0$
$$ \abs{\mathcal{R}^\eps (l_*^\eps+\phi+\Phi^\eps,\hat\varpi)}(\tilde{t}_\eps,\tilde{s}_\eps,\tilde{z}_\eps,\tilde{\xi}_\eps)\leq \eps\frac{{\rm Const}}{(T-\tilde t_\eps)^{1-\eta}}.$$

Now by definition of $ \mathcal{A}^g$ and $ \Phi^\eps$, we observe easily that
$$ \mathcal{A}^g (l_*^\eps+\phi+\Phi^\eps)(\tilde{t}_\eps,\tilde{s}_\eps,\tilde{z}_\eps)=\frac{\overline{p} l_0}{(T-\tilde{t}_\eps)^{1-\overline{p}}}+r^\eps,  $$
where $ r^\eps $ is bounded near 0, so that by \eqref{eq remainder estimate interm} and Assumption \ref{assumption 1 smoothness}(i), we obtain
$$ -\frac{\overline{p} l_0}{(T-\tilde{t}_\eps)^{1-\overline{p}}} +\frac{{\rm Const}}{(T-\tilde t_\eps)^{1-\mu}}+\eps\frac{{\rm Const}}{(T-\tilde t_\eps)^{1-\eta}}+\tilde{r}^\eps \geq 0,$$
where $ \tilde{r}^\eps:=-r^\eps+\Tr{ \alpha^g (\alpha^g)^T},$ is bounded near 0. Choosing $\overline{p}=(\eta\wedge\mu)/2$, this leads to a contradiction for $ \eps>0$ small enough, since $\tilde t_\eps$ goes to $T$.

\vspace{-2em}

\begin{flushright}$\Box$\end{flushright}

\subsection{Viscosity supersolution}
\label{sect:supersol interior}

We are interested in this section in the supersolution part. We first note that since $ \bar{u}^{\eps,g} \geq 0$, the supersolution property on $ \big\{T\big\}\times \R^d \times \R_+$ is indeed trivial. We then only focus on the interior of the domain. Our aim is then to show:

\begin{Proposition}\label{Prop supersol}
Under Assumptions \ref{assumption 1 smoothness}, \ref{assumption 4 bound on derivatives} and \ref{assumption 2 local bound}, $ u^g_*$ is a viscosity supersolution of the second corrector equation \reff{eq:corr2} on $ [0,T)\times(0,+\infty)^{d+1}$.
\end{Proposition}

We first recall some crucial properties proved in \cite{pst}, that we shall use in the proof of Proposition \ref{Prop supersol}. The first one concerns a regular approximation of $ \tilde w^g$ by convolution. Consider $ \upsilon:\R^d \longrightarrow \R$ a positive, even, $ C^\infty $ kernel with support in $ B_1(0)$. We then define for $ m>0$:
$$ \tilde w^{g,m}(\cdot,\xi):=\int_{\R^d}\upsilon^m \left( \zeta\right)\tilde w^g(\cdot,\xi-\zeta) d\zeta,  $$
where $\upsilon^m(x):=m^{-d}\upsilon(x/m)$. The proof of the following lemma can be found in \cite{pst}:
\begin{Lemma}\label{Lem properties of wm}
Under Assumtion \ref{assumption 4 bound on derivatives}, we have for any $ m>0$ that:

$\rm{(i)}$ $\tilde w^{g,m}$ is $C^2$, convex in $\xi$ and for any $(t,s,z,\xi)\in[0,T)\times(0,+\infty)^{d+1}\times\R^d$,
$$0\leq \tilde w^{g,m}(t,s,z,\xi)\leq Lv^g_z(t,s,z)(1+m)(1+\abs{\xi}).$$

\vspace{0.25em}
$\rm{(ii)}$ $\tilde w^{g,m}$ is smooth in $(t,s,z)\in[0,T)\times(0,+\infty)^{d+1}$, and satisfies the following, uniformly in $m$, 
	\begin{align}\label{estim1h}
\nonumber&\left(\abs{\tilde w^{g,m}}+\abs{\tilde w^{g,m}_t}+\abs{\tilde w^{g,m}_s}+\abs{\tilde w^{g,m}_{ss}}+\abs{\tilde w^{g,m}_z}+\abs{\tilde w^{g,m}_{sz}}+\abs{\tilde w^{g,m}_{zz}}\right)(\cdot,\xi)\leq C(\cdot)(1+m)\left(1+\abs{\xi}\right)\\
\nonumber&\left(\abs{\tilde w^{g,m}_\xi}+\abs{\tilde w^{g,m}_{s\xi}}+\abs{\tilde w^{g,m}_{z\xi}}\right)(\cdot,\xi)\leq C(\cdot)\\
&\abs{\tilde w^{g,m}_{\xi\xi}}(\cdot,\xi)\leq C(\cdot){\bf1}_{\xi\in \overline{\mathcal B}^g(\cdot)},
\end{align}
where $C(t,s,z)$ is a continuous function depending on the Merton value function and its derivatives, and $\mathcal B^g(t,s,z)$ is some ball with a continuous radius, centered at $0$. 

\vspace{0.25em}
$\rm{(iii)}$ For every $(i,j)\in\mathcal I$ and every $(t,s,z,\xi)\in[0,T)\times(0,+\infty)^{d+1}\times\R^d$
$$-\lambda^{i,j}v^g_z(t,s,z)+\tilde w^{g,m}_{\xi_i}(t,s,z,\xi)-\tilde w^{g,m}_{\xi_j}(t,s,z,\xi)\leq 0.$$

$\rm{(iv)}$ For every $(t,s,z,\xi)\in[0,T)\times(0,+\infty)^{d+1}\times\R^d$, we have
 $$
 \frac12v^g_{zz}(t,s,z)\int_{\mathbb R^d}
                   \upsilon^m(\zeta)\abs{\sigma(t,s)^T(\xi-\zeta)}^2d\zeta
 -\frac12\Tr{\alpha^g(\alpha^g)^T\tilde w^{g,m}_{\xi\xi}}(t,s,z,\xi)+a^g(t,s,z)\leq 0.
 $$
\end{Lemma}

To build a test function in the proof of Proposition \ref{Prop supersol} we will also use the following result.

\begin{Lemma}\label{Lem properties of hdelta}
For any $\delta\in(0,1)$ and any $\nu>0$, there exists $a^\delta:=a^{\delta,\nu}>1$ and a function $h^{\delta,\nu}:\mathbb R^d\longrightarrow [0,1]$ such that $h^{\delta,\nu}$ is $C^\infty$, $h^{\delta,\nu}=1$ on $B_1(0)$ and $h^{\delta,\nu}=0$ on $B_{a^\delta}(0)^c$. Moreover, for any $1\leq i,j\leq d$ and for any $\xi\in\mathbb R^d$
 $$
\abs{h^{\delta,\nu}_{\xi_i}(\xi)}\leq \frac{\nu\delta}{3}, \ \abs{\xi} |h^{\delta,\nu}_{\xi_i}|
 \leq 
 \nu\delta,
 \text{ and }
 \abs{\xi}|h^{\delta,\nu}_{\xi\xi}|+|h^{\delta,\nu}_{\xi\xi}|\leq C^*,
 $$
for some constant $C^*$ independent of $\delta$.
\end{Lemma}
This Lemma and its proof can be found in \cite{pst}. We conclude these preliminary results with the following useful lemma, which follows directly from Lemmas \ref{Lem properties of wm} and \ref{Lem properties of hdelta}.

\begin{Lemma}\label{Lem prop of wH}
For any $\delta\in(0,1)$, $\nu>0$ and $ m>0$, the map $ \Upsilon:=\tilde{w}^{g,m} h^{\delta,\nu}$ is smooth and satisfies the following estimates 
 \begin{align}\label{estim2h}
\nonumber&\left(\abs{\Upsilon}+\abs{\Upsilon_t}+\abs{\Upsilon_s}+\abs{\Upsilon_{ss}}+\abs{\Upsilon_z}+\abs{\Upsilon_{sz}}+\abs{\Upsilon_{zz}}\right)(t,s,z,\xi)\leq C(t,s,z)(1+m)\left(1+\abs{\xi}\right){\scriptstyle\1_{|\xi|\leq a^\delta}}\\
\nonumber&\left(\abs{\Upsilon_\xi}+\abs{\Upsilon_{s\xi}}+\abs{\Upsilon_{z\xi}}\right)(t,s,z,\xi)\leq 4C(t,s,z)\left( 1+(1+m)(1+|\xi|)\frac{\nu \delta \sqrt{d}}{3}\right){\bf1}_{|\xi|\leq a^\delta}\\
&\abs{\Upsilon_{\xi\xi}(t,s,z,\xi)}\leq C(t,s,z)\left(1+2\frac{\nu \delta \sqrt{d}}{3}+C^*(1+m)(1+|\xi|) \right) {\bf1}_{|\xi|\leq a^\delta},
\end{align}
where $C(t,s,z)$ and $C^*$ were introduced in Lemmas \ref{Lem properties of wm} and \ref{Lem properties of hdelta}.
\end{Lemma}

\textbf{Proof of Proposition \ref{Prop supersol}}.

\vspace{0.25em}
Let $ (t_0,s_0,z_0)\in [0,T)\times (0,+\infty)^{d+1}$ and $ \phi$, $ C^2$ s.t., $\forall  (t,s,z)\in [0,T)\times (0,+\infty)^{d+1} \backslash \left\{ (t_0,s_0,z_0) \right\}$:
$$ 0=(u^{g}_*-\phi)(t_0,s_0,z_0)<(u^g_*-\phi)(t,s,z).$$
We want to show that $ \mathcal{A}^g\phi(t_0,s_0,z_0)-a^g(t_0,s_0,z_0) \geq 0$. Assume on the contrary that:
\begin{align}\label{eq supersol interm 1}
 \mathcal{A}^g\phi(t_0,s_0,z_0)-a^g(t_0,s_0,z_0) < 0,
\end{align}
Then there exists $ r_0>0$ such that $  \mathcal{A}^g\phi(t,s,z)-a^g(t,s,z) \le 0$ on $ B_{r_0}(t_0,s_0,z_0)$.

\vspace{0.25em}
We proceed in 5 steps. The first two steps consist in defining a test function for the dynamic programming equation \eqref{e.dpp}. The third one is devoted to prove that the gradient constraint for this test function is not binding, so that the parabolic part is. The last two steps lead to the required contradiction of \eqref{eq supersol interm 1}.

\vspace{0.25em}
\textit{Step 1:}
By Lemma \ref{Lem: u*g dep of (t,s,z)}, there exists a sequence $ (t^\eps,s^\eps,z^\eps)\longrightarrow(t_0,s_0,z_0)$ when $ \eps\longrightarrow 0$ such that
$$ \hat{u}^{\eps,g}(t^\eps,s^\eps,z^\eps,0) \underset{\eps \longrightarrow 0}{\longrightarrow} u^{g}_*(t_0,s_0,z_0).$$
We have $ l^\eps_*:= \hat{u}^{\eps,g}(t^\eps,s^\eps,z^\eps,0)-\phi(t^\eps,s^\eps,z^\eps) \longrightarrow 0$ and $(x^\eps,y^\eps)\longrightarrow (x_0,y_0)$, as $\eps$ goes to $0$, where
\begin{align*}
(x^\eps,y^\eps)&:=\left(z^\eps-\ybf^g(t^\eps,s^\eps,z^\eps)\cdot \1_d,\ybf^g(t^\eps,s^\eps,z^\eps)\right),\\(x_0,y_0)&:=\left(z_0-\ybf^g(t_0,s_0,z_0)\cdot \1_d,\ybf^g(t_0,s_0,z_0)\right).
\end{align*}

 We then consider $ \eps_0>0$ such that for all $ \eps \le \eps_0$
 $$ | t^\eps-t_0| + |s^\eps-s_0|+  |z^\eps-z_0|\le \frac{r_0}{4}, \ \text{and} \  | l^\eps_*| \le 1. $$
 
Consider next a constant $ q_0>0$ such that:
$$ \sup_{(t,s,z)\in B_{r_0/2}(t_0,s_0,z_0)} \big\{ \phi(t,s,z)+C(t,s,z) \big\}+3 \le q_0 \left( \frac{r_0}{12}\right)^4,$$
where $C(t,s,z)$ is the continuous function appearing in \reff{estim1h}. We then introduce:
$$ \phi^\eps(t,s,z):=\phi(t,s,z)-q_0 \left(|t-t^\eps|^4+|z-z^\eps|^4+|s-s^\eps|^4 \right).$$
We then have for $ \eps \le \eps_0$ and $ (t,s,z)\in \partial B_{r_0/2}(t_0,s_0,z_0)$ that $|t-t^\eps|+|z-z^\eps|+|s-s^\eps|\geq r_0/4,$
and thus that 
$$|t-t^\eps|^4+|z-z^\eps|^4+|s-s^\eps|^4\geq \frac{1}{81}\left(|t-t^\eps|+|z-z^\eps|+|s-s^\eps|\right)^4\geq \left(\frac{r_0}{12}\right)^4.$$
Then on $ \partial  B_{r_0/2}(t_0,s_0,z_0) $, we have:
\begin{align}\label{eq:3bis}
\nonumber
\phi^\eps(t,s,z)+l^*_\eps +C(t,s,z)
=&\ 
\phi(t,s,z)+C(t,s,z)-q_0 \left(|t-t^\eps|^4+|z-z^\eps|^4+|s-s^\eps|^4 \right)+l^*_\eps
\\
\nonumber
\leq & \
q_0\left(\frac {r_0}{12}\right)^4-3-q_0\left(\abs{t-t^\eps}^4+\abs{z-z^\eps}^4+\abs{s-s^\eps}^4\right)+1
\\
\leq & \ q_0\left(\left(\frac {r_0}{12}\right)^4-\abs{t-t^\eps}^4-\abs{z-z^\eps}^4-\abs{s-s^\eps}^4\right)-2\leq -2.
\end{align}

Consider next the function $ \Phi^\eps:=\phi^\eps-\phi+l^*_\eps$. By linearity of the operator $ \mathcal{A}^g$, and Assumption \ref{assumption 1 smoothness}, we have that there exists $\eps^0>0$ such that $b<\infty$, where
\begin{equation}\label{eq:b} b:=\sup \big\{ \left| \mathcal{A}^g \Phi^\eps\right| (t,s,z), \ \eps\le \eps^0, \ (t,s,z)\in \bar{B}_{r_0/2}(t_0,s_0,z_0) \big\}.\end{equation}

Throughout the rest of the proof, we let $m\in(0,1]$. Now for any $ \delta \in (0,1)$ and $ \nu>0$, let $ h^{\delta,\nu}$ be the function defined by Lemma \ref{Lem properties of hdelta}, and introduce a parameter $ \xi^*:= 1 \vee \tilde{\xi}_0 \vee \tilde{\xi}^*$, where $ \tilde{\xi}_0>0$ is greater than $\eta^g(t_0,s_0,z_0)$ times the diameter of $ \mathcal{O}^g(t_0,s_0,z_0)$ and large enough so that for every $|\xi|\geq \tilde \xi_0$, $\tilde w^{g,m}_{\xi\xi}(t,s,z,\xi)=0$, for every $(t,s,z)\in B_{r_0/2}(t_0,s_0,z_0)$, and $ \tilde{\xi}^*$ is such that for any $ \xi \in B_{\tilde{\xi}^*}(0)^c$ and $ (t,s,z)\in \bar{B}_{r_0/2}(t_0,s_0,z_0)$, we have
\begin{align}\label{eq interm step 1}
-\frac12 v^{g}_{zz}\abs{\sigma^T \xi}^2
-\mathcal A^g(\phi+\Phi^\eps)  > \frac{1}{2}\Tr{\alpha^g(\alpha^g)^T}C(t,s,z)\left( C^*+\frac{ \sqrt{d} \delta \underline{\lambda}}{4L}  \right)+1,
\end{align} 
where $ C^*$ is the constant introduced in Lemma \ref{Lem properties of hdelta} and $ C(t,s,z)$ is the function introduced in Lemma \ref{Lem properties of wm} (we remind the reader that they are both uniform in $m$). Define then $H(\xi):=h^{\delta,\nu}\left(\frac{\xi}{\xi^*} \right)$ and the test function
 $$ \hat{\psi}^{\eps,\delta,m}(t,s,z,\xi):=v^g(t,s,z)-\eps^2 (\phi^\eps(t,s,z)+ l_\eps^*)-\eps^4(1-\delta) \tilde w^{g,m}(t,s,z,\xi)H(\xi).$$

\textit{Step 2:} In this part, we introduce a second modification of the test function. Introduce
\begin{align*}
I^{\eps,\delta,m}(t,s,z,\xi):=\eps^{-2} \left(\hat{v}^{\eps,g}-\hat{\psi}^{\eps,\delta,m} \right)(t,s,z,\xi),
\end{align*}
we want to show that $ I^{\eps,\delta,m}$ has a local maximizer on the interior of the domain. By definition,
\begin{align*}
 I^{\eps,\delta,m}(t,s,z,\xi)=&\ \phi^\eps(t,s,z)-\hat{u}^{\eps,g}(t,s,z,\xi)+l^\eps_*-\eps^2w^g(t,s,z,\xi)+\eps^2(1-\delta)H(\xi)\tilde w^{g,m}(t,s,z,\xi). 
 \end{align*}
 Recall that for $ \xi=0$, $ w^g(\cdot,\cdot,\cdot,0)=0$, so that by definition of $ l^\eps_*$, we have $$ I^{\delta,\eps,m}(t^\eps,s^\eps,z^\eps,0)=\eps^2(1-\delta)\tilde w^{g,m}(t^\eps,s^\eps,z^\eps,0),$$
 which goes to $0$ as $\eps$ goes to $0$, uniformly in $m\in(0,1]$, because of the uniform bounds given by Lemma \ref{Lem properties of wm}. Hence, there exists $\eps_1$ such that for any $\eps\leq \eps_0\wedge\eps^0\wedge\eps_1$,
 \begin{equation}\label{eq:Ieps}
 I^{\delta,\eps,m}(t^\eps,s^\eps,z^\eps,0)\geq -1.
 \end{equation}
 
Using successively that $ v^{\eps,g} \le v^g$, $ 0 \le \tilde w^{g,m}(t,s,z,\xi) \le 2C(t,s,z)(1+ |\xi|) $ (with $C(t,s,z)$ still being the continuous function appearing \reff{estim1h} and where we used the fact that $m\in(0,1]$) and $ 0 \le H(\xi) \le \1_{|\xi| \le a^\delta \xi^*}$, we have
\begin{align*}
 I^{\eps,\delta,m}(t,s,z,\xi) &\le \phi^\eps(t,s,z)+l^\eps_*+\eps^2(1-\delta)H(\xi)\tilde w^{g,m}(t,s,z,\xi) \\
 & \le \phi^\eps(t,s,z)+l^\eps_*+2\eps^2(1-\delta)C(t,s,z)H(\xi) (1+|\xi |) \\
 & \le  \phi^\eps(t,s,z)+l^\eps_*+2\eps^2  C(t,s,z)(1+a^\delta \xi^* ),
\end{align*}
where $ a^\delta$ is the constant introduced in Lemma \ref{Lem properties of hdelta}. Then for any $ \eps \le \eps^\delta:= \sqrt{2(1+a^\delta  \xi^*)}$, we have
$$  I^{\eps,\delta,m}(t,s,z,\xi) \le \phi^\eps(t,s,z)+l^\eps_*+ C(t,s,z). $$
Introduce then $ Q_{(t_0,s_0,z_0)}:=\big\{(t,s,z,\xi), \ (t,s,z)\in \bar{B}_{r_0/2}(t_0,s_0,z_0) \big\}$. The above implies in particular that for $\eps\leq  \eps_0 \wedge\eps^0\wedge\eps_1\wedge \eps^\delta$
$$ \mathfrak{I}(\eps,\delta,m):=\sup_{(t,s,z,\xi)\in Q_{(t_0,s_0,z_0)}} I^{\eps,\delta,m}(t,s,z,\xi) <\infty.$$
Moreover, using \eqref{eq:3bis}, we deduce that for $ (t,s,z,\xi)\in \partial Q_{(t_0,s_0,z_0)}$ and for any $\eps\leq  \eps_0 \wedge\eps^0\wedge\eps_1\wedge \eps^\delta$
\begin{equation}\label{eq:-2} I^{\eps,\delta,m}(t,s,z,\xi)\le -2.\end{equation}
We can now consider a sequence $ (\hat{t}_n,\hat{s}_n,\hat{z}_n,\hat{\xi}_n)$ in $ \text{int} \left(Q_{(t_0,s_0,z_0)}\right)$ such that
$$ I^{\eps,\delta,m} (\hat{t}_n,\hat{s}_n,\hat{z}_n,\hat{\xi}_n) \geq \mathfrak I(\eps,\delta,m)-\frac{1}{2n}. $$
It is now time to penalize the test function to obtain the existence of an interior maximiser, which is not obvious with our previous construction. We consider $ f:\R\longrightarrow[0,1]$, smooth such that $ f(0)=1$ and $ f(x)=0$ if $ x \geq 1$. Define
$$ \hat{\psi}^{\eps,\delta,m,n}(t,s,z,\xi):=\hat{\psi}^{\eps,\delta,m}(t,s,z,\xi)-\frac{\eps^2}{n}f\left( |\xi-\hat{\xi}_n | \right).   $$
Consider then
$$
 I^{\eps,\delta,m,n}(t,s,z,\xi)
 :=
 \eps^{-2}\left(\hat v^{\eps,g}-\hat\psi^{\eps,\delta,m,n}\right)(t,s,z,\xi)
 =
 I^{\eps,\delta,m}(t,s,z,\xi)+\frac1n f\left(|\xi-\hat\xi_n|\right).
 $$
By definition of $ (\hat{t}_n,\hat{s}_n,\hat{z}_n,\hat{\xi}_n)$, we have for any $ (t,s,z,\xi)\in Q_{(t_0,s_0,z_0)}$
$$  I^{\eps,\delta,m,n} (\hat{t}_n,\hat{s}_n,\hat{z}_n,\hat{\xi}_n)=I^{\eps,\delta,m}(\hat{t}_n,\hat{s}_n,\hat{z}_n,\hat{\xi}_n)+\frac{1}{n}\geq I^{\eps,\delta,m}(t,s,z,\xi)+\frac{1}{2n}.$$
Notice that for $ | \xi- \hat{\xi}_n|\geq 1$, we have $ I^{\eps,\delta,m,n}(t,s,z,\xi)=I^{\eps,\delta,m}(t,s,z,\xi).$ For $ n$ large enough, we then have that
$$ \sup_{(t,s,z,\xi)\in Q_{(t_0,s_0,z_0) }} I^{\eps,\delta,m,n}(t,s,z,\xi)=\sup_{(t,s,z,\xi)\in Q^n_{(t_0,s_0,z_0)}}I^{\eps,\delta,m,n}(t,s,z,\xi),$$
where $ Q^n_{(t_0,s_0,z_0)}:=\big\{ (t,s,z,\xi), |\xi-\hat{\xi}_n |\le 1, \ (t,s,z)\in Q_{(t_0,s_0,z_0)} \big\}$ is compact. Then since $ I^{\eps,\delta,m,n}$ is continuous, this implies the existence of $ (t_n,s_n,z_n,\xi_n)\in  Q^n_{(t_0,s_0,z_0)} $ which maximises $ I^{\eps,\delta,m,n}$. We also observe that $  (t_n,s_n,z_n,\xi_n)\in  \text{int}\left( Q_{(t_0,s_0,z_0)} \right) $. Indeed, it is clear that we have for $\eps\leq  \eps_0\wedge\eps^0 \wedge\eps_1\wedge \eps^\delta$
$$  I^{\eps,\delta,m,n} (t_n,s_n,z_n,\xi_n) \geq  I^{\eps,\delta,m,n}(t^\eps,s^\eps,z^\eps,0)\geq I^{\eps,\delta,m}(t^\eps,s^\eps,z^\eps,0)\geq -1, $$
and for $ (t,s,z,\xi)\in \partial Q_{(t_0,s_0,z_0)}$, we have using \reff{eq:-2}
$$  I^{\eps,\delta,m,n}(t,s,z,\xi) \le  I^{\eps,\delta,m}(t,s,z,\xi)+\frac{1}{n} \le -2+\frac{1}{n} <-1 , \ \text{for} \ n> 1. $$
We then have for $n>1$ and $\eps\leq  \eps_0 \wedge\eps^0\wedge\eps_1\wedge \eps^\delta$, that 
 by the viscosity subsolution property of $ v^{\eps,g}$ at the point $ (t_n,s_n,z_n,\xi_n)$ (with corresponding $ (t_n,s_n,x_n,y_n)$)
\begin{align}\label{interm visco sub sol veps}
\underset{(i,j)\in\mathcal I}{\min}
 \left\{\ k \psi^{\eps,\delta,m,n}- \Lc \psi^{\eps,\delta,m,n} - \Ut_1(\psi^{\eps,\delta,m,n}_{x}) , 
        \ \Lambda_{i,j}^\eps\cdot(\psi^{\eps,\delta,m,n}_x,\psi^{\eps,\delta,m,n}_y) 
 \right\} 
\le
 0.
\end{align}

\vspace{0.3em}
\textit{Step 3.} We now show that for $ \eps$ small enough, and $ n$ large enough,
$$ D^{i,j}:=\Lambda_{i,j}^\eps\cdot(\psi^{\eps,\delta,m,n}_x,\psi^{\eps,\delta,m,n}_y)(t_n,s_n,x_n,y_n)>0 \ \text{for all} \ (i,j)\in\mathcal I. $$
It is easy to compute that for $ ( i , j)\in\mathcal I $, we have $ D^{i,j}= \eps^3 G^\eps
 -E^\eps-F^{\eps,n},$ with
\begin{align*} G^\eps:=&\ \left[\lambda^{i,j}v_z^g(t_n,s_n,z_n)
             -(1-\delta)(\tilde w^{g,m}H)_{\xi}(t_n,s_n,z_n,\xi_n).(e_i-e_j)
       \right],\\
 E^\eps:= &\ \lambda^{i,j}\eps^7(1-\delta)(\tilde w^{g,m}_zH)(t_n,s_n,z_n,\xi_n)+\lambda^{i,j}\left[\eps^5(\phi_z(t_n,s_n,z_n)-4q_0 |z_n-z^\eps|^3)\right.\\
 & \left.+\eps^6(1-\delta)(\tilde w^{g,m}H)_{\xi}(t_n,s_n,z_n,\xi_n).(e_i-\ybf^g_z(t_n,s_n,z_n))\right],\\
 F^{\eps,n}:=&\ \frac\eps n\frac{f'(|\xi_n-\hat\xi_n|)}{|\xi_n-\hat\xi_n|}\left( \xi_n- \hat{\xi}_n \right) \cdot \left(e_i-e_j+\lambda^{i,j}\eps^3(e_i-\ybf^g_z(t_n,s_n,z_n))\right).
 \end{align*}
By the properties of $ h^{\delta,\nu}$ obtained in Lemma \ref{Lem properties of hdelta} and the estimates of Lemmas \ref{Lem properties of wm}, \ref{Lem prop of wH}, we have:
\begin{align*}
\left| E^\eps \right| \le &\ \lambda^{i,j} \eps^5 \Big[ |\phi_z|(t_n,s_n,z_n)+4q_0 \left| z_n-z^\eps \right|^3+ 2\eps^2 C(t_n,s_n,z_n)(1+ |\xi_n |) \1_{|\xi_n|< a^\delta \xi^* } \\
&+  \eps C_1(t_0,s_0,z_0)C(t_n,s_s,z_n)\left( 1+2(1+|\xi_n|)\frac{\nu \delta \sqrt{d}}{3\xi^*}\right){\bf1}_{|\xi_n|\leq a^\delta\xi^*} \Big]\\
 \le &\ C_2(t_0,s_0,z_0) \eps^5 \left[1+\eps \left(1+\nu \delta a^\delta  \right) +\eps^2 a^\delta \xi^* \right]
\end{align*}
for some functions $ C_1(t_0,s_0,z_0)$, $ C_2(t_0,s_0,z_0)$ which depend on $\ybf^g_z$, $\phi_z$ and the function $C$.

\vspace{0.25em}
Similarly, recalling that $| \xi_n-\hat{\xi}_n| \le 1$, we obtain easily for some constant $ C_3(t_0,s_0,z_0)$, which depends on $\ybf^g_z$
$$ \left|  F^{\eps,n} \right|\le C_3(t_0,s_0,z_0) \frac{\eps}{n}.$$
We then study $ G^\eps$. By Lemma \ref{Lem properties of wm}(i) and (iii), we have
\begin{align*}
G^\eps &=\lambda^{i,j}v^g_z-(1-\delta)(\tilde w^{g,m}_{\xi_i}-\tilde w^{g,m}_{\xi_j})H-(1-\delta)w^m(H_{\xi_i}-H_{\xi_j})\\
&\geq \lambda^{i,j}v^g_z-\lambda^{i,j}(1-\delta)v^g_z-(1-\delta)Lv^g_z(1+\abs{\xi_n})\left(\abs{H_{\xi_i}}+\abs{H_{\xi_j}}\right)\\
& \geq \lambda^{i,j}v^g_z\left(\delta-\frac{L(1+\abs{\xi_n})}{\underline\lambda}(1-\delta)\left(\abs{H_{\xi_i}}+\abs{H_{\xi_j}}\right)\right).
\end{align*}
We then fix the value of $ \nu$ of $ h^{\nu,\delta}$ by $ \nu:= \frac{\underline{3\lambda}}{8L}$, so that by Lemma \ref{Lem properties of hdelta}, we have for all $0 \le i \le d $, $ (1+\left| \xi\right|) H_{\xi_i}(\xi) \le \frac{\underline{\lambda} \delta}{2L}$, and $ G^\eps \geq \lambda^{i,j} \delta^2 v^g_z.$
Notice that the choice of $ \nu$ only depends on $ \underline{\lambda}$ and $ L$, so that the previous constants are not affected by our choice of $ \nu$. Since the sequence $ (t_n,s_n,z_n)$ lives in a compact set and $ v_z^g>0$ and is continuous, we obtain for some constant $ C_4(t_0,s_0,z_0)>0$ that
$$ D^{i,j}\geq \eps^3 \lambda^{i,j} \delta^2 C_4(t_0,z_0,z_0)-  C_2(t_0,s_0,z_0) \eps^5 \left[1+\eps \left(1+\nu \delta a^\delta  \right) +\eps^2 a^\delta \xi^* \right]- C_3(t_0,s_0,z_0) \frac{\eps}{n}. $$
Then for some constant $ \tilde{C}$ and some $\tilde\eps^\delta$, we have for all $ \eps \le \tilde{\eps}^\delta$ and all $ n \geq \tilde{C}\eps^{-5}$ that 
$$ \eps^3 \lambda^{i,j} \delta^2 C_4(t_0,z_0,z_0)-  C_2(t_0,s_0,z_0) \eps^5 \left[1+\eps \left(1+\nu \delta a^\delta \xi^*  \right) +\eps^2 a^\delta \xi^* \right]- C_3(t_0,s_0,z_0) \frac{\eps}{n}>0,$$
so that $ D^{i,j}>0$. By the arbitrariness of the pair $ (i,j)\in\mathcal I$, we obtain that for $ \eps$ small enough and $ n$ large enough
\begin{align}\label{eq interm 5}
  \eps^2 J^{\eps,\delta,m,n}:=\left( k \psi^{\eps,\delta,m,n}- \Lc \psi^{\eps,\delta,m,n} - \Ut_1(\psi^{\eps,\delta,m,n}_{x}) \right)(t_n,s_n,x_n,y_n) \le 0.
\end{align}

\vspace{0.3em}

\textit{Step 4:} We estimate the remainder associated to \reff{eq interm 5}. By the calculations of Lemma \ref{lem remainder estimate}
\begin{align*}
J^{\eps,\delta,m,n}= &-\frac12 v^{g}_{zz}(t_n,s_n,z_n)\abs{\sigma^T(t_n,s_n)\xi_n}^2-\mathcal A^g(\phi+\Phi^\eps)(t_n,s_n,z_n)+\mathcal{R}^{\eps,\delta,m,n}\\
&+\frac12(1-\delta)
\Tr{\alpha^g(\alpha^g)^T(t_n,s_n,z_n)\left( \tilde w^{g,m} H \right)_{\xi\xi}(s_n,z_n,\xi_n)},
&,
\end{align*}
where for all $ 0<\eps \le 1$,
\begin{align}\label{eq interm 6}
 |\mathcal{R}^{\eps,\delta,m,n}| \le C^\delta\left( \eps(|\xi_n|+\eps|\xi_n|^2 +\eps^4+\eps^3+\eps^2n^{-1}\right),
\end{align}
where $ C^\delta$ is a uniform constant depending only on $ \delta$. Indeed, we know by Lemma \ref{Lem prop of wH} that $ \tilde{w}^{g,m}H$ has the required estimates for the evaluation of the remainder estimate in Lemma \ref{lem remainder estimate} (it is easy to do the correspondance between $ \tilde{w}^{g,m}h^{\delta,\nu}$ and $ \tilde{w}^{g,m}H$). Then by Lemma \ref{lem remainder estimate} and since $\eps \xi=y-\ybf^g$ is bounded on the ball $B_{r_0/2}(t_0,s_0,z_0)$, we see that, uniformly in $ m$ and $ \delta$, we have
$$ \mathfrak H^\eps(t,s,z,\xi) \le |v^g_z|+\eps^2\left( \left| \phi_z \right|+\left| \Phi^\eps_z \right|+\eps^2 \left( |\tilde{w}^{g,m}_z|+ \frac{1}{\eps} |\ybf^g_z| |\left(\tilde{w}^{g,m} H\right)_\xi | +\frac{1}{n\eps^3} |\ybf^g_z|K \right)\right) ,$$
where $ \mathfrak H^\eps$ was introduced in the proof of Lemma \ref{lem remainder estimate}, and $ K_1$ is a uniform constant. Then the quantity $ \zeta^\eps(t,s,z,\xi)$ is uniformly bounded on $B_{r_0/2}(t_0,s_0,z_0)$, uniformly in $ m$, but not in $ \delta$. Similarly, we obtain easily that the quantity $ \eps^4\mathfrak R^\eps \left( \left(\tilde{w}^{g,m}H\right)+\frac{1}{\eps^2n}f \right)$ is bounded by some constant depending only on $ (t_0,s_0,z_0)$ and $ \delta$.
 Hence we have for some function $ \tilde{K}^\delta$ depending on $ \ybf^g$ (and its derivatives), the constant $ C(t,s,z)$ introduced in Lemma \ref{Lem properties of wm}, $ K_1$ and $ \delta$ that, uniformly in $ m$:
$$ \left| \mathcal{R}^{\eps,\delta,m,n} \right| \le \tilde K^\delta(t_0,s_0,z_0)(\eps |\xi_n|+\eps^2|\xi_n|^2)+\eps^4\tilde K^\delta(t_0,s_0,z_0)(1+\eps^{-1}+n^{-1}\eps^{-2}).$$

 We then need to show that $ (\xi_n)$ remains bounded as $\eps$ goes to $0$. Indeed, we have from Lemma \ref{Lem properties of hdelta} and Lemma \ref{Lem properties of wm}:
\begin{align*}
\left| \left( \tilde w^{g,m} H \right)_{\xi \xi}\right|(t,s,z,\xi) &\le C(t,s,z)\left(1+2\frac{\nu \delta \sqrt{d}}{3\xi^*}+\frac{2C^*}{(\xi^*)^2}(1+|\xi|) \right) {\bf1}_{|\xi|\leq a^\delta\xi^*}\leq \check C(t,s,z),
\end{align*}
where $ \check{C}$ is a continuous functions depending on the function $C$ appearing in Lemma \ref{Lem properties of wm} and $\xi^*$. Now since we only consider elements $ (t,s,z)\in Q_{(t_0,s_0,z_0)}$ compact, we have by \eqref{eq interm 5}:
\begin{align*}
(-v^g_{zz}(t_n,s_n,z_n))\frac{\abs{\sigma^T\xi_n}^2}{2}+\mathcal{R}^{\eps,\delta,m,n} &\leq\left(\mathcal A^g(\phi+\Phi^\eps)+\abs{\alpha^g(\alpha^g)^T}\check{C}\right)(t_n,s_n,z_n)\\
&\leq \underset{(t,s,z)\in Q_{(t_0,s_0,z_0)}}{\sup}\left(\mathcal A^g(\phi+\Phi^\eps)+\abs{\alpha^g (\alpha^{g})^T}\check{C}\right) \leq {\rm{Const}}.
\end{align*}
Then with \eqref{eq interm 6}, we obtain for $$ \check{C}_1:= \underset{(t,s,z)\in Q_{(t_0,s_0,z_0)}}{\sup} \check{C}(t,s,z) \text{ and } C_0:=\underset{(t,s,z)\in Q_{(t_0,s_0,z_0)}}{\inf}(-v^g_{zz})(t,s,z)>0,$$ since $ v^g$ is strictly concave in $ z$ and smooth, that
$$ C_0 |\xi_n|^2-C^\delta\left( \eps(|\xi_n|+\eps|\xi_n|^2 +\eps^4+\eps^3+\frac{\eps^2}{n}\right) \le  \check{C}_1 .$$
Assume next that $ |\xi_n|$ goes to $+ \infty$ (up to a subsequence) as $\eps$ goes to $0$. Then, the left-hand side above would go to $ +\infty$, which contradicts the fact that it is bounded. Then the sequence $ (\xi_n)$ is bounded by some $ \hat{\xi}^\delta$, depending only on $ \delta$, and not on $m$ since we then know that $ \mathcal{R}^{\eps,\delta,m,n}\longrightarrow 0$ when $ \eps \longrightarrow 0 $ and $ n \longrightarrow \infty$, uniformly in $ m\in(0,1]$.

\vspace{0.25em}
\textit{Step 5:} Following \cite{pst}, we show that $ |\xi_n| \le \xi^*$ for $ n $ large enough (uniformly in $ m$, but not in $ \delta$), where $ \xi^*$ was introduced in Step 1 of the proof. We consider now $ n$ large enough and $ \eps$ small enough so that for all  $ m\in(0,1]$, $|\mathcal{R}^{\eps,\delta,m,n} | \le 1$. Assume on the contrary that $ | \xi_n|>\xi^*$.  By our choice of $ \xi^*$, we know that $ \tilde{w}^{g,m}_{\xi\xi}(\cdot,\xi_n)=0$ on $ Q_{(t_0,s_0,z_0)}$. In the following we omit the dependance w.r.t. $ (t_n,s_n,z_n,\xi_n)$. From \eqref{eq interm 5}, we have, using in particular the fact that $\xi^*\geq 1$
\begin{align*}
-\frac12 v^{g}_{zz}\abs{\sigma^T\xi_n}^2
-\mathcal A^g(\phi+\Phi^\eps) &\le -\frac{1-\delta}{2}
\Tr{\alpha^g (\alpha^g)^T\left( \tilde{w}^{g,m} H \right)_{\xi\xi}} -\mathcal{R}^{\eps,\delta,m,n}\\
&\le  \Tr{\alpha^g(\alpha^g)^T}C(t_n,s_n,z_n)\left(C^*+\frac{ \sqrt{d} \delta \underline{\lambda} }{4 L} \right)+1,
\end{align*}
where we used Lemmas \ref{Lem properties of wm} and \ref{Lem properties of hdelta}.  Furthermore, we remind the reader that the function $ C$ is the one introduced in Lemma \ref{Lem properties of wm}(ii), and the constant $ C^*$ is the one introduced in Lemma \ref{Lem properties of hdelta}. 

\vspace{0.25em}
This last inequality is in contradiction with \eqref{eq interm step 1}, so that we actually have $ | \xi_n| \le \xi^* $. Now since $ (\xi_n)$ is a bounded sequence, we have by classical results of the theory of viscosity solution that, up to extraction, there exists $ \bar{\xi}^{\delta,m}$, with $ \left| \bar{\xi}^{\delta,m}\right|\le 1$, such that $  (t_n,s_n,z_n,\xi_n)\longrightarrow (t_0,s_0,z_0,\bar{\xi}^{\delta,m}) $ when $ n \rightarrow \infty$. Recalling that $ H=1$ on $ \bar{B}_{\xi^*}(0)$, and $ H$ is $ C^2$, we obtain by \eqref{eq interm 5}, \eqref{eq interm 6} and Lemma \ref{Lem properties of wm}(iv) that at the point $ (t_0,s_0,z_0,\bar{\xi}^{\delta,m})$:
\begin{align*}
0 &\geq -\frac12 v^{g}_{zz}\abs{\sigma^T \bar{\xi}^{\delta,m}}^2+\frac{1-\delta}{2}
\Tr{\alpha^g(\alpha^g)^T  \tilde{w}^{g,m}_{\xi\xi}} 
-\mathcal A^g\phi \\
&\geq -\mathcal{A}^g \phi+(1-\delta) a^g-\delta v^g_{zz}\frac{|\sigma^T \bar{\xi}^{\delta,m}|^2}{2}+\frac{1-\delta}{2}\int_{\R^d} \upsilon^m(\zeta) \left|\sigma \zeta \right|^2 d\zeta.
\end{align*}
Now since $ |\bar{\xi}^{\delta,m}|$ is bounded by $ \xi^* $ independent of $ \delta$ and $ m$, we obtain by sending $ \delta$ and $ m$ to 0 that $\mathcal{A}^g \phi (t_0,s_0,z_0)-a^g(t_0,s_0,z_0) \geq 0,$
which contradicts \eqref{eq supersol interm 1}.

\vspace{-2em}

\begin{flushright}$\Box$\end{flushright}

\end{document}